\input amssym.def
\input amssym.tex

\def\item#1{\vskip1.3pt\hang\textindent {\rm #1}}


\newskip\litemindent
\litemindent=0.7cm  
\def\Litem#1#2{\par\noindent\hangindent#1\litemindent
\hbox to #1\litemindent{\hfill\hbox to \litemindent
{\ninerm #2 \hfill}}\ignorespaces}
\def\litem{\Litem1}

\tolerance=300
\pretolerance=200
\hfuzz=1pt
\vfuzz=1pt

\hoffset=0in
\voffset=0.5in

\hsize=5.8 true in 
\vsize=9.2 true in
\parindent=25pt
\mathsurround=1pt
\parskip=1pt plus .25pt minus .25pt
\normallineskiplimit=.99pt

\countdef\revised=100
\mathchardef\emptyset="001F 
\chardef\ss="19
\def\3{\ss}
\def\anf{$\lower1.2ex\hbox{"}$}
\def\frac#1#2{{#1 \over #2}}
\def\>{>\!\!>}
\def\<{<\!\!<}

\def\into{\hookrightarrow}
\def\onto{\to\mskip-14mu\to} 
\def\ssssarr{\hbox to 15pt{\rightarrowfill}}
\def\sssarr{\hbox to 20pt{\rightarrowfill}}
\def\ssarr{\hbox to 30pt{\rightarrowfill}}
\def\sarr{\hbox to 40pt{\rightarrowfill}}
\def\arr{\hbox to 60pt{\rightarrowfill}}
\def\larr{\hbox to 60pt{\leftarrowfill}}
\def\Arr{\hbox to 80pt{\rightarrowfill}}
\def\mapdown#1{\Big\downarrow\rlap{$\vcenter{\hbox{$\scriptstyle#1$}}$}}

\def\sssmapright#1{\smash{\mathop{\sssarr}\limits^{#1}}}
\def\ssmapright#1{\smash{\mathop{\ssarr}\limits^{#1}}}

\def\Alt{\mathop{\rm Alt}\nolimits}
\def\ad{\mathop{\rm ad}\nolimits}

\def\Ad{\mathop{\rm Ad}\nolimits}

\def\Aut{\mathop{\rm Aut}\nolimits}

\def\Alt{\mathop{\rm Alt}\nolimits}

\def\coker{\mathop{\rm coker}\nolimits}

\def\deg{\mathop{\rm deg}\nolimits}
\def\der{\mathop{\rm der}\nolimits}

\def\Diff{\mathop{\rm Diff}\nolimits}

\def\ev{\mathop{\rm ev}\nolimits}

\def\Ext{\mathop{\rm Ext}\nolimits}

\def\Gau{\mathop{\rm Gau}\nolimits}
\def\GL{\mathop{\rm GL}\nolimits}

\def\Hom{\mathop{\rm Hom}\nolimits}%
\def\id{\mathop{\rm id}\nolimits} 
\def\im{\mathop{\rm im}\nolimits}


\def\Int{\mathop{\rm int}\nolimits}

\def\Lin{\mathop{\rm Lin}\nolimits}%

\def\per{\mathop{\rm per}\nolimits}

\def\rk{\mathop{\rm rank}\nolimits}

\def\rk{\mathop{\rm rk}\nolimits}
\def\sgn{\mathop{\rm sgn}\nolimits}
\def\SL{\mathop{\rm SL}\nolimits}
\def\SO{\mathop{\rm SO}\nolimits}
\def\span{\mathop{\rm span}\nolimits}

\def\Sp{\mathop{\rm Sp}\nolimits}



\def\trile{\trianglelefteq}
\def\tor{\mathop{\rm tor}\nolimits}

\def\0{{\bf 0}}
\def\1{{\bf 1}}

\def\a{{\frak a}}

\def\gau{{\frak {gau}}}
\def\b{{\frak b}}

\def\g{{\frak g}}

\def\h{{\frak h}}

\def\n{{\frak n}}

\def\sp{{\frak {sp}}}

\def\so{{\frak {so}}}
\def\sL{{\frak {sl}}}
\def\t{{\frak t}}

\def\z{{\frak z}}

\def\L{\mathop{\bf L{}}\nolimits}

\def\C{{{\Bbb C}{\mskip+1mu}}} 
\def\K{{{\Bbb K}{\mskip+2mu}}} 

\def\R{{\Bbb R}} 
\def\Z{{\Bbb Z}} 
\def\N{{\Bbb N}}

\def\K{{\Bbb K}}

\def\P{{\Bbb P}} 
 
\def\SS{{\Bbb S}} 
\def\T{{\Bbb T}} 

\def\:{\colon}  
\def\.{{\cdot}}
\def\|{\Vert}
\def\bsk{\bigskip}

\def\giantskip{\vskip2\bigskipamount}
\def\gsk{\giantskip}
\def \la {\langle}

\def\msk{\medskip}
\def \ra {\rangle}
\def \res {\!\mid\!\!}

\def\bbr{\bigbreak}
\def\giantbreak{\par \ifdim\lastskip<2\bigskipamount \removelastskip
         \penalty-400 \giantskip\fi}

\def\nin{\noindent}
\def\cen{\centerline}
\def\pagebreak{\vskip 0pt plus 0.0001fil\break}
\def\linebreak{\break}

\def\hat{\widehat}

\def\derat#1{{d \over dt} \hbox{\vrule width0.5pt 
                height 5mm depth 3mm${{}\atop{{}\atop{\scriptstyle t=#1}}}$}}

\def\epsilon{\varepsilon}
\def\eset{\emptyset}

\def\nin{\noindent}
\def\oline{\overline}

\def\pder#1,#2,#3 { {\partial #1 \over \partial #2}(#3)}
\def\pde#1,#2 { {\partial #1 \over \partial #2}}
\def\phi{\varphi}


\def\subeq{\subseteq}

\def\Rarrow{\Rightarrow}

\def\tilde{\widetilde}

\font\ninerm=cmr9
\font\eightrm=cmr8

\font\eightbf=cmbx8


\font\smc=cmcsc10
\font\bfone=cmbx10 scaled\magstep1 
\font\bftwo=cmbx10 scaled\magstep2 

\def\qed{{\unskip\nobreak\hfil\penalty50\hskip .001pt \hbox{}\nobreak\hfil
          \vrule height 1.2ex width 1.1ex depth -.1ex
           \parfillskip=0pt\finalhyphendemerits=0\medbreak}\rm}

\def\qeddis{\eqno{\vrule height 1.2ex width 1.1ex depth -.1ex} $$
                   \medbreak\rm}

\def\Lemma #1. {\bigbreak\vskip-\parskip\noindent{\bf Lemma #1.}\quad\it}

\def\Sublemma #1. {\bigbreak\vskip-\parskip\noindent{\bf Sublemma #1.}\quad\it}

\def\Proposition #1. {\bigbreak\vskip-\parskip\noindent{\bf Proposition #1.}
\quad\it}

\def\Corollary #1. {\bigbreak\vskip-\parskip\nin{\bf Corollary #1.}
\quad\it}

\def\Theorem #1. {\bigbreak\vskip-\parskip\noindent{\bf Theorem #1.}
\quad\it}

\def\Definition #1. {\rm\bigbreak\vskip-\parskip\noindent
{\bf Definition #1.}
\quad}

\def\Remark #1. {\rm\bigbreak\vskip-\parskip\noindent{\bf Remark #1.}\quad}

\def\Example #1. {\rm\bigbreak\vskip-\parskip\noindent{\bf Example #1.}\quad}
\def\Examples #1. {\rm\bigbreak\vskip-\parskip\noindent{\bf Examples #1.}\quad}

\def\Problems #1. {\bigbreak\vskip-\parskip\noindent{\bf Problems #1.}\quad}
\def\Problem #1. {\bigbreak\vskip-\parskip\noindent{\bf Problem #1.}\quad}
\def\Exercise #1. {\bigbreak\vskip-\parskip\noindent{\bf Exercise #1.}\quad}

\def\Conjecture #1. {\bigbreak\vskip-\parskip\noindent{\bf Conjecture #1.}\quad}

\def\Proof#1.{\rm\par\ifdim\lastskip<\bigskipamount\removelastskip\fi\smallskip
            \noindent {\bf Proof.}\quad}

\def\Axiom #1. {\bigbreak\vskip-\parskip\noindent{\bf Axiom #1.}\quad\it}

\def\Satz #1. {\bigbreak\vskip-\parskip\noindent{\bf Satz #1.}\quad\it}

\def\Korollar #1. {\bbr\vskip-\parskip\nin{\bf Korollar #1.} \quad\it}

\def\Folgerung #1. {\bbr\vskip-\parskip\nin{\bf Folgerung #1.} \quad\it}

\def\Folgerungen #1. {\bbr\vskip-\parskip\nin{\bf Folgerungen #1.} \quad\it}

\def\Bemerkung #1. {\rm\bigbreak\vskip-\parskip\noindent{\bf Bemerkung #1.}
\quad}

\def\Beispiel #1. {\rm\bigbreak\vskip-\parskip\noindent{\bf Beispiel #1.}\quad}
\def\Beispiele #1. {\rm\bigbreak\vskip-\parskip\noindent{\bf Beispiele #1.}\quad}
\def\Aufgabe #1. {\rm\bigbreak\vskip-\parskip\noindent{\bf Aufgabe #1.}\quad}
\def\Aufgaben #1. {\rm\bigbreak\vskip-\parskip\noindent{\bf Aufgabe #1.}\quad}

\def\Beweis#1. {\rm\par\ifdim\lastskip<\bigskipamount\removelastskip\fi
           \smallskip\noindent {\bf Beweis.}\quad}

\nopagenumbers

\def\date{\ifcase\month\or January\or February \or March\or April\or May
\or June\or July\or August\or September\or October\or November
\or December\fi\space\number\day, \number\year}

\def\title{Title ??}
\def\author{Author ??}

\def\thanks#1{\footnote*{\eightrm#1}}

\def\rightheadline{\hfil{\eightrm\title}\hfil\tenbf\folio}
\def\leftheadline{\tenbf\folio\hfil{\eightrm\author}\hfil}
\headline={\vbox{\line{\ifodd\pageno\rightheadline\else\leftheadline\fi}}}

\def\firstheadline{}
\def\firstfootline{\cen{\rm\folio}}

\def\seite #1 {\pageno #1
               \headline={\ifnum\pageno=#1 \firstheadline
               \else\ifodd\pageno\rightheadline\else\leftheadline\fi\fi}
               \footline={\ifnum\pageno=#1 \firstfootline\else{}\fi}}

\newdimen\dimenone
 \def\checkleftspace#1#2#3#4{
 \dimenone=\pagetotal
 \advance\dimenone by -\pageshrink   
 \ifdim\dimenone>\pagegoal          
   \else\dimenone=\pagetotal
        \advance\dimenone by \pagestretch
        \ifdim\dimenone<\pagegoal
          \dimenone=\pagetotal
          \advance\dimenone by#1         
          \setbox0=\vbox{#2\parskip=0pt                
                     \hyphenpenalty=10000
                     \rightskip=0pt plus 5em
                     \noindent#3 \vskip#4}    
        \advance\dimenone by\ht0
        \advance\dimenone by 3\baselineskip   
        \ifdim\dimenone>\pagegoal\vfill\eject\fi
          \else\eject\fi\fi}


\def\subheadline #1{\nin\bigbreak\vskip-\lastskip
      \checkleftspace{0.9cm}{\bf}{#1}{\medskipamount}
          \indent\vskip0.7cm\centerline{\bf #1}\medskip}
\def\subsection{\subheadline} 

\def\lsubheadline #1 #2{\bigbreak\vskip-\lastskip
      \checkleftspace{0.9cm}{\bf}{#1}{\bigskipamount}
         \vbox{\vskip0.7cm}\cen{\bf #1}\msk \cen{\bf #2}\bsk}

\def\sectionheadline #1{\bigbreak\vskip-\lastskip
      \checkleftspace{1.1cm}{\bf}{#1}{\bigskipamount}
         \vbox{\vskip1.1cm}\cen{\bfone #1}\bsk}
\def\section{\sectionheadline} 

\def\lsectionheadline #1 #2{\bigbreak\vskip-\lastskip
      \checkleftspace{1.1cm}{\bf}{#1}{\bigskipamount}
         \vbox{\vskip1.1cm}\cen{\bfone #1}\msk \cen{\bfone #2}\bsk}

\def\lchapterheadline #1 #2{\bigbreak\vskip-\lastskip\indent\vskip3cm
                       \cen{\bftwo #1} \msk \cen{\bftwo #2} \gsk}
\def\llsectionheadline #1 #2 #3{\bigbreak\vskip-\lastskip\indent\vskip1.8cm
\cen{\bfone #1} \msk \cen{\bfone #2} \msk \cen{\bfone #3} \nobreak\bsk\nobreak}


\newtoks\literat
\def\[#1 #2\par{\literat={#2\unskip.}%
\hbox{\vtop{\hsize=.15\hsize\nin [#1]\hfill}
\vtop{\hsize=.82\hsize\nin\the\literat}}\par
\vskip.3\baselineskip}

\def\references{
\sectionheadline{\bf References}
\frenchspacing

\entries\par}

\mathchardef\emptyset="001F 
\def\address{Author: \tt$\backslash$def$\backslash$address$\{$??$\}$}

\def\abstract #1{{\narrower\baselineskip=10pt{\noindent
\eightbf Abstract.\quad \eightrm #1 }
\bigskip}}

\def\firstpage{\nin
{\obeylines \parindent 0pt }
\vskip2cm
\centerline{\bfone\title}
\gsk
\centerline{\bf\author}
\vskip1.5cm \rm}

\def\lastpage{\par\vbox{\vskip1cm\nin
\line{
\vtop{\hsize=.5\hsize{\parindent=0pt\baselineskip=10pt\nin\address}}
\hfill} }}

\def\Box #1 { \msk\par\nin 
\centerline{
\vbox{\offinterlineskip
\hrule
\hbox{\vrule\strut\hskip1ex\hfil{\smc#1}\hfill\hskip1ex}
\hrule}\vrule}\msk }

\def\adots{\mathinner{\mkern1mu\raise1pt\vbox{\kern7pt\hbox{.}}
                        \mkern2mu\raise4pt\hbox{.}
                        \mkern2mu\raise7pt\hbox{.}\mkern1mu}}


\pageno=1

\def\title{Abelian extensions of infinite-dimensional Lie groups} 
\def\author{Karl-Hermann Neeb}
\def\date{February 18, 2004}
\def\rightheadline{\tenbf\folio\hfil{\tt abelext.tex}\hfil\eightrm\date}
\def\leftheadline{\tenbf\folio\hfil{\rm\title}\hfil\eightrm\date}

\def\Alt{\mathop{\rm Alt}\nolimits}
\def\pr{\mathop{\rm pr}\nolimits}

\def\address
{Karl-Hermann Neeb

Technische Universit\"at Darmstadt 

Schlossgartenstrasse 7

D-64289 Darmstadt 

Deutschland

neeb@mathematik.tu-darmstadt.de}

\firstpage 

\abstract{In the present paper we study abelian extensions of connected 
Lie groups ${\scriptstyle G}$ modeled on locally convex spaces by 
smooth ${\scriptstyle G}$-modules ${\scriptstyle A}$. We 
parametrize the extension classes by a suitable cohomology group ${\scriptstyle H^2_s(G,A)}$ 
defined by locally smooth cochains and construct an exact sequence that describes 
the difference between ${\scriptstyle H^2_s(G,A)}$ and the corresponding continuous Lie algebra 
cohomology space ${\scriptstyle H^2_c(\g,\a)}$. 
The obstructions for the integrability of a Lie algebra 
extensions to a Lie group extension are described in terms of period and 
flux homomorphisms. We also characterize the extensions with global smooth 
sections resp.\ those given by global smooth cocycles. Finally we apply 
the general theory to extensions of several types of diffeomorphism groups. }

\section{Introduction} 

The main point of the present paper is a detailed analysis of abelian 
extensions of Lie groups $G$ which might be infinite-dimensional, 
a main point being to derive criteria for abelian extensions of Lie algebras to integrate 
to extensions of corresponding connected groups. This is of particular interest for infinite-dimensional Lie algebras because not every infinite-dimensional Lie algebra 
can be `integrated' to a global Lie group. 

The concept of a (not necessarily finite-dimensional) Lie group used here 
is that a {\it Lie group} $G$ is a manifold modeled on a locally convex space  
endowed with a group structure for which the group operations are smooth 
(cf.\ [Mi83]; see also [Gl01] for non-complete model spaces). 
An {\it abelian extension} is an exact sequence of Lie groups 
$A \into \hat G \onto G$ which defines a locally trivial smooth 
principal bundle with the abelian structure group $A$ over the Lie group $G$. 
Then $A$ inherits the structure of a {\it smooth $G$-module} in the sense that 
the conjugation action of $\hat G$ on $A$ factors through 
a smooth map $G \times A \to A$. The extension is called {\it central} if this  
action is trivial. 

The present paper is a sequel to [Ne02] which deals with the 
case of central extensions. Fortunately it was possible to use some of the 
constructions from [Ne02] quite directly in the present paper, but a substantial 
part of the machinery used for central extensions had to be generalized 
and adapted to deal with abelian extensions. In [Ne04] it is shown that the results 
on abelian extensions can in turn be used to classify general extensions. 

A typical class of examples that illustrate the difference between abelian 
and central extensions of Lie groups arises from 
abelian principal bundles. If $q \: P \to M$ is a smooth principal bundle 
with the abelian structure group $Z$ over the compact connected manifold $M$, 
then the group $\Diff(P)^Z$ of all diffeomorphisms of $P$ commuting with 
$Z$ (the automorphism group of the bundle) is an extension of 
an open subgroup of $\Diff(M)$ by the gauge group 
$\Gau(P) \cong C^\infty(M,Z)$ of the bundle. Here the conjugation action of 
$\Diff(M)$ on $\Gau(P)$ is given by composing functions with diffeomorphisms. 
Central extensions corresponding to the bundle $q \: P \to M$ are obtained by 
choosing a principal connection $1$-form 
$\theta \in \Omega^1(P,\z)$. Let $\omega \in \Omega^2(M,\z)$ denote the corresponding 
curvature form. Then the subgroup 
$\Diff(P)^Z_\theta$ of those elements of $\Diff(P)^Z$ preserving $\theta$ 
is a central extension of an open subgroup of $\Diff(M)_\omega$, which is substantially 
smaller that $\Diff(M)$. This example shows that the passage from central 
extensions to abelian extensions is similar to the passage from 
symplectomorphism groups to diffeomorphism groups. 

As the examples of principal bundles over compact manifolds show, abelian 
extensions of Lie groups occur naturally in geometric contexts and in particular 
in symplectic geometry, where the prequantization problem is to find for a 
symplectic manifold $(M,\omega)$ a $\T$-principal bundle with curvature $\omega$, 
which leads to an abelian extension of $\Diff(M)_0$ by the group $C^\infty(M,\T)$. 
Conversely, every abelian extension $q \: \hat G \to G$ 
of a Lie group $G$ by an abelian Lie group $A$ is in particular an $A$-principal 
bundle over $G$, so that there is a close interplay between abelian extensions 
of infinite-dimensional groups and abelian principal bundles over (finite-dimensional) 
manifolds. 

In the representation theory of infinite-dimensional Lie groups abelian extensions 
occur naturally if a connected Lie group $G$ acts on a smooth manifold $M$ which is 
endowed with a $Z$-principal bundle $q \: P \to M$, each element of $G$ lifts 
to an automorphism of the bundle, but there is no principal connection 
$1$-form preserved by the lifts of the elements of $G$ to diffeomorphisms of $P$. 
We refer to [Mi89] for a detailed discussion of the case where $M$ is a restricted 
Gra\3mannian of a polarized Hilbert space and the groups are 
restricted operator groups of Schatten class $p > 2$, resp., mapping groups 
$C^\infty(M,K)$, where $K$ is finite-dimensional and $M$ is a compact manifold 
of dimension $\geq 2$ (see also [PS86] for a discussion of related points). 
Since representations of abelian extensions of vector field Lie algebras occur naturally 
in mathematics physics (cf.\ [La99] and also [AI95] for more general applications 
of Lie group cohomology in physics), 
the question arises whether this picture 
has a global analog in terms of abelian extensions of the corresponding diffeomorphism groups.
Some first results in this direction have been obtain by Y.~Billig in [Bi03], where he 
introduces natural analogs of the Virasoro group which are abelian extensions of 
$\Diff(M)$. 

Another motivation for a general study of abelian extensions comes from 
the fact that for the group $G := \Diff(M)_0$, where $M$ is a compact orientable  
manifold, one has natural modules given by tensor densities and spaces of 
tensors on $M$.  The corresponding abelian extensions can be used to interprete 
certain partial differential equations as geodesic equations on a Lie group, 
which leads to important information on the behavior of their solutions ([Vi02], 
[AK98]). 
An important special case discussed in some detail in Section X is the 
group of diffeomorphisms of the circle and its modules of $\lambda$-densities for real $\lambda$. 
For the identity component 
$D(M,\mu)$ of the group $\Diff(M,\mu)$ of volume preserving diffeomorphisms 
(for a given volume form $\mu$) one obtains 
a Lie algebra cocycle from each closed $2$-form $\omega$ on $M$ 
(Lichnerowicz cocycle) which is obtained by composing the integration map 
with $\omega$, interpreted as a $2$-cocycle for ${\cal V}(M)$ with values 
in the smooth module $C^\infty(M,\R)$. The existence of corresponding 
central extensions is addressed for special cases in Section XI, where we use 
relevant information on the associated abelian extensions of $\Diff(M)_0$ obtained in 
Section IX. 
For more references dealing specifically with central extensions we refer to 
[Ne02]. See in particular [CVLL98] which is a nice survey of central $\T$-extensions of 
Lie groups and their role in quantum physics. That paper also contains 
a description of the universal central extension for finite-dimensional groups. 
For infinite-dimensional groups universal central extensions are constructed 
in [Ne03b] and for root graded Lie algebras in [Ne03a]. 

\msk 

As one would expect from general homological algebra, the natural context  
to deal with abelian extensions of Lie groups is provided by a suitable 
Lie group cohomology with values in smooth modules:  
If $G$ is a Lie group, then we call an abelian Lie group $A$ a {\it
  smooth $G$-module} if it is a $G$-module and the action map 
$G \times A \to A$ is smooth. In Appendix~B we describe a natural
adaptation of the group cohomology complex to the Lie group setting. 
Here the space of $n$-cochains $C^n_s(G,A)$ consists of maps 
$G^n \to A$ which are smooth in an identity neighborhood and 
vanish on all tuples of the form $(g_1,\ldots, \1, \ldots, g_n)$. 
We thus obtain a cochain complex $(C^n_s(G,A), d_G)$ with the
cohomology groups $H^n_s(G,A)$. 
If $G$ and $A$ are discrete, these groups coincide with the
standard cohomology groups of $G$ with values in $A$. We 
refer to [Mo64] and [Mo76] for an appropriate definition of topological 
group cohomology which fits well for locally compact groups.   
Since the cohomology groups 
$H^n_s(G,A)$ can be considered as rather complicated objects, it is
desirable to relate them to the corresponding Lie algebra cohomology
groups $H^n_s(\g,\a)$. Passing to the derived representation, 
the Lie algebra $\a$ of $A$ is a module of the Lie algebra $\g$ of
$\a$ which is {\it topological} in the sense that the module structure
is a continuous bilinear map $\g \times \a \to \a$. Then 
the continuous alternating maps $\g^n \to \a$ form the (continuous)
Lie algebra cochain complex $(C^n_c(\g,\a),d_\g)$, and its cohomology
spaces are denoted $H^n_c(\g,\a)$. 

In Appendix B we show that for $n \geq 2$ there is a natural {\it derivation map} 
$$ D_n \: H^n_s(G,A) \to H^n_c(\g,\a) $$
from locally smooth Lie group cohomology to continuous Lie algebra cohomology. This map 
is based on the isomorphism 
$$ H^n_c(\g,\a) \cong H^n_{\rm dR,eq}(G,\a) $$
between Lie algebra cohomology and the de Rham cohomology of the complex of 
equivariant $\a$-valued differential forms on $G$ (cf.\ [CE48] for finite-dimensional groups). For $n = 1$ we only have a map $D_1 \: Z^1_s(G,A) \to Z^1_c(\g,\a)$, 
and if $A \cong \a/\Gamma_A$ holds for a discrete subgroup $\Gamma_A$ of
$\a$, then this map factors to a map on the level of cohomology. Since 
the Lie algebra cohomology spaces $H^n_c(\g,\a)$ are
much better accessible by algebraic means than those of $G$, it is important to understand 
the amount of information lost by the map $D_n$. More concretely, one
is interested in kernel and cokernel of $D_n$. A determination of
the cokernel can be considered as describing integrability conditions
on cohomology classes $[\omega] \in H^n_c(\g,\a)$ which have to be
satisfied to ensure the existence of $f \in Z^n_s(G,A)$ with $D_n f =
\omega$.

In the present paper we completely solve this problem for 
the important case $n =2$, a connected Lie group $G$ and connected 
smooth modules $A$ of the form 
$\a/\Gamma_A$, where $\Gamma_A$ is a discrete subgroup of $\a$. 
We also describe the solution for $n = 1$
which is much simpler, but already reflects the spirit of the
problem. We plan to return in a subsequent paper to this problem for
non-connected groups $G$, which, in view of the present results, 
means to obtain accessible criteria for
the extendibility of a $2$-cocycle on the identity component $G_0$ of
$G$ to the whole group $G$. 

The special importance of the group $H^2_s(G,A)$ stems from the fact
that for connected groups $G$ 
it classifies all Lie group extensions $q \: \hat G \to G$ of $G$
by $A \cong \ker q$, where the action of $G$ on $A$ induced by the
conjugation action of $\hat G$ on the abelian normal subgroup $A$ coincides
with the original $G$-module structure. This was our original
motivation to study the cohomology groups $H^2_s(G,A)$. If $G$ is not
connected, then we have to consider an appropriate subgroup 
$H^2_{ss}(G,A) \subeq H^2_s(G,A)$ which then classifies the 
extensions of $G$ by $A$. 

The second cohomology groups do not only classify abelian extensions
of $G$, they also play an equally important role 
in the classification of general extensions:  
Let $N$ be a Lie group and $Z(N)$ its center. Suppose further that $Z(N)$ is a
smooth $G$-module such that every smooth map $M \to N$ with values in
$Z(N)$ defines a smooth map $M \to Z(N)$. 
Then the group $H^2_s(G,Z(N))$ parameterizes the equivalence classes of
extensions of $G$ by $N$ corresponding to a given smooth outer action
of $G$ on $N$ (see [Ne04] for the details and the definition of a smooth outer action). 
If $N = Z(N)$ is abelian,
then a smooth outer action of $G$ on $N$ is the same as a smooth
module structure. 

Taking the derivation maps $D_n$ into account, we obtain for connected
groups $G$ and $A \cong \a/\Gamma_A$ the following
commutative diagram with an exact second row (see Proposition~D.8 and the subsequent 
discussion) and exact columns 
(Proposition~III.4 and Theorem VII.2): 

$$\!\!\!\!\!\!\!\!\! \matrix{ 
\0&& \0 &&  && \Hom(\pi_1(G),A^G) & & \0\cr 
\mapdown{}&&\mapdown{} &&  && \mapdown{\delta}& & \mapdown{}\cr 
H^1_s(G, A) &\!\!\!\!\!\!\!\!\!\!\!\!\sssmapright{I} & H^1_s(\tilde
G,A) & \!\!\!\!\!\!\sssmapright{R}\!\!\!\!
& \Hom(\pi_1(G), A^G) &\!\!\!\!\!\! \sssmapright{\delta} \!\!\!& H^2_s(G,A) & 
\!\!\!\!\!\!\!\sssmapright{I} \!\!\!\!\!\!& H^2_s(\tilde G,A) \cr 
\mapdown{D_1} & & \mapdown{D_1} &  &\mapdown{\id}& & \mapdown{D_2} & &\mapdown{D_2}\cr 
H^1_c(\g,\a) &\!\!\!\! \!\!\!\!\!\!\!\!\sssmapright{\id}
& H^1_c(\g,\a) &\!\!\!\!\!\!\sssmapright{P_1}\!\!\!\!\!\!& \Hom(\pi_1(G),A^G)&& 
H^2_c(\g,\a)& \!\!\!\!\!\!\!\sssmapright{\id} \!\!\!& H^2_c(\g,\a)\cr  
\mapdown{P_1} &&\mapdown{} &&  && \mapdown{P_2}& &\mapdown{P_2} \cr 
\,\,\,\,\,\,\,\,\,\,\Hom(\pi_1(G), A^G)\!\!\!\!\!\!\!\! && \0 &&  &&{\Hom(\pi_2(G),A^G) \oplus \atop 
\Hom(\pi_1(G), H^1_c(\g,\a))} & &\!\!\!\!\Hom(\pi_2(G),A^G)\,\,\, \cr} $$

Here $I$ denotes natural inflation maps, 
$\delta$ assigns to $\gamma \: \pi_1(G) \to A^G$ the quotient of 
the semi-direct product $A \rtimes \tilde G$ modulo the graph of $\gamma$ 
and. For a Lie algebra cocycle $\alpha \in Z^1_c(\g,\a)$ 
the homomorphism $P_1([\alpha]) \: \pi_1(G) \to A^G$ 
is obtained by integrating the corresponding equivariant $1$-form 
$\alpha^{\rm eq} \in \Omega^1(G,\a)$ over loops and then interpreting the result 
an element of the quotient group $A = \a/\Gamma_A$. 
For a $2$-cocycle $\omega \in Z^2_c(\g,\a)$ and 
the corresponding equivariant $2$-form $\omega^{\rm eq} \in \Omega^2(G,\a)$ 
on $G$, the first component of $P_2([\omega])$ 
is the period homomorphism $\pi_2(G) \to A^G$ obtained by integrating 
$\omega^{\rm eq}$ over smooth maps $\SS^2 \to G$ and then interpreting the result 
modulo $\Gamma_A$ as an element of $A$. The second component of 
$P_2([\omega])$ is the flux homomorphism 
$F_\omega \: \pi_1(G) \to H^1_c(\g,\a)$ which can be viewed, in a certain sense, 
as $P_1([f_\omega])$ for the Lie algebra flux cocycle 
$$f_\omega \: \g \to C^1_c(\g,\a)/d_\g \a, \quad  x \mapsto [i_x \omega] $$
(in Section VI we give a direct definition which does not require to 
topologize the space $C^1_c(\g,\a)$ and its quotient space module $B^1_c(\g,\a) = d_\g \a$). 

If $G$ is simply connected, things become much simpler and the 
criterion for the integrability of a Lie algebra cocycle 
$\omega$ to a group cocycle is that all periods of $\omega^{\rm eq}$ are contained 
in $\Gamma_A \subeq \a$. 
Similar conditions arise in the theory of abelian principal bundles 
on  smoothly paracompact presymplectic manifolds $(M,\Omega)$ ($\Omega$ is a 
closed $2$-form on $M$). 
Here the integrality of the cohomology class $[\Omega]$ 
is equivalent to the existence of a 
so-called pre-quantum bundle, i.e., a $\T$-principal bundle $\T \into
\hat M \onto M$ whose curvature $2$-form is $\Omega$ (cf.\ [Bry90]). 

For finite-dimensional Lie groups the integrability 
criteria also simplify significantly because $\pi_2(G)$ vanishes 
([Ca52]). This in turn has been used by \'E.~Cartan 
to construct central extensions and thus 
to prove Lie's Third Theorem that each finite-dimensional Lie algebra belongs to a 
global Lie group.  We generalize Cartan's construction in Section~VIII 
to characterize abelian extensions with global smooth sections. 

We emphasize that our results hold for Lie groups which are not necessarily 
smoothly paracompact, so that one cannot use smooth partitions of unity 
to construct bundles for prescribed curvature forms and de Rham's Theorem is 
not available (cf.\ [KM97,  Th.~16.10]). This point is important because 
many interesting Banach--Lie groups are not smoothly paracompact which comes from the fact 
that their model spaces do not permit smooth bump functions (cf.\ [KM97]). 

\msk 

The contents of the present paper is as follows. 
In Section~I we briefly discuss the relation between abelian extensions of 
topological Lie algebras and the continuous cohomology space 
$H^2_c(\g,\a)$ (see [CE48] for the case of abstract Lie algebras). 
The parameterization of the class of all Lie group extensions of 
a connected Lie group $G$ by $A$ via the cohomology group $H^2_s(G,A)$ is obtained 
in Section~II. In Section~III we briefly discuss the relation between locally smooth 
$1$-cocycles on Lie groups and the corresponding continuous Lie algebra cocycles. 
This is instructive for the understanding of the flux cocycle occurring below as an 
obstruction to the existence of global group extensions. 
In Section~IV we briefly discuss the period homomorphism 
$\per_\omega \: \pi_2(G) \to \a^G$ associated to a Lie algebra cocycle 
$\omega \in Z^2_c(\g,\a)$. To integrate Lie algebra cocycles on simply 
connected groups in Section~V we use a slight adaptation of the method used in [Ne02] 
for central extensions. In Section~VI we eventually turn to the refinements needed 
for non-simply connected groups which leads to the flux cocycle. This part is 
considerably  
more complicated than for central extensions where the flux cocycle simplifies to a 
homomorphism with values in a space of homomorphisms of Lie algebras 
and not only in a cohomology space. 
In Section VII all pieces are put together to obtain the exactness of rows 
and columns in the big diagram above. Abelian extensions with smooth global 
sections are characterized in Section~VIII. 

The remaining 
Sections IX-XI contain examples and some discussion of special cases. 
In Section IX we turn to the special situation of diffeomorphism groups on 
compact manifolds and the special class of $2$-cocycles on 
the Lie algebra ${\cal V}(M)$ given by closed $2$-forms on $M$. 
In this case we explain how interesting information on period 
map and flux cocycle can be calculated in geometrical terms. 
In Section X we consider the situation where 
$G$ is the diffeomorphism group of the circle and $\a$ is the module of 
$\lambda$-densities for some $\lambda \in \R$. The corresponding cocycles 
for $\Diff(\SS^1)$ have been discussed by Ovsienko and Roger in 
[OR98]. In particular we describe how their results can be extended 
to Lie algebra cocycles not integrable on $\Diff(\SS^1)_0$ 
which integrate to group cocycles of the 
universal covering group $\tilde\Diff(\SS^1)$. As  a byproduct of this construction, 
we obtain a non-trivial abelian 
extension of the group $\tilde\SL_2(\R)$ by an infinite-dimensional 
Fr\'echet space. Section XI contains 
some information on the integration of Lichnerowicz cocycles to 
central group extensions. In particular we show that for a compact connected 
Lie group $G$ each Lichnerowicz cocycle on $D(G,\mu)$ can be integrated 
to a corresponding group cocycle on the covering group 
$\tilde D(G,\mu)$ acting as a group of diffeomorphisms on the universal 
covering group $\tilde G$. 

We conclude this paper with several appendices dealing 
with the relation between differential forms and Alexander--Spanier 
cohomology (Appendix A), cohomology of Lie groups and Lie algebras 
(Appendix B), constructing global Lie groups from local data (Appendix C), 
the exact Inflation-Restriction Sequence for Lie group cohomology 
(Appendix D), the long exact sequence in Lie group cohomology induced 
from a topologically split 
exact sequence of smooth modules (Appendix E), and 
multiplication of Lie group and Lie algebra cocycles (Appendix F). 

We are grateful to S.~Haller for providing a crucial argument concerning the 
flux homomorphism for the group of volume preserving diffeomorphisms 
(cf.\ Section XI). We also thank C.~Vizman for many inspiring discussions 
on the subject, G.~Segal for suggesting a different type of obstructions to 
the ingrability of abelian extensions in [Se02], and to 
A.~Dzhumadildaev for asking for global central extensions 
of groups of volume preserving diffeomorphisms which correspond to the 
cocycles he studied on the Lie algebra level ([Dz92]). This led us to the 
results in Section~XI. 


\section{0. Preliminaries and notation} 

In this paper $\K \in \{\R,\C\}$ denotes the field of real or complex
numbers. Let $X$ and $Y$ be topological $\K$-vector spaces, $U
\subeq X$ open and $f \: U \to Y$ a map. Then the {\it derivative
  of $f$ at $x$ in the direction of $h$} is defined as 
$$ df(x)(h) := \lim_{t \to 0} {1 \over t} \big( f(x + t h) - f(x)\big)
$$
whenever the limit exists. The function $f$ is called {\it differentiable at
  $x$} if $df(x)(h)$ exists for all $h \in X$. It is called {\it
  continuously differentiable or $C^1$} if it is continuous and differentiable at all
points of $U$ and 
$$ df \: U \times X \to Y, \quad (x,h) \mapsto df(x)(h) $$
is a continuous map. It is called a $C^n$-map if $f$ is $C^1$ and $df$ is a
$C^{n-1}$-map, and $C^\infty$ ({\it smooth}) if it is $C^n$ for all $n \in \N$. 
This is the notion of differentiability used in [Mil83], and
[Gl01], where the latter reference deals with the modifications
necessary for incomplete spaces. 

Since we have a chain rule for $C^1$-maps between locally convex 
spaces ([Gl01]), we can define smooth manifolds $M$ as in
the finite-dimensional case. 
A Lie group $G$ is a smooth manifold modeled on a locally convex space 
$\g$ for which the group multiplication and the 
inversion are smooth maps. We write $\1 \in G$ for the identity
element, $\lambda_g(x) = gx$ for left multiplication, 
$\rho_g(x) = xg$ for right multiplication, and 
$c_g(x) := gxg^{-1}$ for conjugation. Then each $x \in T_\1(G)$ corresponds to
a unique left invariant vector field $x_l$ with 
$x_l(g) := d\lambda_g(\1).x, g \in G.$
The space of left invariant vector fields is closed under the Lie
bracket of vector fields, hence inherits a Lie algebra structure. In
this sense we obtain on $\g := T_\1(G)$ a continuous Lie bracket which
is uniquely determined by $[x,y]_l = [x_l, y_l]$. 
We call a Lie algebra $\g$ which is a topological vector space 
such that the Lie bracket is continuous a {\it topological Lie
  algebra $\g$}. In this sense the Lie algebra $\g = \L(G)$ of a Lie
group $G$ is a locally convex topological Lie algebra. If $G$ is a connected 
Lie group, then we write $q_G \: \tilde G \to G$ for its universal covering Lie 
group and identify $\pi_1(G)$ with the kernel of $q_G$.

Throughout this paper we write abelian groups $A$ additively with $\0$
as identity element. 
If $G$ is a Lie group, then a {\it smooth $G$-module} is an abelian
Lie group $A$, endowed with a smooth $G$-action $\rho_A \: G \times A \to A$ by 
group automorphisms. We sometimes write $(A,\rho_A)$ to include the
notation $\rho_A$ for the action map. 
If $\a$ is the Lie algebra of $A$, then the
smooth action induces a smooth action on $\a$, so that $\a$ also 
is a smooth $G$-module, hence also a module of the Lie algebra $\g$ of
$G$.  In the following we shall mostly assume that the identity component 
$A_0$ of $A$ is of the form $A_0 \cong \a/\Gamma_A$, where 
$\Gamma_A \subeq \a$ is a discrete subgroup. Then the quotient map 
$q_A \: \a \to A_0$ is the universal covering map of $A_0$ and $\pi_1(A) 
\cong \Gamma_A$. 

A linear subspace $W$ of a topological vector space $V$ is called 
{\it (topologically) split} if it is closed and there is a continuous linear map 
$\sigma \: V/W \to V$ for which the map 
$$ W \times V/W \to V, \quad (w,x) \mapsto w + \sigma(x) $$
is an isomorphism of topological vector spaces. Note that the closedness of 
$W$ guarantees that the quotient topology turns $V/W$ into a Hausdorff space which 
is a topological vector space with respect to 
the induced vector space structure. 
A continuous linear map $f \: V \to W$ between topological vector spaces  
is said to be {\it (topologically) split} if the subspaces 
$\ker(f) \subeq V$ and $\im(f) \subeq W$ are topologically split.

\sectionheadline{I. Abelian extensions of topological Lie algebras} 

For the definition of the cohomology of a topological Lie algebra $\g$
with values in a topological $\g$-module $\a$ we refer to Appendix~B. 

\Definition I.1. Let $\g$ and $\n$ be topological Lie algebras. A topologically split 
short exact sequence 
$$ \n \into \hat\g \onto \g $$
is called a {\it (topologically split) extension of $\g$ by $\n$}. 
We identify $\n$ with its image 
in $\hat\g$, and write $\hat\g$ as a direct sum $\hat\g = \n \oplus
\g$ of topological vector spaces. Then $\n$ is a topologically split ideal
and the quotient map $q \: \hat\g \to \g$ corresponds to $(n,x)
\mapsto x$. If $\n$ is abelian, then the extension is called {\it abelian}. 

Two extensions 
$\n \into \hat\g_1 \onto \g$ and 
$\n \into \hat\g_2 \onto \g$ are called {\it equivalent} if there exists a 
morphism $\phi \: \hat\g_1 \to \hat\g_2$ of topological Lie algebras 
such that the diagram 
$$ \matrix{ 
\n & \into & \hat \g_1 & \onto & \g \cr 
\mapdown{\id_\n} &  & \mapdown{\phi} & & \mapdown{\id_\g} \cr 
\n & \into & \hat \g_2 & \onto & \g \cr } $$
commutes. It is easy to see that this implies that 
$\phi$ is an isomorphism of topological Lie algebras, hence
defines an equivalence relation. We write $\Ext(\g,\n)$
for the set of equivalence classes of extensions of $\g$ by $\n$. 

We call an extension $q \: \hat\g \to \g$ with $\ker q = \n$ 
{\it trivial}, or say that the extension 
{\it splits}, if there exists a continuous Lie algebra homomorphism 
$\sigma \: \g \to \hat\g$ with $q \circ \sigma = \id_\g$. 
In this case the map 
$$ \n \rtimes_S \g \to \hat\g, \quad 
(n,x) \mapsto n + \sigma(x) $$
is an isomorphism, where the semi-direct sum is defined by the 
homomorphism 
$$ S \: \g \to \der(\n), \quad S(x)(n) := [\sigma(x),n]. 
\qeddis 

\Definition I.2. Let $\a$ be a topological $\g$-module. 
To each continuous $2$-cocycle $\omega \in Z^2_c(\g,\a)$ 
we associate a topological Lie algebra $\a \oplus_\omega \g$ as the topological 
product vector space $\a \times \g$ endowed with the Lie bracket 
$$ [(a,x), (a', x')] := (x.a' - x'.a + \omega(x,x'),[x,x']). $$
The quotient map $q \: \a \oplus_\omega \g \to \g, (a,x) \mapsto x$ is a continuous homomorphism 
of Lie algebras with kernel $\a$, hence defines an $\a$-extension of $\g$. 
The map $\sigma \: \g \to \a \oplus_\omega \g, x \mapsto (0,x)$ is a continuous 
linear section of $q$. 
\qed

\Proposition I.3. Let $(\a,\rho_\a)$ be a topological $\g$-module and
write $\Ext_{\rho_\a}(\g,\a)$ for the set of all equivalence classes of
$\a$-extensions $\hat\g$ of $\g$ for which the adjoint action of
$\hat\g$ on $\a$ induces the given $\g$-module structure on $\a$. 
Then the map 
$$ Z^2_c(\g,\a) \to \Ext_{\rho_\a}(\g,\a), \quad 
\omega \mapsto [\a \oplus_\omega \g] $$
factors through a bijection 
$$ H^2_c(\g,\a) \to \Ext_{\rho_\a}(\g,\a), \quad 
[\omega] \mapsto [\a \oplus_\omega \g]. $$

\Proof. Suppose that $q \: \hat\g \to \g$ is an $\a$-extension of
$\g$ for which the induced $\g$-module structure on $\a$ coincides
with $\rho_\a$. Let 
$\sigma \: \g \to \hat\g$ 
be a continuous linear section, so that $q \circ \sigma = \id_\g$. Then 
$$ \omega(x,y) := [\sigma(x), \sigma(y)] - \sigma([x,y]) $$
has values in the subspace $\a = \ker q$ of $\hat\g$ and the map 
$\a \times \g \to \hat\g, (a,x) \mapsto a + \sigma(x)$
defines an isomorphism of topological Lie algebras $\a \oplus_\omega \g \to \hat\g$. 

It is easy to verify that $\a \oplus_\omega \g \sim \a\oplus_\eta \g$ if and only if 
$\omega - \eta \in B^2_c(\g,\a)$. Therefore the quotient space 
$H^2_c(\g,\a)$ classifies the equivalence classes of $\a$-extensions of $\g$ 
by the assignment $[\omega] \mapsto [\a \oplus_\omega \g]$. 
\qed

\sectionheadline{II. Abelian extensions of Lie groups} 

Let $A$ be a smooth $G$-module. 
In this section we explain how to assign to a cocycle $f \in
Z^2_s(G,A)$ (satisfying some additional smoothness condition if $G$ is
not connected) a Lie group $A \times_f G$ which is an extension of $A$
by $G$ for which the induced action of $G$ on $A$ coincides with the original
one. We shall see that this assignment leads to a bijection between a
certain subgroup $H^2_{ss}(G,A)$ of $H^2_s(G,A)$ with the set of
equivalence classes of extensions of $G$ by the smooth $G$-module $A$.
If $G$ is connected, then $H^2_{ss}(G,A) = H^2_s(G,A)$. 
We also show that the assingment $f \mapsto A \times_f G$ is
compatible with the derivation map $D \: Z^2_s(G,A) \to Z^2_c(\g,\a)$
in the sense that $\a \oplus_{Df} \g$ is the Lie algebra of $A
\times_f G$ (cf.\ Appendix~B for definitions). 

\Lemma II.1. Let $G$ be a group, $A$ a $G$-module and $f \: G \times G
\to A$ a normalized $2$-cocycle, i.e., 
$$ f(g,\1) = f(\1,g) = \0, \quad 
f(g,g') + f(gg',gg') = g.f(g',g'') + f(g,g'g''), \quad g,g',g''\in
G.$$
Then we obtain a group $A \times_f G$ by endowing the product set 
$A \times G$ with the multiplication 
$$ (a,g) (a',g') := (a + g.a' + f(g,g'), gg'). \leqno(2.1) $$
The unit element of this group is $(\0,\1)$, inversion is given by 
$$ (a,g)^{-1} = (-g^{-1}.(a + f(g,g^{-1})),g^{-1}), \leqno(2.2)  $$
and conjugation by the formula 
$$ (a,g)(a',g') (a,g)^{-1} 
= \big(a + g.a' - gg'g^{-1}.a + f(g,g')
-f(gg'g^{-1},g),gg'g^{-1}\big). \leqno(2.3) $$
The map $q \: A \times_f G \to G, (a,g) \mapsto g$ is a surjective
homomorphism whose kernel $A \times \{\1\}$ is isomorphic to $A$. 
The conjugation action of $A \times_f G$ on the normal subgroup $A$
factors through the original action of $G$ on $A$. 

\Proof. The condition $f(\1,g) = f(g,\1) = \0$ implies that 
$(\0,\1)$ is an identity element in $A \times_f G$, and the associativity of
the multiplication is equivalent to the cocycle condition. 
The formula for the inversion is easily verified. 
Conjugation in $A \times_f G$ is given by 
$$ \eqalign{ (a,g)(a',g') (a,g)^{-1} 
&= \big(a + g.a' + f(g,g'), gg'\big) \big(-g^{-1}.(a + f(g,g^{-1})),g^{-1}\big)\cr
&= \big(a + g.a' + f(g,g')- gg'g^{-1}.(a + f(g,g^{-1})) +
f(gg',g^{-1}),gg'g^{-1}\big).\cr} $$
To simplify this expression, we use 
$$ f(g,g^{-1}) = f(g,g^{-1}) + f(\1,g) = f(g,\1) + g.f(g^{-1},g) = g.f(g^{-1},g) $$
and 
$$ f(gg',g^{-1}) + f(gg'g^{-1},g) = f(gg',\1) + gg'.f(g^{-1},g) = gg'.f(g^{-1},g) $$
to obtain 
$$ \eqalign{ (a,g)(a',g') (a,g)^{-1} 
&= \big(a + g.a' + f(g,g')- gg'g^{-1}.a - gg'g^{-1}.f(g,g^{-1}) +f(gg',g^{-1}),gg'g^{-1}\big)\cr
&= \big(a + g.a' + f(g,g')- gg'g^{-1}.a - gg'.f(g^{-1}.g) +f(gg',g^{-1}),gg'g^{-1}\big)\cr
&= \big(a + g.a' + f(g,g')- gg'g^{-1}.a -
f(gg'g^{-1},g),gg'g^{-1}\big).\cr} $$
In particular we obtain 
$$ (\0,g)(a,\1) (\0,g)^{-1} 
= (g.a,\1).  $$
This means that the action of $G$ on $A$ given by 
$q(g).a := gag^{-1}$ for $g \in A \times_f G$ 
coincides with the given action of $G$ on $A$. 
\qed

\Definition II.2. An {\it extension of Lie groups} is a 
surjective morphism $q \: \hat G \to G$ of Lie groups with a smooth
local section for which $N := \ker q$ has a natural Lie group
structure such that the map $N \times \hat G \to \hat G, (n,g) \mapsto
ng$ is smooth. Then the existence of a smooth local section
implies that $\hat G$ is a smooth $N$-principal bundle, so that 
$N$ is a split Lie subgroup of $G$ in the sense of Definition~C.4.

We call two extensions 
$N \into \hat G_1 \onto G$ and 
$N \into \hat G_2 \onto G$ of the Lie group $G$ by the 
Lie group $N$ {\it equivalent} if there exists a Lie group morphism 
$\phi \: \hat G_1 \to \hat G_2$ such that the following diagram commutes: 
$$ \matrix{
 N & \into& \hat G_1 & \onto & G \cr 
\mapdown{\id_N} & & \mapdown{\phi} & & \mapdown{\id_G} \cr 
 N & \into& \hat G_2 & \onto & G. \cr } $$
It is easy to see that any such $\phi$ is an isomorphism of 
 group and that its inverse is smooth. Thus $\phi$ is an isomorphism
 of Lie groups, and we obtain indeed an equivalence relation. 
We write $\Ext(G,N)$ for the set of equivalence classes of
Lie groups extensions of $G$ by~$N$. 
\qed

\Lemma II.3. If $A \into \hat G_1 \sssmapright{q_1} G$ and 
$A \into \hat G_2 \sssmapright{q_2} G$ are equivalent abelian extensions of $G$ by 
the Lie group $A$,
then the induced actions of $G$ on $A$ coincide. 

\Proof. There exists a morphism of Lie groups $\phi \: \hat G_1\to \hat
G_2$ with $\phi\res_A = \id_A$ and $q_2 \circ \phi = q_1$. 
For $g \in G$ and $a \in A$ the extension $\hat G_1$ defines an action
of $G$ on $A$ by $g *_1 a := g_1ag_1^{-1}$, where $q_1(g_1) = g$. We
likewise obtain from the extension $\hat G_2$ an action of $G$ on $A$
by $g *_2 a := g_2 a g_2^{-1}$ for $q_2(g_2) = g$. 
We then have 
$$ g *_1 a =  g_1ag_1^{-1} = \phi(g_1ag_1^{-1}) = \phi(g_1)a\phi(g_1)^{-1} = 
q_2(\phi(g_1)) *_2 a = q_1(g_1) *_2 a = g *_2 a. 
\qeddis 

\Definition II.4. If $(A,\rho_A)$ is a smooth $G$-module, then an {\it
extension of $G$ by $A$} is always understood to be an abelian Lie
group extension $q \: \hat G \to G$ with kernel $A$ for which the
natural action of $G$ on $A$ induced by the conjugation action
(Lemma C.5) coincides with $\rho_A$. 
In view of Lemma II.3, it makes sense to write 
$\Ext_{\rho_A}(G,A) \subeq \Ext(G,A)$ for the 
subset of equivalence classes of those extensions of $G$ by $A$ for
which the induced action of $G$ on $A$ coincides with~$\rho_A$. 
\qed

\Definition II.5. Let $G$ be a Lie group and $A$ a smooth
$G$-module. For $f \in Z^2_s(G,A)$ (cf.\ Definition~B.2) 
and $g \in G$ we consider the function 
$$ f_g \: G \to A, \quad f_g(g') := f(g,g') - f(gg'g^{-1},g) $$
and write 
$$ Z^2_{ss}(G,A) := \{ f \in  Z^2_s(G,A) \: (\forall g \in G) f_g
\in C^1_s(G,A)\} $$
for those locally smooth normalized $2$-cocycles $f$ on $G$ for which
all functions $f_g$ are smooth in an identity neighborhood of $G$. 

If $\ell \in C^1_s(G,A)$ and $f(g,g') = (d_G \ell)(g,g') 
= \ell(g) + g.\ell(g') - \ell(gg')$, then 
$$ \eqalign{ f_g(g') 
&= \ell(g) + g.\ell(g') - \ell(gg') 
- \big(\ell(gg'g^{-1}) + (gg'g^{-1}).\ell(g) - \ell(gg')\big)\cr
&= \ell(g) + g.\ell(g') - \ell(gg'g^{-1}) - (gg'g^{-1}).\ell(g)  \cr} $$
is smooth in an identity neighborhood of $G$ for each $g \in G$. Therefore 
$B^2_s(G,A) \subeq Z^2_{ss}(G,A)$ and 
$$ H^2_{ss}(G,A) :=  Z^2_{ss}(G,A)/B^2_s(G,A) $$
is a subgroup of $H^2_s(G,A)$. 
\qed

\Proposition II.6. Let $G$ be a Lie group and $(A, \rho_A)$ a smooth 
$G$-module. Then for each $f \in Z^2_{ss}(G,A)$ the group 
$A \times_f G$ carries the structure of a Lie group such that the map 
$q \: A \times_f G \to G, (a,g) \mapsto g$ is a Lie group extension of
$G$ by the smooth $G$-module $A$. Conversely, every Lie group
extension of $G$ by the smooth $G$-module $A$ is equivalent to one of this form. The
assignment 
$$ Z^2_{ss}(G,A) \to \Ext_{\rho_A}(G,A), \quad f \mapsto [A \times_f G] $$
factors through a bijection 
$$ H^2_{ss}(G,A) \to \Ext_{\rho_A}(G,A). $$

If $G$ is connected, then $Z^2_{ss}(G,A) = Z^2_s(G,A)$ and we obtain a
bijection 
$$ H^2_{s}(G,A) \to \Ext_{\rho_A}(G,A). $$

\Proof. {\rm(1)} Let $f \in Z^2_{ss}(G,A)$ and form the group $\hat G
:= A \times_f G$ (Lemma~II.1). First we construct the Lie group structure on $\hat G$. 
Let $U_G \subeq G$ be an open symmetric $\1$-neighborhood such that $f$
is smooth on $U_G \times U_G$, and consider the subset 
$$ U := A \times U_G = q^{-1}(U_G) \subeq \hat G = A \times_f G. $$
Then $U = U^{-1}$. We endow $U$ with the product manifold structure
from $A \times U_G$. Since the multiplication $m_G\res_{U_G \times U_G} \: U_G \times U_G \to
G$ is continuous, there exists an open identity neighborhood $V_G
\subeq U_G$ with $V_G V_G \subeq U_G$. Then the set 
$V := A \times V_G$ is an open subset of $U$ such that the multiplication map 
$$ V \times V \to U, \quad \big((a,x),(a',x')\big) \to (a + x.a' + f(x,x'), xx') $$
is smooth. The inversion 
$$ U \to U, \quad (a,x) \mapsto \big(-x^{-1}.(a + f(x,x^{-1})), x^{-1}\big) $$
(Lemma II.1) is also smooth. 

For $(a,g) \in \hat G$ let $V_g \subeq U_G$ be an open identity
neighborhood such that the conjugation map $c_g(x) = gxg^{-1}$
satisfies  
$c_g(V_g) \subeq U_G$. Then $c_{(a,g)}(q^{-1}(V_g)) \subeq U$ and the conjugation 
map 
$$ c_{(a,g)} \: q^{-1}(V_g) \to U, \quad 
(a',g') \mapsto (a + g.a' - gg'g^{-1}.a + f_g(g'), gg'g^{-1}) $$ 
(Lemma II.1) is smooth in an identity neighborhood because 
$f \in Z^2_{ss}(G,A)$. 

Now Theorem C.2 implies that $\hat G$ 
carries a unique Lie group structure for
which the inclusion map $U = A \times U_G \into \hat G$ is a {\bf local} 
diffeomorphism onto an identity neighborhood. 
It is clear that with respect to this Lie group structure on $\hat G$, the
map $q \: \hat G \to G$ defines a smooth $A$-principal bundle because the 
map $V_G \to \hat G, g \mapsto (0,g)$ defines a section of $q$ which
is smooth on an identity neighborhood in $G$ which might be smaller
than $V_G$. 

\par\nin (2) Assume, conversely, that $q \: \hat G \to G$ is an
extension of $G$ by the smooth $G$-module $A$. 
Then there exists an open $\1$-neighborhood $U_G \subeq G$ and a
smooth section $\sigma \: U_G \to \hat G$ of the map $q \: \hat G
\to G$. We extend $\sigma$ to a global section $G \to \hat G$. Then 
$$ f(x,y) := \sigma(x)\sigma(y)\sigma(xy)^{-1} $$
defines a $2$-cocycle $G \times G \to A$ which is smooth in a
neighborhood of $(\1,\1)$, and the map 
$$ A \times_f G \to \hat G, \quad (a,g) \mapsto a \sigma(g) $$
is an isomorphism of groups. The functions $f_g \: G \to A$ are given by 
$$ \eqalign{ f_g(g') 
&= f(g,g') - f(gg'g^{-1},g)= \sigma(g)\sigma(g')\sigma(gg')^{-1} 
- \sigma(gg'g^{-1})\sigma(g)\sigma(gg')^{-1} \cr
&= \sigma(g)\sigma(g')\sigma(gg')^{-1} 
\sigma(gg')\sigma(g)^{-1} \sigma(gg'g^{-1})^{-1} 
= \sigma(g)\sigma(g')\sigma(g)^{-1} \sigma(gg'g^{-1})^{-1}, \cr} $$
hence smooth near $\1$. This shows that $f \in Z^2_{ss}(G,A)$.  
In view of (1),  the group 
$A \times_f G$ carries a Lie group structure for which there exists an 
identity neighborhood $V_G \subeq G$ for which the product map 
$$A \times V_G \to A \times_f G, \quad (a,v) \mapsto (a,\1)(0,v) = (a,v) $$
is smooth. This implies that the group isomorphism 
$A \times_f G \to \hat G$ is a local diffeomorphism, hence an
isomorphism of Lie groups.

\par\nin (3) Steps (1) provides a map 
$$ Z^2_{ss}(G,A) \to \Ext_{\rho_A}(G,A), \quad f \mapsto [A \times_f
G], $$
and (2) shows that it is surjective. 
Assume that two extensions
of the form $A \times_{f_i} G$ for $f_1, f_2 \in Z^2_{ss}(G,A)$ are
equivalent as Lie group extensions. An isomorphism 
$A \times_{f_1} G \to A \times_{f_2} G$
inducing an equivalence of abelian extensions  must be of the form 
$$ (a,g) \mapsto (a + h(g), g), \leqno(2.4) $$
where $h \in C^1_s(G,A)$. 
The condition that (2.4) is a group homomorphism implies that 
$$ (h(gg') + f_1(g,g'),gg') = (h(g), g)(h(g'),g') = (h(g) + g.h(g') + f_2(g,g'), gg'), $$
which means that 
$$ (f_1 - f_2)(g.g') = g.h(g') - h(gg') + h(g)  = (d_G h)(g,g'),
\leqno(2.5) $$
so that $f_1 - f_2 \in B^2_s(G,A)$. 

If, conversely, $h \in C^1_s(G,A)$ and $f_1 = f_2 = d_G h$, then it is
easily verified that (2.4) defines a group
isomorphism for which there exists an open identity neighborhood mapped 
diffeomorphically onto its image. Hence (2.5) is an isomorphism of Lie
groups. We conclude that the map $Z^2_{ss}(G,A) \to \Ext_{\rho_A}(G,A)$
factors through a bijection $H^2_{ss}(G,A) \to \Ext_{\rho_A}(G,A)$. 

\par\nin (4) Assume now that $G$ is connected and that $f \in
Z^2_s(G,A)$.  In the context of (1), the conjugation map 
$c_{(a,g)} \: q^{-1}(V_g) \to U$
is smooth in an identity neighborhood if and only if the function 
$f_g$ is smooth in an identity neighborhood. As $f \in Z^2_s(G,A)$,
the set $W$ of all $g \in G$ for which this condition is satisfied 
is an identity neighborhood. On the other hand, the set $W$ is closed
under multiplication. In view of the
connectedness of $G$, we have $G = \bigcup_{n \in \N} W^n = W$. 
This means that $f \in Z^2_{ss}(G,A)$, and therefore that 
$Z^2_s(G,A) = Z^2_{ss}(G,A)$. 
\qed

\Problem II. Do the two spaces $Z^2_s(G,A)$ and $Z^2_{ss}(G,A)$
also coincide if $G$ is not connected? We do not know any cocycle 
$f \in Z^2_{s}(G,A) \setminus Z^2_{ss}(G,A)$. 
\qed

The following lemma shows that the derivation map 
$$ D \: Z^2_s(G,A) \to Z^2_c(\g,\a), \quad 
(Df)(x,y) = d^2f(\1,\1)(x,y) - d^2f(\1,\1)(y,x)  $$
from Theorem B.6 and Lemma B.7 
is compatible with the construction in Proposition~II.6. 
In the following proof we use the notation $d^2 f$ introduced 
in Appendix~A. 

\Lemma II.7. Let $A \cong \a/\Gamma_A$, where $\Gamma_A \subeq \a$ is a discrete subgroup,  
$f \in Z^2_{ss}(G,A)$ and $\hat G = A \times_f G$ the corresponding 
extension of $G$ by $A$. Then the Lie algebra cocycle 
$Df$ satisfies $\hat \g \cong \g \oplus_{Df} \a$. 

\Proof. Let $U_\a \subeq \a$ be an open $0$-neighborhood such that the
restriction $\phi_A \: U_\a \to U_\a + \Gamma_A \subeq A$ of the
quotient map $q_A \: \a \to A$ is a diffeomorphism onto an open identity
neighborhood in $A$ and $\phi_G \: U_\g \to G$ a local chart of $G$,
where $U_\g \subeq \g$ is an open $0$-neighborhood, $\phi_G(0) = \1$
and $d\phi_G(0) =\id_\g$. After shrinking $U_\g$ further, we 
obtain a chart of $A \times_f G$ by the map 
$$ \phi \: U_\a \times U_\g \to A \times_f G, \quad 
(a,x)  \mapsto (\phi_A(a), \phi_G(x)). $$
Moreover, we may assume that $U_\g$ is so small that 
$f(\phi_G(U_\g) \times \phi_G(U_\g)) \subeq \phi_A(U_\a))$, which
implies that there exists a smooth function $f_\a \: U_\g \times U_\g \to U_\a$
with $\phi_A \circ f_\a = f \circ (\phi_G \times \phi_G).$ 

Writing $x * x' := \phi_G^{-1}(\phi_G(x)\phi_G(x'))$ for $x,x' \in
U_\g$ with $\phi_G(x)\phi_G(x') \in \phi_G(U_\g)$, 
the multiplication 
$$ (a,g)(a',g') = (a + g.a' + f(g,g'), gg') $$
in $A \times_f G$ can be expressed in local coordinates for
sufficiently small $a,a' \in \a, x,x' \in \g$ by  
$$ \eqalign{ \phi(a,x) \phi(a',x') 
&= (\phi_A(a) + \phi_G(x).\phi_A(a') + f(\phi_G(x), \phi_G(x')), \phi_G(x)\phi_G(x'))\cr
&= (\phi_A(a + \phi_G(x).a' + f_\a(x,x')), \phi_G(x * x'))\cr
&= \phi(a + \phi_G(x).a' + f_\a(x,x'), x * x'). \cr} $$
Here the identity element has the coordinates $(0,0) \in \a \times\g$. 

For the multiplication in $G$ we have 
$$ x * x' = x + x' + b(x,x') + \cdots $$
where $\cdots$ stands for the terms of order at least three in the
Taylor expansion of the product map and the quadratic term is 
bilinear. The Lie bracket in $\g$ is given by 
$$ [x,x'] = b(x,x') - b(x',x) $$
([Mil83, p.1036]). Therefore the Lie bracket in the Lie algebra 
$\L(A \times_f G)$ of $A \times_f G$ can be obtained from 
$$ \eqalign{ 
&\ \ \ \ (a + \phi_G(x).a' + f_\a(x,x'), x * x') \cr
&= (a + a' + x.a' + d^2 f_\a(0,0)(x,x') + \cdots, x + x' +b(x,x') +
\cdots) \cr
&= (a + a' + x.a' + d^2 f(\1,\1)(x,x') + \cdots, x + x' +b(x,x') + \cdots),  \cr} $$
as 
$$ [(a,x), (a',x')] 
= (x.a' - x'.a + Df(x,x'), [x,x']). 
\qeddis 

\sectionheadline{III. Locally smooth $1$-cocycles} 

Let $G$ be a Lie group and $A$ a smooth $G$-module. In this section we
take a closer look at the space $Z^1_s(G,A)$ of locally smooth
$A$-valued $1$-cocycles on $G$. We know from Appendix B that there is
a natural map 
$$ D_1 \: Z^1_s(G,A) \to Z^1_c(\g,\a), \quad 
D_1(f)(x) := df(\1)(x). $$
If $A \cong \a/\Gamma_A$ holds for a discrete subgroup $\Gamma_A$ of
$\a$ and $q_A \: \a \to A$ is the quotient map, then we have for 
$a \in \a$ the relation 
$$ D_1(d_G(q_A(a))) = d_\g(a) $$
and hence $D_1(B^1_s(G,A)) = B^1_c(\g,\a)$. Hence $D_1$ 
induces a map 
$$ D_1 \: H^1_s(G,A) \to H^1_c(\g,\a), $$
and it is of fundamental interest to have a good description of 
kernel and cokernel of $D_1$ on the level of cocycles and cohomology classes. 

We shall see that the integration problem
for Lie algebra $1$-cocycles has a rather simple solution, the only
obstruction coming from $\pi_1(G)$. 

\Lemma III.1. Each $f \in Z^1_s(G,A)$ is a smooth function and its
differential $df \in \Omega^1(G,\a)$ is an equivariant $1$-form. 

\Proof. Let $g \in G$. The cocycle condition 
$$ f(gh) = g.f(h) + f(g) \leqno(3.1) $$
shows that the smoothness of $f$ in an identity neighborhood implies
the smoothness in a neighborhood of $g$. 

Formula (3.1) means that 
$f \circ \lambda_g = \rho_A(g) \circ f +f(g),$
so that $df$ satisfies 
$\lambda_g^* df = \rho_A(g) \circ df,$
i.e., $df$ is equivariant. 
\qed

\Lemma III.2. Let $G$ be a Lie group with identity component $G_0$ and $A$ a smooth $G$-module. 
Then for a smooth function 
$f \: G \to A$ with $f(\1) = 0$ the following are equivalent: 
\litem{(1)} $df$ is an equivariant $\a$-valued $1$-form on $G$. 
\litem{(2)} $f(gn) = f(g) + g.f(n)$ for $g \in G$ and $n \in G_0$. 

\par\nin If, in addition, $G$ is connected, then $df$ is equivariant if and
only if $f$ is a cocycle. 

\Proof. We write $g.a = \rho_\a(g).a$ for the action of $G$ on $\a$
and $g.a = \rho_A(g).a$ for the action of $G$ on $A$. 

\par\nin (1) $\Rarrow$ (2): Let $g \in G$. 
In view of 
$d(\rho_A(g) \circ f) = \rho_\a(g) \circ df,$
we have 
$$ d(f \circ \lambda_g - \rho_A(g) \circ f - f(g)) 
 = \lambda_g^* df - \rho_\a(g) \circ df. $$
Hence (1) means that all the functions $f \circ \lambda_g - \rho_A(g)
 \circ f - f(g)$ are locally constant. 
Since the value of this function in $\1$ is $0$, 
all these functions are constant $0$ on $G_0$, which is (2). 

\par\nin (2) $\Rarrow$ (1): If (2) is satisfied, then 
$df(g) d\lambda_g(\1) = \rho_\a(g) \circ df(\1)$
holds for each $g \in G$, and this means that $df$ is equivariant. 
\qed

\Definition III.3. Suppose that $\a$ is sequentially complete. 
If $\alpha \in Z^1_c(\g,\a)$ and $\alpha^{\rm eq}$ is the corresponding closed equivariant 
$1$-form on $G$ (cf.\ Definition B.4), 
then we obtain a morphism of abelian groups, called
the {\it period map of $\alpha$:}   
$$ \per_\alpha \: \pi_1(G) \to \a, 
\quad 
[\gamma] \mapsto \int_\gamma \alpha^{\rm eq} 
=  \int_0^1 \alpha^{\rm eq}_{\gamma(t)}(\gamma'(t))\, dt 
=  \int_0^1 \gamma(t).\alpha(\gamma(t)^{-1}.\gamma'(t))\, dt, $$
where $\gamma \: [0,1] \to G$ is a piecewise smooth loop based in
$\1$. The map 
$$ C^\infty(\SS^1, G) \to \a, \quad \gamma \mapsto \int_\gamma \alpha^{\rm eq} $$
is locally constant, so that the connectedness of $G$ implies in particular that for 
$g \in G$ we have 
$$ \int_{\gamma} \alpha^{\rm eq}  
= \int_{\lambda_g \circ \gamma} \alpha^{\rm eq}  
= \int_{\gamma} \lambda_g^*\alpha^{\rm eq}  
= \rho_\a(g).\int_{\gamma} \alpha^{\rm eq} $$
which leads to 
$$ \im(\per_\alpha) \subeq \a^G. $$

If $\Gamma_A \subeq \a^G$ is a discrete subgroup, then $A :=
\a/\Gamma_A$ is a smooth $G$-module with respect to the induced
action. Let $q_A \: \a \to A$ denote the quotient map. 
We then obtain a group homomorphism 
$$ P \: Z^1_c(\g,\a) \to \Hom(\pi_1(G),A^G), \quad 
P(\alpha) := q_A \circ \per_\alpha. $$
The importance of the period map stems from the fact that 
the $1$-form $\alpha^{\rm eq}$ is the
differential of a smooth function $f \: G \to A$ if and only if
$P(\alpha)=0$ ([Ne02, Prop.~3.9]). 
\qed

\Proposition III.4. If $G$ is a connected Lie group and $A_0 \cong
\a/\Gamma_A$, where $\Gamma_A \subeq \a^G$ is a discrete subgroup 
and $\a$ is sequentially complete, then the sequence 
$$ \0 \to Z^1_s(G,A) \sssmapright{D_1} Z^1_c(\g,\a) \sssmapright{P}
\Hom(\pi_1(G), A^G) \leqno(3.2) $$
is exact. If $A \cong \a/\Gamma_A$, then it induces an exact sequence 
$$ \0 \to H^1_s(G,A) \sssmapright{D_1} H^1_c(\g,\a) \sssmapright{P}  
\Hom(\pi_1(G), A^G). \leqno(3.3) $$

\Proof. If $f \in Z^1_s(G,A)$ satisfies $D_1 f = 0$, then Lemma III.2 implies that 
$df = 0$ because $df$ is equivariant, and hence that 
$f$ is constant, and we get $f(g) = f(\1) = 0$ for each $g \in G$. Therefore $D_1 $ is injective on 
$Z^1_s(G,A)$. The kernel of 
$P \: Z^1_c(\g,\a) \to \Hom(\pi_1(G),A)$ 
consists of those $1$-cocycles $\alpha$ for which $\alpha^{\rm eq}$ is
the differential of a smooth function $f \: G \to A$ with $f(\1) =0$ ([Ne02, Prop.~3.9]), 
which means that $\alpha = D_1 f$ for 
some $f \in Z^1_s(G,A)$ (Lemma III.2).  This proves the exactness of the first sequence. 

Now we assume that $A \cong \a/\Gamma_A$.  
If $\alpha \in B^1_c(\g,\a)$, then $\alpha^{\rm eq}$ is exact
(Lemma~B.5), 
so that $P(\alpha)= 0$. Therefore 
$P$ factors through a map $H^1_c(\g,\a) \to \Hom(\pi_1(G),\a)$. 
The exactness of (3.3) follows from the observation that 
$D_1 (B^1_s(G,A)) = B^1_c(\g,\a)$ and the
exactness of (3.2). 
\qed

\Remark III.5. For each $\alpha \in Z^1_c(\g,\a)$ the corresponding
equivariant $1$-form $\alpha^{\rm eq}$ is closed and it is exact if
$\alpha \in B^1_c(\g,\a)$, so that we obtain a map 
$$ H^1_c(\g,\a)  \to H^1_{\rm dR}(G,\a), \quad [\alpha] \mapsto
[\alpha^{\rm eq}]. $$
Proposition III.4, applied to the case $A = \a$ now means that the
sequence 
$$ \0 \to H^1_s(G,\a) \sssmapright{D_1} H^1_c(\g,\a) \sssmapright{} 
H^1_{\rm dR}(G,\a) $$
is exact. Let $\Gamma_A \subeq \a$ be a discrete subgroup and consider
$A := \a/\Gamma_A$. For 
$$ H^1_{\rm dR}(G,\Gamma_A) := \Big\{ [\alpha] \in H^1_{\rm dR}(G,\a) \: 
(\forall \gamma \in C^\infty(\SS^1,G))\, \int_\gamma \alpha
\in\Gamma_A\Big\}, $$
we then have 
$$ H^1_{\rm dR}(G,\Gamma_A) = dC^\infty(G,A)/dC^\infty(G,\a) $$
([Ne02, Prop.~3.9]), and we obtain an exact sequence 
$$ \0 \to H^1_s(G,A) \sssmapright{D_1} H^1_c(\g,\a) \sssmapright{} 
\Omega^1(G,\a)/dC^\infty(G,A), \leqno(3.4) $$
because for $\alpha \in Z^1_c(\g,\a)$ the condition 
$[\alpha^{\rm eq}] \in dC^\infty(G,A)$ is equivalent to $P([\alpha])
= 0$ (Proposition~III.4). 
\qed

\Definition III.6. Let $A$ be a smooth $G$-module for the connected Lie group $G$ 
and assume that $A_0 \cong \a/\Gamma_A$ holds for the identity component of $A$. 
Then for each $a \in A$ we obtain a smooth cocycle 
$$ d_G'(a) \in Z^1_s(G,A_0), \quad d_G'(a)(g) := g.a - a. $$
Taking derivatives in $\1$ leads to homomorphisms 
$$ \theta_A := D_1 \circ d_G' \: A \to Z^1_c(\g,\a) 
\quad \hbox{ and } \quad \oline\theta_A  \: \pi_0(A) \to H^1_c(\g,\a). $$
The map $\oline\theta_A$ is called the {\it characteristic homomorphism of the 
$G$-module $A$}. 
\qed

\Lemma III.7. Let $A$ and $B$ be smooth modules of the connected Lie group 
$G$ and assume that $A_0 = B_0 \cong \a/\Gamma_A$ as $G$-modules, where 
$\Gamma_A \subeq \a$ is a discrete subgroup. Then there exists an 
isomorphism $\psi \: A \to B$ of $G$-modules with $\psi\res_{A_0} = \id_{A_0}$ 
if and only if there exists a homomorphism 
$\gamma \: \pi_0(A) \to \pi_0(B)$ such that the characteristic homomorphisms 
of $A$ and $B$ are related by 
$$ \oline\theta_B \circ \gamma = \oline\theta_A. $$

\Proof. If $\psi \: A \to B$ is an isomorphism of $G$-modules restricting 
to the identity on $A_0$, then $\psi$ induces an isomorphism 
$\gamma := \pi_0(\psi) \: \pi_0(A) \to \pi_0(B)$, and it follows directly from the 
definitions that $\oline\theta_B \circ \gamma = \oline\theta_A.$

Suppose, conversely, that $\gamma \: \pi_0(A) \to \pi_0(B)$ is an isomorphism with 
$\oline\theta_B \circ \gamma = \oline\theta_A$. 
Since $A_0$ is an open divisible subgroup of $A$, we have 
$A \cong A_0 \times \pi_0(A)$ as abelian Lie groups and likewise 
$B \cong A_0 \times \pi_0(B)$. For each homomorphism 
$\phi_0 \: \pi_0(A) \to A_0$ we then obtain a Lie group isomorphism 
$$ \phi \: A \to B, \quad (a_0, a_1) \mapsto (a_0 + \phi_0(a_1), \gamma(a_1)). 
\leqno(3.5) $$
Since $G$ acts on $A \cong A_0 \times \pi_0(A)$ by 
$$ g.(a_0, a_1) = (g.a_0 + d_G'(a_1)(g), a_1), $$
the isomorphism $\phi$ is $G$-equivariant if and only if 
$$ \phi_0(a_1) + d_G'(a_1)(g) = g.\phi_0(a_1) + d_G'(\gamma(a_1))(g) \leqno(3.6) $$
for $g \in G$, $a_1 \in \pi_0(A)$, which means that 
$$ d_G(\phi_0(a_1)) = d_G'(a_1) - d_G'(\gamma(a_1)) =: \beta. $$
To see that a homomorphism $\phi_0$ with the required properties exists, 
we first observe that our assumption implies that $\beta$ 
is a homomorphism $\pi_0(A) \to Z^1_s(G,A_0)$ with 
$\im(D_1 \circ \beta) \subeq d_\g \a$. In view of the divisibility of $\a$, 
there exists a homomorphism 
$\delta \: \pi_0(A) \to \a$ with $D_1 \circ \beta 
= d_\g \circ \delta = D_1 \circ d_G \circ q_A \circ \delta$.
Since $D_1$ is injective on cocycles (Proposition III.4), we obtain 
$\beta = d_G \circ q_A \circ \delta$.  We may therefore put 
$\phi_0 := q_A \circ \delta$ to obtain an isomorphism $\phi$ of $G$-modules 
as in (3.5). 
\qed

\sectionheadline{IV. Period homomorphisms for abelian groups} 

\nin In this section $G$ denotes a connected Lie group, 
$\a$ is a smooth sequentially complete $G$-module, and $\omega \in Z^2_c(\g,\a)$ is a
continuous Lie algebra cocycle. We shall define a homomorphism of
abelian groups 
$$ \per_\omega \: \pi_2(G) \to \a, $$
called the period map of $\omega$. 

Suppose that $q \: \hat G \to G$ is an extension of $G$
by the smooth $G$-module $A$ whose Lie algebra is isomorphic to 
$\a \oplus_\omega \g$ and $A_0 \cong \a/\Gamma_A$ holds
for a discrete subgroup $\Gamma_A \cong \pi_1(A)$ of $\a$. Then 
we show that the period map is, up to
sign, the connecting map of the long exact homotopy sequence of the 
principal $A$-bundle 
$A \into \hat G \onto G$, whose range is contained in the
subgroup $\Gamma_A \subeq \a$. 

\Definition IV.1. In the following 
$\Delta_p = \{ (x_1,\ldots, x_p) \in \R^p \: x_i \geq 0, \sum_j x_j \leq 1\}$
denotes the {\it {$p$-dimensional} standard simplex} in $\R^p$. 
We also write $\la v_0, \ldots, v_p \ra$ for the affine simplex in a vector space 
spanned by the points $v_0, \ldots, v_p$. In this sense 
$\Delta_p = \la 0, e_1, \ldots, e_p\ra$, where $e_i$ denotes the
$i$-th canonical basis vector in $\R^p$. 

Let $Y$ be a smooth manifold. 
A continuous map $f \: \Delta_p  \to Y$ is called {\it a $C^1$-map} 
if it is differentiable in the interior $\Int(\Delta_p)$ and in each
local chart of $Y$ all directional derivatives $x \mapsto df(x)(v)$ 
of $f$ extend continuously to the boundary $\partial \Delta_p$ of
$\Delta_p$. For $k \geq 2$ we call $f$ a $C^k$-map if it is $C^1$ and 
all maps $x \mapsto df(x)(v)$
are $C^{k-1}$, and we say that $f$ is {\it smooth} if $f$ is $C^k$ for
every $k \in \N$. 
We write $C^\infty(\Delta_p, Y)$ for the set of smooth maps $\Delta_p \to Y$. 

If $\Sigma$ is a simplicial complex, then we call a map $f \: \Sigma \to Y$ 
{\it piecewise smooth} if it is continuous and its restrictions to all simplices in 
$\Sigma$ are smooth. We write $C^\infty_{pw}(\Sigma, Y)$ for the set of piecewise smooth maps 
$\Sigma \to Y$. There is a natural topology on this space inherited from the 
natural embedding of $C^\infty_{pw}(\Sigma,Y)$ into the space 
$\prod_{S \subeq \Sigma} C^\infty_{pw}(S,Y)$, where $S$ runs through all simplices 
of $\Sigma$ and the topology on $C^\infty_{pw}(S,Y)$ is defined as in 
[Ne02, Def.~A.3.5] as the topology of uniform convergence of all directional derivatives 
of arbitrarily high order. 
\qed

The equivariant form $\omega^{\rm eq}$ is a closed $2$-form on $G$, and we obtain with 
[Ne02, Lemma 5.7] a period map 
$$ \per_\omega \: \pi_2(G) \to \a $$
which is given on piecewise smooth representatives 
$\sigma \: \SS^2\to G$ of free homotopy classes by the integral 
$$ \per_\omega([\sigma]) = \int_{\SS^2} \sigma^*\omega^{\rm eq} 
= \int_{\sigma} \omega^{\rm eq}. $$
If $\omega$ is a coboundary, then Lemma B.5 implies that $\omega^{\rm eq}$
is exact, so that the period map is trivial by Stoke's Theorem. We
therefore obtain a homomorphism 
$$ H^2_c(\g,\a) \to \Hom(\pi_2(G), \a), \quad [\omega] \mapsto
\per_\omega. $$ 
The image $\Pi_\omega := \per_\omega(\pi_2(G))$ is called the {\it
period group} of $\omega$. 
Since the group $G$ is connected, the group $\pi_0(C^\infty(\SS^2,G))$ 
of connected components of
the Lie group $C^\infty(\SS^2,G)$ is isomorphic to $\pi_2(G)$, and we
may think of $\per_\omega$ as the map on $\pi_2(G)$ obtained by
factorization of the map 
$$ C^\infty(\SS^2,G) \to \a, \quad \sigma \mapsto \int_\sigma
\omega^{\rm eq} $$
which is locally constant ([Ne02, Lemma 5.7]). 

\Lemma IV.2. The image of the period map is fixed pointwise by $G$,
i.e., $\Pi_\omega \subeq \a^G$. 

\Proof. In view of [Ne02, Th.~A.3.7], each homotopy class in
$\pi_2(G)$ has a smooth representative $\sigma \: \SS^2\to G$. 
Since $G$ is connected, and the
map $G \to C^\infty(\SS^2,G), g \mapsto \lambda_g \circ \sigma$ is
continuous, we have for each $g \in G$: 
$$ \per_\omega([\sigma]) 
= \int_{\SS^2} \sigma^*\omega^{\rm eq}
= \int_{\SS^2} \sigma^*\lambda_g^*\omega^{\rm eq}
= \int_{\SS^2} \sigma^*(\rho_\a(g) \circ \omega^{\rm eq}) 
= \rho_\a(g).\per_\omega([\sigma]). $$
We conclude that the image of $\per_\omega$ is fixed pointwise by $G$. 
\qed

\subheadline{Period maps as connecting homomorphisms} 

Let $A \into \hat G \sssmapright{q} G$ be an abelian Lie group
extension of $A$. Then the Lie algebra $\hat\g$ of $\hat G$ has the
form $\a \oplus_\omega \g$ because the existence of a smooth local
section implies that $\hat\g \to \g$ has a continuous linear section (Proposition~I.3). 
In this subsection we will relate the period homomorphism 
$\per_\omega$ to the connecting homomorphism $\delta \: \pi_2(G) \to \pi_1(A)$ 
from the long exact homotopy sequence of the bundle $A \into \hat G
\sssmapright{q} G$. 

\Definition IV.3. We recall the definition of {\it relative homotopy groups}. 
Let $I^n := [0,1]^n$ denote the $n$-dimensional cube. Then the
boundary $\partial I^n$ of $I^n$ can be written as 
$I^{n-1} \cup J^{n-1}$, where $I^{n-1}$ is called the {\it initial
face} and $J^{n-1}$ is the union of all other faces of $I^n$. 

Let $X$ be a topological space, $Y \subeq X$ a subspace, and $x_0 \in
Y$. A {\it map}
$$f \: (I^n, I^{n-1}, J^{n-1}) \to (X,Y,x_0) $$
of space triples is a continuous map $f \: I^n \to X$ satisfying 
$f(I^{n-1}) \subeq Y$ and $f(J^{n-1}) = \{x_0\}$. We write 
$F^n(X,Y,x_0)$ for the set of all such maps and 
$\pi_n(X,Y,x_0)$ for the homotopy classes of such maps, i.e., the
arc-components of the topological space $F^n(X,Y,x_0)$ endowed
with the compact open topology (cf.\ [Ste51]). 
We define $F^n(X,x_0) := F^n(X,\{x_0\},x_0)$ and 
$\pi_n(X,x_0) := \pi_n(X,\{x_0\},x_0)$ and observe 
that we have a canonical map 
$$ \partial \: \pi_n(X,Y,x_0) \to \pi_{n-1}(Y,x_0), \quad 
[f] \mapsto [f\res_{I^{n-1}}]. 
\qeddis 

\Example IV.4. Let $q \: P \to M$ be a (locally trivial) $H$-principal bundle, 
$y_0 \in P$ a base point, $x_0 := q(y_0)$ and identity $H$ with the fiber $q^{-1}(x_0)$. 
Then the maps 
$$ q_* \: \pi_k(P,H) := \pi_k(P,H,y_0) \to \pi_k(M) := \pi_k(M,x_0), 
\quad [f] \mapsto [q \circ f] $$
are isomorphisms ([Ste51, Cor.\ 17.2]), so that we obtain 
connecting homomorphisms 
$$ \delta := \partial \circ (q_*)^{-1} \: \pi_k(M) \to \pi_{k-1}(H). $$  
The so obtained sequence 
$$  \to \pi_k(P) \to \pi_k(M) \to \pi_{k-1}(H) \to 
\ldots \to \pi_1(P) \to \pi_1(M) \to \pi_0(H) \to \pi_0(P) \onto \pi_0(M) $$
is exact, where the last two maps cannot be considered as group homomorphisms. 
This sequence is called the {\it long exact homotopy sequence of the principal bundle $P \to M$}. 
\qed
 
\Proposition IV.5. Let $q \: \hat G \to G$ be an abelian extension 
of not necessarily connected 
Lie groups with kernel $A$ satisfying 
$A_0 \cong \a/\Gamma_A$, where $\a$ is a sequentially complete locally convex space. 
Then $q$ defines in
particular the structure of an $A$-principal bundle on $\hat G$. 
If $\omega \in Z^2_c(\g,\a)$ is a Lie algebra 
$2$-cocycle with $\hat\g \cong \a \oplus_\omega \g$, then 
$\delta \: \pi_2(G) \to \pi_1(A)$ and the period map 
$\per_\omega \: \pi_2(G) \to \a$ are related by 
$$ \delta = - \per_\omega \: \pi_2(G) \to \pi_1(A) \subeq \a.  $$

\Proof. Let $\Theta \in \Omega^1(\hat G,\a)$ be a $1$-form with the property that
for each $g \in \hat G$ the orbit map 
$\eta_g \: A \into \hat G, a \mapsto ga$ satisfies 
$\eta_g^*\Theta = \theta_A$, where $\theta_A \in \Omega^1(A,\a)$ is the 
invariant $1$-form on $A$ with $\theta_{A,\0} = \id_\a$, i.e., the Maurer-Cartan form on $A$. 
We have seen in [Ne02, Prop.~5.11] that if 
$\Omega \in \Omega^2(G,\a)$ satisfies  
$q^*\Omega = - d\Theta$, then $\delta = -\per_\Omega$. 

To apply this to our situation, we consider the action of $\hat G$ on $A$ given by 
$g.a := q(g).a$. Then $q^*\omega^{\rm eq}$ is an equivariant closed
$2$-form on $\hat G$ with 
$(q^*\omega^{\rm eq})_\1 = \L(q)^*\omega.$
Let $p_\a \: \hat\g\cong \a \oplus_\omega \g\to \a, (a,x) \mapsto x$ denote the 
projection onto $\a$. Then 
$$ \eqalign{ dp_\a((a,x),(a',x')) 
&= (a,x).p_\a(a',x') - (a',x').p_\a(a,x) - p_\a([(a,x),(a',x')])\cr
&= x.a' - x'.a  - (x.a' - x'.a + \omega(x,x')) = - \omega(x,x') \cr
&= -(\L(q)^*\omega)((a,x),(a',x')). \cr} $$
In view of Lemma B.5, this implies 
$$d(p_\a^{\rm eq}) = (d_\g p_\a)^{\rm eq}  = - (\L(q)^*\omega)^{\rm
  eq} = -q^*\omega^{\rm eq}.$$ 
Applying the preceding remarks with $\Theta = p_\a^{\rm eq}$, we 
obtain $\delta = - \per_\omega.$
\qed

\Remark IV.6. Let $A \into \hat G \onto G$ be an abelian 
extension of connected Lie groups and assume that $A \cong \a/\Gamma_A$
holds for a discrete subgroup $\Gamma_A \subeq \a$, that we 
identify with $\pi_1(A)$. In view of $\pi_2(A) \cong \pi_2(\a) =
\1$, the long exact homotopy
sequence of the bundle $\hat G \to G$ leads to an exact sequence 
$$ \0 \to \pi_2(\hat G) \into \pi_2(G) \ssmapright{\per_\omega} 
\pi_1(A) \to \pi_1(\hat G) \onto \pi_1(G) \to \0.$$
This implies that 
$$ \pi_2(\hat G) \cong \ker \per_\omega \subeq \pi_2(G) 
\quad \hbox{ and } \quad  \pi_1(G) \cong \pi_1(\hat G)/\coker
\per_\omega. $$
These relations show how the period homomorphism controls how the first two
homotopy groups of $G$ and $\hat G$ are related. 
\qed

\sectionheadline{V. From Lie algebra cocycles to group cocycles} 

\nin In Sections V and VI we describe the image of the map 
$$D \: H^2_s(G,A) \to H^2_c(\g,\a)$$ for a 
connected Lie group $G$, and an abelian Lie group $A$ of the form $\a/\Gamma_A$. 
In the present section we deal with the special case where $G$ is
simply connected. 

Let $G$ be a connected simply connected Lie group 
and $\a$ a sequentially complete locally convex smooth $G$-module.
Further let $\Gamma_A \subeq \a^G$ be a subgroup and write 
$A := \a/\Gamma_A$ for the quotient group, that carries a natural
$G$-module structure. We write $q_A \: \a \to A$ for the quotient map. 
If, in addition, $\Gamma_A$ is discrete, then 
$A$ carries a natural Lie group structure and the action of $G$ on $A$ is smooth. 

Let $\omega \in Z^2_c(\g,\a)$ and $\Pi_\omega \subeq \a^G$ be the corresponding 
period group (Lemma~IV.2). In the following we shall assume that 
$$\Pi_\omega \subeq \Gamma_A. $$
The main result of the present section is the existence of a
locally smooth group cocycle $f \in Z^2_s(G,A)$ with $Df = \omega$ if
$\Gamma_A$ is discrete (Corollary V.3). 

A special case of the 
following construction has also been used in [Ne02] in the context
of central extensions. 
For $g \in G$ we choose a smooth path $\alpha_{\1,g} \: [0,1] \to G$ 
from $\1$ to $g$. 
We thus obtain a left invariant system of smooth arcs 
$\alpha_{g,h} := \lambda_g \circ \alpha_{\1,g^{-1}h}$ from $g$ to $h$, where
$\lambda_g(x) := gx$ denotes left translation. 
For $g,h,u \in G$ we then obtain a singular smooth cycle 
$$ \alpha_{g,h,u} := \alpha_{g,h} + \alpha_{h,u} - \alpha_{g,u}, $$
that corresponds to the piecewise smooth map 
$\alpha_{g,h,u} \in C^\infty_{pw}(\partial \Delta_2, G)$ with 
$$ \alpha_{g,h,u}(s,t) = 
\cases{ 
\alpha_{g,h}(s), & for $t = 0$ \cr 
\alpha_{h,u}(1-s), & for $s + t = 1$ \cr 
\alpha_{g,u}(t), & for $s = 0$. \cr} $$

For a simplicial complex $\Sigma$ we write $\Sigma_{(j)}$ for the $j$-th 
{\it barycentric subdivision of $\Sigma$}. 
According to [Ne02, Prop.~5.6], each map $\alpha_{g,h,u}$ can be 
obtained as the restriction of a piecewise smooth map 
$\sigma \: (\Delta_2)_{(1)} \to G$. 
Let $\sigma' \: (\Delta_2)_{(1)} \to G$ be another piecewise smooth map 
with the same boundary values as $\sigma$. We claim that  
$\int_\sigma \omega^{\rm eq} - \int_{\sigma'} \omega^{\rm eq} \in \Pi_\omega$. 
In fact, we consider the sphere $\SS^2$ as an oriented simplicial complex $\Sigma$ 
obtained by gluing 
two copies $D$ and $D'$ of $\Delta_2$ along their boundary, where the inclusion of 
$D$ is orientation preserving and the inclusion on $D'$ reverses orientation. 
Then $\sigma$ and $\sigma'$ combine to a piecewise smooth map 
$\gamma \: \Sigma \to G$  with $\gamma\res_{D} = \sigma$ and 
$\gamma\res_{D'} = \sigma'$, and we get with [Ne02, Lemma 5.7]  
$$ \int_\sigma \omega^{\rm eq} - \int_{\sigma'} \omega^{\rm eq} = \int_\gamma \omega^{\rm eq} 
\in \Pi_\omega \subeq \Gamma_A.$$
We thus obtain a well-defined map 
$$ F{} \: G^3 \to A, \quad (g,h,u) \mapsto q_A\Big(\int_\sigma \omega^{\rm eq}\Big), $$
where $\sigma \in C^\infty_{pw}((\Delta_2)_{(1)}, G)$ is a piecewise smooth map whose boundary values
coincide with $\alpha_{g,h,u}$. 

\Lemma V.1. The function
$$f{} \: G^2 \to A, \quad (g,h) \mapsto F(\1,g,gh) $$
is a group cocycle with respect to the action of $G$ on $A$. 

\Proof. First we show that for $g,h \in G$ we have  
$$ f(g,\1) = F(\1,g,g) = 0 \quad \hbox{ and } \quad 
f(\1,h) = F(\1,\1,h)= 0. $$
If $g = h$ or $h = u$, then
we can choose the map $\sigma\: \Delta_2 \to G$ extending 
$\alpha_{g,h,u}$ in such a way that $\rk(d\sigma) \leq 1$ in every point, so that 
$\sigma^*\omega^{\rm eq} =0$. In particular we obtain $F(g,h,u) = 0$ in these cases. 

From $\alpha_{g,h,u} = \lambda_g \circ \alpha_{\1,g^{-1}h, g^{-1}u}$ 
we see that for every extensions $\sigma \: (\Delta_2)_{(1)} \to G$ of
$\alpha_{\1,g^{-1}h,g^{-1}u}$ the map $\lambda_g \circ \sigma$ is an
extension of $\alpha_{g,h,u}$. In view of $\lambda_g^*\omega^{\rm eq} =
\rho_\a(g) \circ \omega^{\rm eq}$, we obtain 
$$ \int_{\SS^2} (\lambda_g \circ \sigma)^*\omega^{\rm eq} 
= \int_{\SS^2} \sigma^* \lambda_g^*\omega^{\rm eq} 
= \rho_\a(g).\int_{\SS^2} \sigma^* \omega^{\rm eq}, $$ 
and therefore 
$$F(g,h,u) = \rho_A(g).F(\1, g^{-1}h, g^{-1}u). \leqno(5.1) $$

Let $\Delta_3 \subeq \R^3$ be the standard $3$-simplex. Then we define a piecewise 
smooth map $\gamma$ of its $1$-skeleton to $G$ by 
$$ \gamma(t,0,0) = \alpha_{\1,g}(t), \quad 
\gamma(0,t,0) = \alpha_{\1,gh}(t), \quad 
\gamma(0,0,t) = \alpha_{\1,ghu}(t) $$
and 
$$ \gamma(1-t,t,0) = \alpha_{g,gh}(t), \quad 
\gamma(0,1-t,t) = \alpha_{gh,ghu}(t), \quad 
\gamma(1-t,0,t) = \alpha_{g,ghu}(t). $$
As $G$ is simply connected, we obtain with [Ne02, Prop.~5.6]   
for each face $\Delta_3^j$, $j =0,\ldots, 3$, of $\Delta_3$ a piecewise smooth map
$\gamma_j$ of the first barycentric subdivision 
to $G$, extending the given map on the {$1$-skeleton}. 
These maps combine to a piecewise smooth map 
$\gamma \: (\partial \Delta_3)_{(1)} \to G$. Modulo the period group $\Pi_\omega$ we now have 
$$ \eqalign{ \int_\gamma \omega^{\rm eq} 
&= \int_{\partial \Delta_3} \gamma^*\omega^{\rm eq} 
=  \sum_{i=0}^3 \int_{\gamma_i} \omega^{\rm eq} \cr
&= F(g,gh, ghu) - F(\1, gh,ghu) + F(\1,g,ghu) - F(\1,g,gh) \cr
&= \rho_A(g).f(h,u) - f(gh,u) + f(g,hu) - f(g,h). \cr} $$
Since  $\int_\gamma \omega^{\rm eq} \in \Pi_\omega$, this proves that $f$ is a group cocycle. 
\qed

In the next lemma we will see that for an appropriate choice of paths from 
$\1$ to group elements close to $\1$ the cocycle $f$ will be smooth in  an identity 
neighborhood. The following lemma is a slight generalization of 
Lemma 6.2 in [Ne02].

\Lemma V.2. Let $U \subeq \g$ be an open convex $0$-neighborhood and 
$\phi \: U \to G$ a chart of $G$ with $\phi(0)= \1$ and $d\phi(0) = \id_\g$. 
We then define the arcs $[\1,\phi(x)]$ by $\alpha_{\phi(x)}(t) := \phi(tx)$. 
Let $V \subeq U$ be an open convex $0$-neighborhood with 
$\phi(V) \phi(V) \subeq \phi(U)$ and define 
$x * y := \phi^{-1}(\phi(x)\phi(y))$ for $x,y \in V$. 
If we define $\sigma_{x,y} := \phi \circ \gamma_{x,y}$
with 
$$ \gamma_{x,y} \: \Delta_2 \to U, \quad (t,s) \mapsto t(x * sy) +
s(x*(1-t)y), $$
then for any closed $2$-form $\Omega \in \Omega^2(G,\z)$, $\z$ a
sequentially complete locally convex space, the function 
$$ f_V \: V \times V \to \z, \quad (x,y) \mapsto \int_{\sigma_{x,y}} \Omega $$
is smooth with $d^2 f_V(0,0)(x,y) = {1\over 2}\Omega_\1(x,y)$ 
{\rm(see the end of Appendix B for the notation)}.

\Proof. First we note that the function $V \times V \to U, (x,y) \mapsto x * y$ is smooth. 
We consider the cycle 
$$ \alpha_{\1, \phi(x), \phi(x)\phi(y)} 
= \alpha_{\1, \phi(x), \phi(x*y)} 
= \alpha_{\1,\phi(x)} + \alpha_{\phi(x), \phi(x*y)} - \alpha_{\1, \phi(x*y)}. $$ 
The arc connecting $x$ to $x*y$ is given by $s \mapsto x * sy$, 
so that we may define $\sigma_{x,y} := \phi \circ \gamma_{x,y}$ with $\gamma_{x,y}$ as above. Then 
$$ f_V \: V \times V \to \z, \quad (x,y) \mapsto \int_{\phi \circ \gamma_{x,y}} \Omega =
\int_{\Delta_2} \gamma_{x,y}^*\phi^*\Omega, $$
and 
$$ f_V(x,y) = \int_{\Delta_2} (\phi^*\Omega)\big(\phi(\gamma_{x,y}(t,s))\big)
\Big({\partial \over \partial t}\gamma_{x,y}(t,s), 
{\partial \over \partial s}\gamma_{x,y}(t,s)\Big)\, dt\, ds \leqno(5.2) $$
implies that $f_V$ is a smooth function in $V \times V$. 

The map $\gamma \: (x,y) \mapsto \gamma_{x,y}$ satisfies 
\par\nin (1) $\gamma_{0,y}(t,s) = sy$ and $\gamma_{x,0}(t,s) = (t+s)x$. 
\par\nin (2) ${\partial \over \partial t} \gamma_{x,y} \wedge {\partial
\over \partial s} \gamma_{x,y} = 0$ for $x = 0$ or $y = 0$. 

\par\nin In particular we obtain $f_V(x,0) = f_V(0,y) = 0$. Therefore the
second order Taylor polynomial 
$$ T_2(f_V)(x,y)= f_V(0,0) + df_V(0,0)(x,0) +  df_V(0,0)(0,y) 
+ {1\over 2} d^{[2]} f_V(0,0)\big((x,y),(x,y)\big) $$
of $f_V$ in $(0,0)$ is bilinear and given by 
$$ T_2(f_V)(x,y)
= {1\over 2} d^{[2]} f_V(0,0)\big((x,0),(0,y)\big)
 + {1\over 2} d^{[2]} f_V(0,0)\big((0,y),(x,0)\big)
= d^2 f_V(0,0)(x,y) $$
(see the end of Appendix B). 

Next we observe that (1) implies that 
${\partial \over \partial t} \gamma_{x,y}$ and ${\partial \over \partial s} \gamma_{x,y}$ 
vanish in $(0,0)$. Therefore the chain rule for Taylor expansions
and (1) imply that for each pair $(t,s)$ the second order term of 
$$ (\phi^*\Omega)(\gamma_{x,y}(t,s))
\Big({\partial \over \partial t}\gamma_{x,y}(t,s), 
{\partial \over \partial s}\gamma_{x,y}(t,s)\Big) $$ 
is given by 
$$ (\phi^*\Omega)(\gamma_{0,0}(t,s))(x,y) 
= (d\phi(0)^*\Omega_\1)(x,y) = \Omega_\1(x,y), $$
and eventually 
$$ d^2 f_V(0,0)(x,y) = T_2(f_V)(x,y) 
= \int_{\Delta_2} \, dt\, ds \cdot \Omega_\1(x,y) ={1\over 2} \Omega_\1(x,y). 
\qeddis 

\Corollary V.3. Suppose that $\Gamma_A$ is discrete with $\Pi_\omega
\subeq \Gamma_A$ and construct for 
$\omega \in Z^2_c(\g,\a)$ the group cocycle $f \in Z^2(G,A)$ as above
from the closed $2$-form $\omega^{\rm eq} \in \Omega^2(G,\a)$. 
If the paths $\alpha_{\1,g}$ for $g \in \phi(U)$ are chosen as in 
{\rm Lemma V.2}, then $f \in Z^2_s(G,A)$ with $D(f) = \omega$. 

\Proof. In the notation of Lemma V.2 we have for $x,y \in V$ the relation 
$$ f(\phi(x), \phi(y)) = q_A(f_V(x,y)), $$
so that $f$ is smooth on $\phi(V) \times \phi(V)$, 
and further 
$$ Df(x,y) = d^2 f_V(\1,\1)(x,y) 
- d^2 f_V(\1,\1)(y,x) = \omega(x,y). 
\qeddis  

The outcome of this section is the following result: 

\Theorem V.4. Let $G$ be a connected simply connected Lie group and 
$A$ a smooth $G$-module of the form $\a/\Gamma_A$, where $\Gamma_A
\subeq \a$ 
is a discrete subgroup. Let $\omega \in Z^2_c(\g,\a)$ be a continuous $2$-cocycle and 
$\Pi_\omega \subeq \a^G$ its period group. Then the following
assertions are equivalent: 
\litem{(1)} The Lie algebra extension 
$\a \into \hat\g := \a \oplus_\omega \g \onto \g$
can be integrated to a Lie group extension 
$A \into \hat G \onto G$. 
\litem{(2)} $[\omega] \in \im(D)$. 
\litem{(3)} $\omega \in \im(D)$. 
\litem{(4)} $\Pi_\omega \subeq \Gamma_A$. 
\litem{(5)} If $q_A \: \a \to A$ is the quotient map and 
$P([\omega]) := q_A \circ \per_\omega \: \pi_2(G) \to A,$
then $P([\omega]) = 0$. 

\Proof. (1) $\Rarrow$ (2): If $\hat G$ is an extension of $G$
by $A$ corresponding to the Lie algebra extension 
$\hat\g = \a \oplus_\omega \g$, then we can write $\hat G$ as 
$A \times_f G$ (Proposition II.6), and Lemma II.7 implies that 
$D[f] = [Df] = [\omega]$. 

\par\nin (2) $\Rarrow$ (3): If $[\omega] = D[f] = [Df]$ for some $f
\in Z^2_s(G,A)$, then $Df - \omega \in B^2_c(\g,\a)$ and there exists
an $\alpha \in C^1_c(\g,\a)$ with $Df - \omega = d_\g\alpha$. 
Then the $2$-form $(d_\g  \alpha)^{\rm eq} = d \alpha^{\rm eq} \in
\Omega^2(G,\a)$ is exact (Lemma~B.5), so that its period group is
trivial, and Corollary V.3 implies the existence of $h \in
Z^2_s(G,A)$ with $D h = d_\g \alpha$. Then 
$f_1 := f - h \in Z^2_s(G,A)$ satisfies $D(f -h) = Df - Dh = \omega.$

\par\nin (3) $\Rarrow$ (1): If $D f = \omega$, then 
the Lie group extension 
$A \times_{f} G \to G$ (Proposition II.6) corresponds to the Lie
algebra extension $\a \oplus_{Df} \g = \a \oplus_\omega \g\to \g$
(Lemma II.7). 

\par\nin (1) $\Rarrow$ (4) follows from Proposition IV.5 which implies that 
if $\hat G$ exists, then the period map coincides up to sign with the connecting homomorphism 
$\delta \: \pi_2(G) \to \pi_1(A) \cong \Gamma_A \subeq \a$ in the long
exact homotopy sequence of the principal $A$-bundle $\hat G$. 

\par\nin (4) $\Rarrow$ (3) follows from Corollary V.3. 

\par\nin (4) $\Leftrightarrow$ (5) is a trivial consequence of the
definitions. 
\qed

\sectionheadline{VI. Abelian extensions of non-simply connected groups} 

\nin We have seen in the preceding section that for a simply connected
Lie group $G$ and a smooth $G$-module of the form $A = \a/\Gamma_A$ the image of
the map 
$D \: H^2_{s}(G,A) \to H^2_{c}(\g,\a)$
can be represented by those cocycles $\omega \in Z^2_c(\g,\a)$ for which 
$\Pi_\omega \subeq \Gamma_A \cong \pi_1(A)$. 

In this section we drop the assumption that $G$ is simply connected. 
We write $q_G \: \tilde G \to G$ for the simply connected covering
group of $G$ and identify $\pi_1(G)$ with the discrete central
subgroup $\ker q_G$ of $\tilde G$. 

Let $\omega \in Z^2_{c}(\g,\a)$. In the following we write 
$\rho_A$ for the action of $G$ on $A$, 
$\rho_\a$ for the action of $G$ on $\a$ and $\dot\rho_\a$ for the derived
representation of $\g$ on $\a$. 

\Remark VI.1. (a) To a $2$-cocycle $\omega \in Z^2_c(\g,\a)$ we associate the linear map 
$$ \tilde\omega \: \g \to C^1_c(\g,\a) = \Lin(\g,\a), \quad 
x \mapsto i_x \omega. $$
We consider $\Lin(\g,\a)$ as a $\g$-module with respect to the action 
$$ (x.\alpha)(y) := \dot\rho_\a(x).\alpha(y) - \alpha([x,y]). $$
We do not consider any topology on
this space of maps. 
The corresponding Lie algebra differential 
$d_\g \: C^1(\g,\Lin(\g,\a)) \to C^2(\g,\Lin(\g,\a))$ then satisfies 
$$ \eqalign{ (d_\g \tilde\omega)(x,y)(z) 
&= (x.i_y\omega - y.i_x \omega - i_{[x,y]}\omega)(z) \cr
&= x.\omega(y,z) - \omega(y,[x,z]) 
-  y.\omega(x,z) + \omega(x,[y,z]) - \omega([x,y],z) \cr
&= -z.\omega(x,y) = -d_\g(\omega(x,y))(z).\cr} $$
Since the subspace $B^1_c(\g,\a) = d_\g \a \subeq C^1_c(\g,\a)$ is 
$\g$-invariant, we can also form the quotient $\g$-module 
$$ \hat H^1_c(\g,\a) := C^1_c(\g,\a)/B^1_c(\g,\a). $$
We then obtain a linear map 
$$ f_\omega \: \g \to \hat H^1_c(\g,\a), \quad 
x \mapsto [i_x\omega], $$
and the preceding calculation shows that this map is a $1$-cocycle. 
We call $f_\omega$ the {\it infinitesimal flux cocycle}. In the following
we are concerned with integrating this cocycle to a group cocycle 
$$ F_\omega \: \tilde G \to \hat H^1_c(\g,\a). $$
This is problematic because the right hand side does not have a
natural topology, so that we cannot directly apply Proposition~III.4. 

\par\nin (b) A first step to globalize the situation is to translate
matters from the Lie algebra to vector fields on $G$. We shall see
that on the level of vector fields the infinitesimal flux cocycle corresponds
to the map 
$$  \g \to \hat H^1_{\rm dR}(G, \a) := \Omega^1(G,\a)/dC^\infty(G,\a),
\quad x \mapsto [i_{x_r} \omega^{\rm eq}],$$ 
where for $x \in \g$ we write $x_r$ for the corresponding right invariant
vector field on $G$ with $x_r(\1) = x$. 
Note that for $v \in T_g(G)$ we have 
$$ (i_{x_r}\omega^{\rm eq})(v) = g.\omega(\Ad(g)^{-1}.x,g^{-1}.v). $$

Formally the linear map $\tilde f_\omega \: \g \to \Lin(\g,\a), x
\mapsto i_x \omega$
defines an equivariant $\Lin(\g,\a)$-valued $1$-form 
$\tilde f_\omega^{\rm eq}$ on $G$ as
follows. For each $x \in \g$ evaluation in $x$ is a linear map 
$\ev_x \: \Lin(\g,\a) \to \a, \alpha \mapsto \alpha(x)$ 
and $\ev_x \circ f_\omega \: \g \to \a$
is a continuous linear map, hence defines an equivariant $\a$-valued
$1$-form $(\ev_x \circ \tilde f_\omega)^{\rm eq}$ 
on $G$. For any piecewise smooth path $\gamma \: [0,1] \to G$
we then have 
$$\eqalign{ \int_\gamma (\ev_x \circ \tilde f_\omega)^{\rm eq} 
&= \int_0^1 \gamma(t).\Big(\tilde
f_\omega(\gamma(t)^{-1}\gamma'(t))\Big)(x)\, dt
= \int_0^1 \gamma(t).\omega(\gamma(t)^{-1}\gamma'(t),
\Ad(\gamma(t))^{-1}.x)\, dt \cr
&= \int_0^1 \gamma(t).\omega(\gamma(t)^{-1}\gamma'(t),\gamma(t)^{-1}.
x_r(\gamma(t)))\, dt 
= -\int_\gamma i_{x_r}.\omega^{\rm eq}. \cr} $$

Next we derive some formulas that will be useful in the following. 
We recall the Lie derivative ${\cal L}_{x_r} = d \circ i_{x_r} +
i_{x_r} \circ d$ as an operator on differential forms. 
The equivariance of $\omega^{\rm eq}$ leads to 
$$ {\cal L}_{x_r}.\omega^{\rm eq} = \dot\rho_\a(x) \circ \omega^{\rm eq} $$
([Ne02, Lemma A.2.4]). In view of the closedness of $\omega^{\rm eq}$, this leads to 
$$ d(i_{x_r}\omega^{\rm eq}) = {\cal L}_{x_r}.\omega^{\rm eq} - i_{x_r}
d\omega^{\rm eq}  = \dot\rho_\a(x)\circ \omega^{\rm eq}. \leqno(6.1) $$
Further the formula 
$[{\cal L}_{x_r}, i_{y_r}] = i_{[x_r, y_r]} = - i_{[x,y]_r}$
implies 
$$ \eqalign{ i_{[x,y]_r}\omega^{\rm eq} 
&= i_{y_r} {\cal L}_{x_r} \omega^{\rm eq} - {\cal L}_{x_r}i_{y_r}
\omega^{\rm eq} 
= i_{y_r}(\dot\rho_\a(x)\circ \omega^{\rm eq}) - (i_{x_r} \circ d 
+ d \circ i_{x_r})i_{y_r} \omega^{\rm eq} \cr
&= \dot\rho_\a(x) \circ i_{y_r}\omega^{\rm eq} - \dot\rho_\a(y) \circ
i_{x_r}\omega^{\rm eq} - d(i_{x_r} i_{y_r}\omega^{\rm eq}). \cr} $$
This means that the $\a$-valued $1$-form 
$$ \dot\rho_\a(x) \circ i_{y_r}\omega^{\rm eq} 
-  \dot\rho_\a(y) \circ i_{x_r}\omega^{\rm eq} 
- i_{[x,y]_r}\omega^{\rm eq} = d(i_{x_r}i_{y_r} \omega^{\rm eq}) \leqno(6.2) $$
is exact, which entails that 
$$ \tilde f_\omega \: \g \to \hat H^1_{\rm dR}(G,\a), \quad 
x \mapsto [i_{x_r}\omega^{\rm eq}] $$
is a $1$-cocycle with respect to the representation of $\g$ on 
$\hat H^1_{\rm dR}(G,\a)$ given by 
$x.[\alpha] := [\dot\rho_\a(x) \circ \alpha]$. Since the form $\omega^{\rm eq}$ is closed, 
the map $\tilde f_\omega$ also is a cocycle with respect to the action given by 
$x.[\alpha] := [-{\cal L}_{x_r}.\alpha]$ 
because $\g \to {\cal V}(G), x \mapsto - x_r$ 
is a homomorphism of Lie algebras (cf.\ Lemma IX.8). 
\qed

\Lemma VI.2. Let $\gamma \: [0,1]  \to G$ be a piecewise smooth
path. Then we obtain a continuous linear map 
$$ \tilde F_\omega(\gamma) \: \g \to \a, \quad x \mapsto 
-\int_\gamma i_{x_r}\omega^{\rm eq}
= \int_0^1
\gamma(t).\omega(\gamma(t)^{-1}\gamma'(t),\Ad(\gamma(t))^{-1}.x)\,
dt $$
with the following properties: 
\litem{(1)} If $\gamma(1)^{-1} \gamma(0)$ is contained in $Z(G)$ and
acts trivially on $\a$, then $\tilde F_\omega(\gamma) \in
Z^1_c(\g,\a)$. 
\litem{(2)} If $\gamma_1$ and $\gamma_2$ are homotopic with fixed
endpoints, then $\tilde F_\omega(\gamma_1) - \tilde F_\omega(\gamma_2)$ is a coboundary. 
\litem{(3)} For a piecewise smooth curve $\eta \: [0,1] \to G$ we have 
$$ \int_\eta \tilde F_\omega(\gamma)^{\rm eq} = \int_H \omega^{\rm eq} $$
for the piecewise smooth map 
$H \: [0,1]^2 \to G, (t,s) \mapsto \eta(s)\cdot \gamma(t).$

\Proof. In view of formula (6.2) above, we find for 
$x,y \in \g$ the relation 
$$ \eqalign{ d_\g(\tilde F_\omega(\gamma))(x,y) 
&= x.\tilde F_\omega(\gamma)(y) - y.\tilde F_\omega(\gamma)(x)
- \tilde F_\omega(\gamma)([x,y]) \cr
&= -\int_\gamma \dot\rho_\a(x) \circ i_{y_r}\omega^{\rm eq} 
-  \dot\rho_\a(y) \circ i_{x_r}\omega^{\rm eq} 
- i_{[x,y]_r}\omega^{\rm eq} 
= -\int_\gamma  d(i_{x_r} i_{y_r}\omega^{\rm eq})  \cr
&= \omega^{\rm eq}(\gamma(0))\big(y_r(\gamma(0)), x_r(\gamma(0))\big) 
-\omega^{\rm eq}(\gamma(1))\big(y_r(\gamma(1)), x_r(\gamma(1))\big) \cr
&= \gamma(0).\omega(\Ad(\gamma(0))^{-1}.y,\Ad(\gamma(0))^{-1}.x)
- \gamma(1).\omega(\Ad(\gamma(1))^{-1}.y, \Ad(\gamma(1))^{-1}.x). \cr} $$

\par\nin (1) If $\gamma(1)^{-1}\gamma(0) \in Z(G) = \ker \Ad$ 
acts trivially on $\a$, then the above formula implies that 
$d_\g\big(\tilde F_\omega(\gamma)\big) = 0$, i.e., that $\tilde
F_\omega(\gamma) \in Z^1_c(\g,\a)$. 

\par\nin (2) For $g \in G$ we first observe that 
$$ \eqalign{ \tilde F_\omega(g\cdot\gamma)(x) 
&= -\int_{\lambda_g \circ \gamma} i_{x_r}.\omega^{\rm eq}
= \int_0^1 g\gamma(t).\omega(\gamma(t)^{-1}.\gamma'(t), \Ad(g\gamma(t))^{-1}.x)\, dt\cr
&= g.\int_0^1 \gamma(t).\omega(\gamma(t)^{-1}.\gamma'(t),\Ad(\gamma(t))^{-1}\Ad(g)^{-1}.x)\, dt\cr
&= g.\tilde F_\omega(\gamma)(\Ad(g)^{-1}.x) = (g.\tilde F_\omega(\gamma))(x). \cr} $$
For the natural action of $G$ on $\Lin(\g,\a)$ by 
$(g.\phi)(x) := g.\phi(\Ad(g)^{-1}.x)$ and the left translation action
on the space $C^1_{pw}(I,G)$ of piecewise smooth maps $I :=[0,1] \to G$,
the preceding calculation shows that the map 
$$ \tilde F_\omega \: C^1_{pw}(I,G) \to \Lin(\g,\a) = C^1_c(\g,\a) $$
is equivariant. 

For the composition 
$$ (\gamma_1 \sharp \gamma_2)(t) := \cases{ 
\gamma_1(2t) & for $0 \leq t \leq {1\over 2}$ \cr 
\gamma_1(1)\gamma_2(0)^{-1}\gamma_2(2t-1) & for ${1\over 2} \leq t \leq 1$ \cr} $$
 of paths we thus obtain the composition  formula 
$$ \tilde F_\omega(\gamma_1 \sharp \gamma_2) 
= \tilde F_\omega(\gamma_1) + \tilde
F_\omega(\gamma_1(1)\gamma_2(0)^{-1}\gamma_2)
= \tilde F_\omega(\gamma_1)
+ \gamma_1(1)\gamma_2(0)^{-1}.\tilde F_\omega(\gamma_2). \leqno(6.3) $$
For the inverse path $\gamma^-(t) := \gamma(1-t)$ we trivially get
$\tilde F_\omega(\gamma^-) = -\tilde F_\omega(\gamma)$
from the transformation formula for one-dimensional integrals. 
If the two paths $\gamma_1$ and $\gamma_2$ have the same start and endpoints,
then the path $\gamma_1 \sharp \gamma_2^-$ is closed, and we derive
with (1) that 
$$ \tilde F_\omega(\gamma_1) - \tilde F_\omega(\gamma_2)
= \tilde F_\omega(\gamma_1) + \gamma_1(1)\gamma_2^-(0)^{-1}.\tilde F_\omega(\gamma_2^-)
= \tilde F_\omega(\gamma_1 \sharp \gamma_2^-) \in Z^1_c(\g,\a). $$

That two paths $\gamma_1$ and $\gamma_2$ with the same endpoints are homotopic with fixed endpoints
implies that the loop $\gamma := \gamma_1 \sharp \gamma_2^-$ is
contractible. It therefore has a closed piecewise smooth 
lift $\tilde \gamma \: [0,1] 
\cong \partial \Delta_2\to \tilde G$ with $q_G \circ \tilde\gamma =
\gamma$. Using Proposition 4.6 in [Ne02], we find a piecewise smooth map 
$\tilde\sigma \: \Delta_2 \to \tilde G$ 
such that $\tilde\sigma\res_{\partial\Delta_2} = \tilde\gamma$. 
Let $\sigma := q_G \circ \tilde\sigma$. Then $\sigma\res_{\partial
\Delta_2} = \gamma$, so that Stoke's Theorem and
formula (6.1) lead to 
$$ \eqalign{ -\tilde F_\omega(\gamma)(x) 
&= \int_\gamma i_{x_r}\omega^{\rm eq} 
= \int_{\partial\Delta_2} \sigma^*(i_{x_r}\omega^{\rm eq})
= \int_{\Delta_2} \sigma^*d (i_{x_r}\omega^{\rm eq})\cr
&= \int_\sigma d(i_{x_r} \omega^{\rm eq})
= \int_\sigma \dot\rho_\a(x) \circ \omega^{\rm eq}
= \dot\rho_\a(x).\int_\sigma \omega^{\rm eq}. \cr} $$
Therefore $\tilde F_\omega(\gamma) \in B^1_c(\g,\a)$, and (2)
follows. 

\par\nin (3) We have 
$$ \eqalign{ 
\int_\eta \tilde F_\omega(\gamma)^{\rm eq} 
&= \int_0^1 \eta(s).\tilde F_\omega(\gamma)(\eta(s)^{-1}.\eta'(s))\, ds \cr
&= \int_0^1 \int_0^1
\eta(s)\gamma(t).\omega(\gamma(t)^{-1}.\gamma'(t), 
\Ad(\gamma(t)^{-1})\circ \eta(s)^{-1}.\eta'(s))\, dt\, ds \cr
&= \int_0^1 \int_0^1
H(t,s).\omega(H(t,s)^{-1}\eta(s).\gamma'(t), 
H(t,s)^{-1}.(\eta'(s).\gamma(t)))\, dt\, ds \cr
&= \int_0^1 \int_0^1
H(t,s).\omega\Big(H(t,s)^{-1}.{\partial H(t,s) \over \partial t}, 
H(t,s)^{-1}.{\partial H(t,s) \over \partial s}\Big)\, dt\, ds \cr
&= \int_{[0,1]^2} H^* \omega^{\rm eq} = \int_H \omega^{\rm eq}. \cr}
$$
\qed

\Proposition VI.3. We have a well-defined map 
$$ F_\omega \: \tilde G \to \hat H^1_c(\g,\a) = \Lin(\g,\a)/B^1_c(\g,\a), \quad 
g \mapsto [\tilde F_\omega(q_G \circ \gamma_g)] := \tilde F_\omega(q_G
\circ \gamma_g) + B^1_c(\g,\a), $$
where $\gamma_g \: [0,1] \to \tilde G$ is piecewise smooth with $\gamma_g(0) =
\1$ and $\gamma_g(1) = g$. The map $F_\omega$ is a $1$-cocycle
with respect to the natural action of $\tilde G$ on $\hat
H^1_c(\g,\a)$. 
Moreover, we obtain by restriction a group homomorphism 
$Z(\tilde G) \cap \ker \rho_\a \to H^1_{c}(\g,\a), 
[\gamma] \mapsto [\tilde F_\omega(\gamma)]$
and further by restriction to $\pi_1(G)$ a homomorphism 
$$ F_\omega \: \pi_1(G) \to H^1_{c}(\g,\a). $$

\Proof. That $F_\omega$ is well-defined follows from 
Lemma VI.1(2) because two different choices of paths $\gamma_g$ and
$\eta_g$ lead to paths $q_G \circ \gamma_g$ and 
$q_G \circ \eta_g$ in $G$ which are homotopic with fixed endpoints. 
Next we note that for paths $\gamma_{g_i}$, $i=1,2$, from $\1$ to
$g_i$ in $\tilde G$ the composed path 
$\gamma_{g_1} \sharp \gamma_{g_2}$ connects $\1$ to $g_1 g_2$. Hence 
the composition formula (6.3) leads to 
$$ F_\omega(g_1 g_2) 
= \tilde F_\omega(\gamma_{g_1} \sharp \gamma_{g_2}) 
= \tilde F_\omega(\gamma_{g_1}) + g_1.\tilde F_\omega(\gamma_{g_2}) 
= F_\omega(g_1) + g_1.F_\omega(g_2), $$
showing that the map
$F_\omega$ is a $1$-cocycle. 

Since $Z(\tilde G) \cap \ker \rho_\a$ acts trivially on $\g$ and $\a$, hence on
$\Lin(\g,\a)$, the restriction of $F_\omega$ to this subgroup is 
a group homomorphism, and Lemma VI.2(1) shows that its values lie 
in the subspace $H^1_c(\g,\a)$. 
\qed

The cocycle $F_\omega \: \tilde G \to \hat H^1_c(\g,\a)$ is called the
{\it flux cocycle} and its restriction to $\pi_1(G)$ the {\it flux
homomorphism} for reasons that will become clear in Definition~IX.9 below. 
Next we relate the flux homomorphism to group extensions. 
Although the following proposition is quite technical, it contains a lot of interesting 
information, even for the case of non-connected groups $A$. 

\Proposition VI.4. Let $A$ be an abelian Lie group whose identity component 
satisfies $A_0 \cong \a/\Gamma_A$, where $\Gamma_A \subeq \a$ is a discrete subgroup. 
Further let $q \: \hat G \to G$ be a Lie group extension of $G$ by $A$
corresponding to the Lie algebra cocycle $\omega \in Z^2_c(\g,\a)$, so
that its Lie algebra is $\hat\g \cong \a \oplus_\omega \g$. In these
terms we write the adjoint action of $\hat G$ on $\hat\g$ as 
$$ \Ad(g).(a,x) = (g.a - \theta(g)(g.x), g.x), \quad g \in \hat G, a
\in \a, x \in \g, \leqno(6.4) $$
where $g.x = \Ad(q(g)).x$ and 
$$ \theta \: \hat G \to C^1_c(\g,\a) = \Lin(\g,\a) $$
is a $1$-cocycle with respect to the action of $\hat G$ on
$\Lin(\g,\a)$ by $(g.\alpha)(x) := g.\alpha(g^{-1}.x)$. Its
restriction $\theta_A := \theta\res_A$ is a homomorphism given by 
$$ \theta_A(a) = D(d_G(a)) \quad \hbox{ with } \quad 
(d_G a)(g) := g.a - a\ \hbox{ and }\ D(d_G a)(x) := x.a 
:= \big(d(d_G a)(\1)\big)(x). $$
This
$1$-cocycle maps $A_0$ to $B^1_c(\g,\a)$ and factors through a
$1$-cocycle 
$$ \oline\theta \: \hat G/A_0 \to \hat H^1_c(\g,\a) = \Lin(\g,\a)/B^1_c(\g,\a),
\quad q(g) \mapsto [\theta(g)]. $$
The map $\oline q \: \hat G/A_0 \to G, gA_0 \mapsto q(g)$ is a covering of
$G$, so that there is a unique
covering morphism $\hat q_G \: \tilde G \to \hat G/A_0$ with 
$\oline q \circ \hat q_G = q_G$, and the following assertions hold: 
\litem{(1)} The coadjoint action of $\hat G$ on $\hat\g$ and the flux cocycle
are related by 
$$ F_\omega = -\oline\theta \circ \hat q_G \: \tilde G \to \hat
H^1_c(\g,\a). $$
\litem{(2)} If $\delta \: \pi_1(G) \to \pi_0(A) \subeq \hat G/A_0$ 
is the connecting homomorphism 
from the long exact homotopy sequence of the principal $A$-bundle $q
\: \hat G \to G$, then 
$$ F_\omega = -\oline\theta_A \circ \delta \: \pi_1(G) \to H^1_c(\g,\a),$$ 
where $\theta_A \: \pi_0(A) \to H^1_c(\g,\a)$ 
is the characteristic homomorphism of the 
smooth $G$-module~$A$. 
\litem{(3)} If $A$ is connected, then $F_\omega(\pi_1(G)) = \{0\}$. 

\Proof. From the description of the Lie algebra $\hat\g$ as 
$\a \oplus_\omega \g$, it is clear that there exists a function 
$\theta \: \hat G \to \Lin(\g,\a)$ for which the map 
$(g,x) \mapsto \theta(g)(x)$ is smooth and the adjoint action of $\hat
G$ on $\g$ is given by (6.4). Since $\Ad$ is a representation of $G$, 
we have $\theta(\1,x) =0$ and 
$$ \theta(g_1 g_2)(g_1 g_2 x) = g_1.\theta(g_2)(g_2.x) + \theta(g_1)(g_1g_2.x), 
\quad g_1, g_2\in \hat G, x\in \g, \leqno(6.5) $$ 
which means that 
$$ \theta(g_1 g_2) = g_1.\theta(g_2) + \theta(g_1), $$
i.e., $\theta$ is a $1$-cocycle. As $A$ acts trivially on $\a$ and
$\g$, the restriction $\theta_A := \theta\res_A$ is a homomorphism 
$$ \theta_A \: A \to Z^1_c(\g,\a) \quad \hbox{ with } \quad 
\Ad(b).(a,x) = (a - \theta_A(b)(x), x), \quad b \in A, a \in \a, x
\in \g. $$
The relation $\theta(b) \in Z^1_c(\g,\a)$ follows directly from $\Ad(b) \in \Aut(\hat\g)$. 

For $\hat g \in \hat G$ with $q(\hat g) = g$ and $b \in A$ we have 
$$ b \hat g b^{-1} 
= (b \hat g b^{-1} \hat g^{-1}) \hat g 
= (b - g.b) \cdot \hat g, $$
which leads to 
$$ \Ad(b).(a,x) = (a - x.b,x) $$
and therefore to $\theta_A(b)(x) = x.b$. 
For $a \in \a$ and $b = q_A(a)$ we have $x.b = x.a$, so that 
$\theta(A_0) = B^1_c(\g,\a)$. Hence $\theta$ factors through
a $1$-cocycle $\oline\theta \: \hat G/A_0 \to \hat H^1_c(\g,\a)$ whose
restriction $\oline\theta_A$ to $\pi_0(A) = A/A_0$ is given by 
$$ \oline\theta_A \: \pi_0(A) \cong A/A_0 \to H^1_c(\g,\a), \quad [a] \mapsto
[\theta_A(a)] = [D(d_G a)]. $$

\par\nin (1) For a fixed $x \in \g$ the cocycle condition (6.5) 
implies for the smooth functions $\theta_x  \:\hat G \to \a, g \mapsto 
\theta(g)(x)$ the relation  
$$ \theta_x(gh)  = g.\theta_{g^{-1}.x}(h) + \theta_x(g). $$
For the differentials we thus obtain 
$$ d\theta_x(g) d\lambda_g(\1) = \rho_\a(g) \circ d \theta_{g^{-1}.x}(\1). \leqno(6.6) $$
From formula (6.4) for the adjoint action, we get in view of $\theta(\1) = 0$ the formula 
$$ (x'.a - x.a' + \omega(x',x), [x',x]) = 
\ad(a',x')(a,x) = (x'.a - d\theta_x(\1)(a',x'), [x',x]), $$
so that $\theta$ and the corresponding Lie algebra cocycle are related by 
$$ d\theta_{x}(\1)(a',x') = \omega(x,x') +  x.a'. $$
With (6.6) this further leads to 
$$ d\theta_x(g) d\lambda_g(\1)(a',x') 
= g.\big(\omega(g^{-1}.x, x') + (g^{-1}.x).a'\big) 
= \omega^{\rm eq}(x_r(q(g)), d\lambda_{q(g)}(\1).x') + x.(g.a'). $$
In $\Omega^1(\hat G,\a)$ we therefore have the relation 
$$ d\theta_x =\dot\rho_\a(x) \circ p_\a^{\rm eq} + q^*(i_{x_r} \omega^{\rm eq}), $$
where $p_\a(a',x') = a'$ is the projection of $\hat\g$ onto $\a$ and $p_\a^{\rm
eq}$ the corresponding equivariant $1$-form on~$\hat G$. 

Let $\gamma \: [0,1] \to G$ be any piecewise smooth loop based in $\1$. Then there exists a piecewise 
smooth map $\hat\gamma \: [0,1] \to \hat G$ with $q \circ \hat\gamma = \gamma$ and 
$\hat\gamma(0) = \1$. 
Then 
$\tilde\gamma := \hat q_G \circ \hat\gamma \: [0,1] \to \tilde G$ is
the unique lift of $\gamma$ to a piecewise smooth path in $\tilde G$
starting in $\1$. We now have 
$$ \eqalign{ 
-\tilde F_\omega(\gamma)(x) 
&=  \int_\gamma i_{x_r} \omega^{\rm eq} 
=  \int_{[0,1]} \gamma^*(i_{x_r} \omega^{\rm eq})
=  \int_{[0,1]} \hat\gamma^*q^*(i_{x_r} \omega^{\rm eq})\cr
&=  \int_{\hat\gamma} q^*(i_{x_r} \omega^{\rm eq}) 
= \int_{\hat\gamma} d\theta_x - \rho_\a(x) \circ p_\a^{\rm eq} \cr
&=  \theta_x(\hat\gamma(1)) - \theta_x(\hat\gamma(0)) 
- \rho_\a(x) .\int_{\hat \gamma} p_\a^{\rm eq} 
=  \theta(\hat\gamma(1))(x) 
- \rho_\a(x) .\int_{\hat \gamma} p_\a^{\rm eq}. \cr} $$ 
This means that 
$$ F_\omega(\tilde\gamma(1)) = [\tilde F_\omega(\gamma)]
= -[\theta(\hat\gamma(1))] 
= -\oline\theta(\hat q_G(\tilde\gamma(1))) $$
and therefore that $F_\omega = -\oline\theta \circ \hat q_G$ because 
$\gamma$ was arbitrary. 

\par\nin (2) If $\gamma \: [0,1] \to G$ is a piecewise smooth loop
based in $\1$, then $\hat\gamma(1) \in \ker q = A$ and 
$\delta([\gamma]) = [\hat\gamma(1)]$, as an element of $\pi_0(A)$.
This means that $\delta$ can be considered as the restriction of 
$\hat q_G \: \tilde G \to \hat G/A_0$ to the subgroup $\pi_1(G) = \ker
q_G$. Therefore (2) follows from (1) by restriction. 

\par\nin (3) If $A$ is connected, then $\delta = 0$, so that (3)
follows from (2). 
\qed 

\Corollary VI.5. If, in addition to the assumptions of {\rm
Proposition VI.4}, the group $G$ is simply connected, then 
$A$ is connected and 
$$ F_\omega = -\oline\theta \: G  \to \hat H^1_c(\g,\a). $$
On the subgroup 
$A^\sharp := q^{-1}(Z(G) \cap \ker \rho_A)$ of $\hat G$ 
the cocycle $\theta$ restricts to a homomorphism 
$$ \theta^\sharp \: A^\sharp \to Z^1_c(\g,\a), \quad a \mapsto D(d_G(a)), 
\leqno(6.7) $$
where for each $a \in A^\sharp$ the smooth cocycle 
$d_G(a) \in Z^1_s(G,A)$ is defined by 
$d_G(a)(q(g)) := gag^{-1}a^{-1}$. 
For two piecewise smooth curves 
$\gamma, \eta \: [0,1] \to G$ with 
$\gamma(0) = \eta(0) = \1$ and $\gamma(1), \eta(1) \in A^\sharp$ we
have for $H \: I^2 \to G, H(t,s) = \gamma(t)\eta(s)$ the formula 
$$ \gamma(1) \eta(1) \gamma(1)^{-1} \eta(1)^{-1} 
= -\int_\gamma \tilde F_\omega(\eta)^{\rm eq} + \Gamma_A 
= \int_H \omega^{\rm eq} + \Gamma_A. \leqno(6.8) $$

\Proof. To derive the first part from Propositions~VI.3 and VI.4, 
we only have to observe 
that for $a \in A^\sharp$ the condition $\rho_A(a) = \id_A$ implies that 
$d_G(a)$ is well-defined on $G$ by $d_G(a)(q(g)) = gag^{-1}a^{-1}$, 
and that this is an element of $A$ because $q(a) \in Z(G)$ implies 
$d_G(a) \in \ker q$. 

For (6.8) we first observe that for 
$x \in \a$ and $q_A(x) = x + \Gamma_A \in A$ 
the map $d_G q_A(x) \: G \to A$ satisfies 
$$ \eqalign{ 0 
&= \rho_A(\gamma(1))(q_A(x)) - q_A(x) = (d_{G} q_A(x))(\gamma(1)) 
= \int_\gamma d(d_{G}(q_A(x))) + \Gamma_A  \cr
&= \int_\gamma (D(d_{G}q_A(x)))^{\rm eq} + \Gamma_A  
= \int_\gamma (d_\g x)^{\rm eq} + \Gamma_A, \cr} $$ 
so that the integration along $\gamma$ yields a well-defined map 
$\hat H^1_c(\g,\a)  \to \a, [\alpha] \mapsto \int_\gamma
\alpha^{\rm eq}.$
We therefore get with Proposition~VI.4, Lemma VI.2(3) (note the sign change) 
and $-\oline\theta = F_\omega$: 
$$ \eqalign{ 
&\ \ \ \ \gamma(1) \eta(1) \gamma(1)^{-1} \eta(1)^{-1} 
= d_G(\eta(1))(\gamma(1)) 
= \int_\gamma d(d_G(\eta(1))) + \Gamma_A  
= \int_\gamma D(d_G(\eta(1)))^{\rm eq}  + \Gamma_A \cr
&= \int_\gamma \theta^\sharp(\eta(1))^{\rm eq}  + \Gamma_A 
= -\int_\gamma F_\omega(\eta(1))^{\rm eq}  + \Gamma_A 
= -\int_\gamma \tilde F_\omega(\eta)^{\rm eq} + \Gamma_A 
= \int_H \omega^{\rm eq} + \Gamma_A. \cr }$$
\qed

\Corollary VI.6. Suppose that $A \cong \a/\Gamma_A$, that  
$q_G \: \tilde G \to G$ is a universal covering homomorphism, 
let $q \: \hat G \to \tilde G$ be an $A$-extension of 
$\tilde G$ corresponding to $\omega \in Z^2_c(\g,\a)$, and  
${\hat\pi_1(G)} := q^{-1}(\pi_1(G))$. 
Then the following are equivalent: 
\litem{(1)} $F_\omega(\pi_1(G)) = 0$. 
\litem{(2)} $\theta({\hat\pi_1(G)})\subeq B^1_c(\g,\a) = \theta(A)$. 
\litem{(3)} ${\hat\pi_1(G)} = A + \ker(\theta\res_{{\hat\pi_1(G)}})$.  
\litem{(4)} $q(\ker(\theta\res_{{\hat\pi_1(G)}})) = \pi_1(G)$. 
\litem{(5)} There exists a group homomorphism 
$\sigma \: \pi_1(G) \to \ker(\theta\res_{{\hat\pi_1(G)}})
 = {\hat\pi_1(G)} \cap Z(\hat G)$ with 
$q \circ \sigma = \id_{\pi_1(G)}$. 

\Proof. The equivalence of (1) and (2) follows from Corollary~VI.5, 
and (2) is clearly equivalent to (3), 
which in turn is equivalent to (4) because $\ker q= A$. 

That (5) implies (4) is trivial. If (4) is satisfied, 
then we first observe that $\ker (\theta \res_{{\hat\pi_1(G)}}) 
= {\hat\pi_1(G)} \cap Z(\hat G)$,
so that (3) implies that 
${\hat\pi_1(G)}$ is abelian. Further (6.7) in Corollary VI.5 leads to 
$$ \ker (\theta\res_{{\hat\pi_1(G)}}) \cap \ker q = \ker (\theta\res_A) = q_A(\a^\g), $$
which is a divisible group. Hence the extension 
$q_A(\a^\g) \into \ker(\theta\res_{{\hat\pi_1(G)}}) \onto \pi_1(G)$
splits, which is~(5). 
\qed

The following theorem is the central result of the present paper. 
\Theorem VI.7. {\rm(Integrability Criterion)} Let $G$ be a connected
Lie group and $A$ be a smooth $G$-module with 
$A_0  \cong \a/\Gamma_A$, 
where $\Gamma_A$ is a discrete subgroup of the sequentially complete locally convex space
$\a$. For each $\omega \in Z^2_c(\g,\a)$ the abelian Lie algebra extension 
$\a \into \hat\g := \a \times_\omega \g \onto \g$
integrates to a Lie group extension $A \into \hat G \onto G$ with a 
connected Lie group $\hat G$ if and only if 
\litem{(1)} $\Pi_\omega \subeq \Gamma_A$, and 
\litem{(2)} there exists a surjective homomorphism $\gamma \: \pi_1(G) \to \pi_0(A)$ 
such that the flux homomorphism $F_\omega \: \pi_1(G) \to H^1_{c}(\g,\a)$ 
is related to the characteristic homomorphism $\oline\theta_A \: \pi_0(A) 
\to H^1_c(\g,\a)$ by 
$$ F_\omega = \oline\theta_A \circ \gamma. $$ 

If $A$ is connected, then {\rm(2)} is equivalent to $F_\omega = 0$. 

\Proof. Suppose first that a Lie group extension $\hat G$ of $G$ by $A$ exists which corresponds 
to the Lie algebra cocycle $\omega$. According to Proposition IV.5, up to sign 
the period map can be interpreted as 
the connecting map $\pi_2(G) \to \pi_1(A) \cong \Gamma_A$. This 
implies (1). That (2) is satisfied follows from Proposition VI.4(2) because 
in view of the connectedness of $\hat G$, the long exact homotopy sequence 
of the $A$-bundle $\hat G$ implies that the connecting homomorphism 
$\delta \: \pi_1(G) \to \pi_0(A)$ is surjective. 

Conversely, suppose that (1) and (2) hold. 
Let $q_G \: \tilde G \to G$ denote the simply connected covering group of $G$ and recall 
that $\pi_2(q_G)$ is an isomorphism $\pi_2(\tilde G) \to \pi_2(G)$. 
We may therefore identify the period maps 
$\per_\omega$ of $G$ and $\tilde G$ and likewise for all quotients of 
$\tilde G$ by subgroups of $\pi_1(G)$. 

From the case of simply connected groups (Proposition V.3) 
we know that there exists an $A_0$-extension 
$q^\sharp \: G^\sharp \to \tilde G$, where $A$ 
carries the natural $\tilde G$-module structure induced by 
the $G$-module structure. The Lie algebra of $G^\sharp$ is 
$\hat\g = \a \oplus_\omega \g$. Let 
$G_1 := \tilde G/\ker\gamma$ and observe that $\pi_1(G_1) \cong \ker\gamma$. 
Condition (2) implies 
$\pi_1(G_1) = \ker \gamma \subeq \ker F_\omega$, so that 
Corollary VI.6 implies that there exists a homomorphism 
$$ \sigma \: \pi_1(G_1) \to \ker(\theta\res_{{\hat\pi_1(G)}}) \subeq Z(G^\sharp) $$ 
with $q^\sharp \circ \sigma = \id_{\pi_1(G_1)}$. 
Then the image of $\sigma$ is a discrete central subgroup of $G^\sharp$, and therefore 
$$ \hat G := G^\sharp/\sigma(\pi_1(G_1)) $$
defines an abelian extension 
$A_0 \into \hat G \sssmapright{q_1} G_1$
corresponding to the given Lie algebra extension $\a \oplus_\omega \g \to \g$. 
If $q_1 \: G_1 \to G$ is the quotient map with kernel 
$\pi_1(G)/\ker \gamma \cong \im \gamma \cong \pi_0(A)$, then 
$B := q_1^{-1}(\pi_1(G)/\ker \gamma)$ is a subgroup of $\hat G$ 
with $B_0 = A_0$ and $\pi_0(B) =B/B_0 \cong \pi_0(A)$, which implies 
that $B \cong B_0 \times \pi_0(B) \cong A_0 \times \pi_0(A) \cong A$ as 
abelian Lie groups. As $\gamma$ factors through an isomomorphism 
$\oline\gamma \: \pi_0(B) \to \pi_0(A)$ and the characteristic maps 
$\oline\theta_A \: \pi_0(A) \to H^1_c(\g,\a)$ and 
$\oline\theta_B \: \pi_0(B) \to H^1_c(\g,\a)$ satisfy 
$$ \oline\theta_A \circ \oline\gamma = \oline\theta_B $$ 
(Proposition VI.4, Corollary VI.5), Lemma III.7 implies that 
$A \cong B$ as smooth $G$-modules. Therefore 
$\hat G$ is an $A$-extension of $G$. 
\qed

\Remark VI.8. (a) Suppose that only (1) in Theorem VI.7 is satisfied, 
and that $A$ is connected. Consider 
the corresponding extension $q^\sharp \: G^\sharp \to \tilde G$ of 
$\tilde G$ by $A \cong \a/\Gamma_A$. 
Then $G \cong G^\sharp/{\hat\pi_1(G)}$, where 
${\hat\pi_1(G)} := (q^\sharp)^{-1}(\pi_1(G))$ 
is a central $A$-extension of $\pi_1(G)$, hence
$2$-step nilpotent. 

We have seen in the proof of Theorem VI.7 that whenever 
an $A$-extension $\hat G$ of $G$ corresponding to $\omega \in Z^2_c(\g,\a)$ exists, 
then it can be obtained as a quotient of $G^\sharp$ by a subgroup 
$\sigma(\pi_1(G))$, where 
$\sigma \: \pi_1(G) \to Z(G^\sharp) \cap \hat\pi_1(G)$ 
a splitting homomorphism for $\hat\pi_1(G)$. 
This implies in particular that $\hat\pi_1(G)$ is abelian. 

Let us take a closer look at the nilpotent group 
$\hat\pi_1(G)$. If this group is abelian, then the divisibility 
of $A_0 \cong \a/\Gamma_A$ implies that $\hat\pi_1(G)$ splits as an 
$A_0$-extension of $\pi_1(G)$. Clearly this condition is weaker than the 
requirement that it splits by a homomorphism 
$\sigma \: \pi_1(G) \to \hat\pi_1(G) \cap Z(G^\sharp)$. 

That $\hat\pi_1(G)$ is abelian is equivalent to the triviality 
of the induced commutator map 
$$ C \: \pi_1(G) \times \pi_1(G) \to A. $$
According to Corollary VI.5, 
$$ C([\gamma],[\eta]) = -\int_\gamma \tilde F_\omega(\eta)^{\rm eq} + \Gamma_A
= -P(F_\omega([\eta]))([\gamma]) + \Gamma_A, \leqno(6.9) $$
where $P \: H^1_c(\g,\a) \to \Hom(\pi_1(G),\a)$ is defined as in
Proposition~III.4. Therefore the commutator map vanishes if and only
if 
$$ P(F_\omega(\pi_1(G)))(\pi_1(G)) \subeq \Gamma_A. \leqno(6.10) $$
This means that for all smooth loops 
$\gamma, \eta \: \SS^1 \to G$ and $H \: \T^2 \to G, (t,s) \mapsto
\gamma(t)\eta(s)$ we have 
$$ \int_{\T^2} H^*\omega^{\rm eq} = P(F_\omega([\eta]))([\gamma]) \in
\Gamma_A. $$
In view of Proposition III.4, Condition (6.10) is equivalent to 
$$ \im(F_\omega) \subeq \im(D_1) \subeq H^1_c(\g,\a), \leqno(6.11) $$
i.e., that the image of the flux homomorphism consists of classes of 
integrable $1$-cocycles. 

In Corollary VI.5 we have seen that we have a homomorphism 
$$ \theta^\sharp = D_1 \circ d_{\tilde G} \: \hat\pi_1(G) \to Z^1_c(\g,\a) $$
which factors through the (negative) flux homomorphism 
$- F_\omega \: \pi_1(G) \to H^1_c(\g,\a).$
The group $\hat\pi_1(G)$ is a smooth $\tilde G$-module which is 
abelian if and only $\pi_1(G)$ acts trivially, which in turn is (6.11). 
If this is the case, then 
$$ - F_\omega \: \pi_0(\hat\pi_1(G)) \cong \pi_1(G) \to H^1_c(\g,\a)$$
is the characteristic homomorphism of the smooth $G$-module 
$\hat\pi_1(G)$. In view of Lemma III.7, it vanishes if and only if 
the identity component $\hat\pi_1(G)_0 \cong A$ has a $G$-invariant complement. 

In Example IX.17 below 
we will see cases where the commutator map vanishes and the flux
homomorphism $F_\omega \: \pi_1(G) \to H^1_c(\g,\a)$ is non-zero. 

\par\nin (b) With similar arguments as in Section IV, resp.\ Section 5 of [Ne02], 
we can define a {\it toroidal period map} by observing that the
integration map 
$$ \tilde\per^\T_\omega \: C^\infty(\T^2, G) \to \a^G, \quad 
[\sigma] \mapsto \int_\sigma \omega^{\rm eq} $$
is constant on the connected components and defines a map 
$$ \per^\T_\omega \: \pi_0(C^\infty(\T^2, G)) 
\cong \pi_1(G) \times \pi_1(G) \times \pi_2(G) \to \a $$
(cf.\ [MN03, Remark I.11(b)], [Ne02, Th.~A.3.7]). 
The restriction to $\pi_2(G)$, which corresponds to homotopy classes of 
maps vanishing on 
$(\T \times \{\1\}) \cup (\{\1\} \times \T)$, is the period map 
$\per_\omega \: \pi_2(G) \to \a$. The map 
$$ \pi_1(G) \times \pi_1(G) \to \pi_0(C^\infty(\T^2, G)) $$ 
is induced by the map 
$$ ([\gamma], [\eta]) \mapsto [\gamma*\eta]\quad \hbox{ with } \quad
(\gamma*\eta)(t,s) = \gamma(t)\eta(s), $$
and we have seen in Corollary VI.5 that the commutator map 
$\pi_1(G) \times \pi_1(G) \to A$
is given by 
$$ ([\gamma], [\eta]) \mapsto \per^\T_\omega([\gamma*\eta]) +\Gamma_A. $$
Note that this map is biadditive and {\sl not} a group homomorphism 
$\pi_1(G) \oplus \pi_1(G) \to A$, which implies that $\per^\T_\omega$ is {\sl not} a
group homomorphism. The condition 
$$ \im(\per^\T_\omega) \subeq \Gamma_A $$
means at the same time that 
$\Pi_\omega \subeq \Gamma_A$ and that the commutator map
$C$ is trivial. 
\qed 

\Remark VI.9. If $A \cong \a/\Gamma_A$, then $\a^G = \a^\g$ is a closed
subspace of $\a$ containing $\Gamma_A$. Therefore 
$$ A/A^G \cong \b := \a/\a^G $$
is a locally convex space which carries a natural smooth $G$-module
structure. Note that the quotient space $\b$ need not be
sequentially complete if $\a$ has this property. 
Nevertheless the construction in Section V leads to a group cocycle 
$f \in Z^2_s(G,\a/\Pi_\omega)$ and since $\Pi_\omega$ is always
contained in $\a^G$ (Lemma IV.2), we obtain a group cocycle 
$$ f_1 \in Z^2_s(G,\b) \quad \hbox{ with } \quad 
Df_1 = \omega^\b := q_\b \circ \omega,  $$
where $q_\b \: \a \to \b$ is the quotient map 
(Corollary V.3). This leads to a Lie group extension 
$$ \b \into \hat G \onto \tilde G $$
with $\hat\g \cong \b \oplus_{\omega^\b} \g$. 
Note that 
$$ \b = \a/\a^G \cong B^1_c(\g,\a) \subeq Z^1_c(\g,\a), $$
so that we may identify the quotient map 
$q_\b$ with the coboundary map 
$d_\g \: \a \to B^1_c(\g,\a).$ 
This makes it easier to identify the corresponding flux cocycle. 

In Proposition X.4 we shall encounter examples of modules 
$\a$ with $\a^\g = \{0\}$ for which the flux cocycle 
is non-trivial (this is the case for the module ${\cal F}_1$ of 
$\Diff(\SS^1)_0$). Therefore one cannot expect $F_{\omega_\b}$ to vanish. 
\qed

\sectionheadline{VII. An exact sequence for abelian Lie group extensions} 

\nin Let $G$ be a connected Lie group and $A$ a smooth $G$-module 
of the form $A \cong \a/\Gamma_A$, where $\Gamma_A \subeq \a$ is a discrete subgroup. 
The main result of the present section is an exact sequence relating 
the group homomorphism 
$$ D \: H^2_s(G,A) \to H^2_c(\g,\a) $$
to the exact Inflation-Restriction Sequence associated to the normal subgroup 
$\pi_1(G) \cong \ker q_G$ of $\tilde G$, where $q_G \: \tilde G \to G$
is the universal covering map (cf.\ Appendix D). The crucial
information on $\im(D)$ has already been obtained in Theorem~VI.7, so
that it essentially remains to show that $\ker D$ coincides with the
image of the connecting homomorphism $\delta \: \Hom(\pi_1(G),A^G) \to
H^2_s(G,A)$. 

In the following we shall always
consider $A$ as a $\tilde G$-module, where $g \in \tilde G$ acts on
$A$ by $g.a := q_G(g).a$, so that $\pi_1(G)$ acts trivially. 

\Proposition VII.1. Let $G$ be a connected Lie group. For 
an abelian Lie group extension $A \into \hat G \sssmapright{q} G$ the
following conditions are equivalent: 
\litem{(1)} There exists an open identity neighborhood $U \subeq G$
and a smooth section $\sigma_U \: U \to
\hat G$ of $q$ with $\sigma_U(xy)= \sigma_U(x)\sigma_U(y)$ for $x,y,xy\in~U$. 
\litem{(2)} $\hat G \cong A \times_f G$, where $f \in Z^2_s(G,A)$
is constant $0$ on an identity neighborhood in $G \times G$. 
\litem{(3)} There exists a homomorphism $\gamma \: \pi_1(G) \to A^G$
and an isomorphism $\Phi \: (A \rtimes \tilde G)/\Gamma(\gamma)
\to \hat G$ with $q\big(\Phi([\1,x])\big) = q_G(x)$, $x \in \tilde G$, where 
$\Gamma(\gamma) = \{ (\gamma(d),d) \: d \in \pi_1(G)\}$ is the graph of $\gamma$. 

\Proof. (1) $\Leftrightarrow$ (2) follows directly from the
definitions and Proposition~II.6. 

\par\nin (1) $\Rarrow$ (3): We may w.l.o.g.\ assume that $U$ is
connected, $U = U^{-1}$,  and that there exists a smooth section $\tilde \sigma \: U \to
\tilde G$ of the universal covering map $q_G$. Then 
$$\sigma_U \circ q_G\res_{\tilde\sigma(U)} \: \tilde\sigma(U) \to \hat G $$
extends uniquely to a smooth homomorphism $\phi \: \tilde G \to \hat
G$ with $\phi \circ \tilde \sigma = \sigma_U$ and $q \circ \phi = q_G$
([Ne02, Lemma 2.1]; see also [HoMo98, Cor.~A.2.26]). 
We define $\psi \: A \rtimes \tilde G \to \hat G, (a,g) \mapsto a
\phi(g)$. Then $\psi$ is a smooth group
homomorphism which is a local diffeomorphism because
$$\psi(a,\tilde\sigma(x)) = a\phi(\tilde\sigma(x)) = a\sigma_U(x) \quad
\hbox{ for } \quad x \in U, a \in A. $$
We conclude that $\psi$ is a covering homomorphism. Moreover, $\psi$
is surjective because its range is a subgroup of $\hat G$ containing $A$ 
and mapped surjectively by $q$ onto $G$. This proves that 
$$\hat G \cong (A \rtimes \tilde G)/\ker \psi, \quad 
\ker \psi = \{ (-\phi(g),g) \: g \in \phi^{-1}(A)\}. $$
On the other hand, $\phi^{-1}(A) = \ker(q \circ \phi) = \ker q_G =
\pi_1(G)$, so that 
$$ \ker \psi = \{ (\gamma(d),d) \: d \in \pi_1(G)\} = \Gamma(\gamma) \quad 
\hbox{ for } \quad  \gamma:= -\phi\res_{\pi_1(G)}. $$

\par\nin (3) $\Rarrow$ (1) follows directly from the fact that the
map $A \rtimes \tilde G \to \hat G$ is a covering morphism. 
\qed

For the following theorem we recall the definition of the period 
map $\per_\omega$ (Section IV) and the flux homomorphism $F_\omega \:
\pi_1(G) \to H^1_c(\g,\a)$ associated to
$\omega \in Z^2_s(\g,\a)$ (Proposition VI.3). 

\Theorem VII.2. Let $G$ be a connected Lie group, $A$ a smooth $G$-module 
of the form $A \cong \a/\Gamma_A$, where $\Gamma_A \subeq \a$ is a
discrete subgroup of the sequentially complete locally convex space $\a$ 
and $q_A \: \a \to A$ the quotient map. The map 
$$ \tilde P \: Z_c^2(\g,\a) \ \to  \Hom\big(\pi_2(G),A\big)
\times \Hom\big(\pi_1(G), H^1_c(\g,\a)\big), \quad 
\tilde P(\omega) = (q_A \circ \per_\omega, F_\omega) $$
factors through a homomorphism 
$$ P \: H_c^2(\g,\a) \ \to  \Hom\big(\pi_2(G),A\big)
\times \Hom\big(\pi_1(G), H^1_c(\g,\a)\big), \quad 
P([\omega]) = (q_A \circ \per_\omega, F_\omega) $$
and the following sequence is exact: 
$$ \eqalign{
\0 &\to  H^1_s(G,A) \sssmapright{I} H^1_s(\tilde  G,A) \sssmapright{R}  
H^1\big(\pi_1(G),A\big)^G \cong \Hom\big(\pi_1(G),A^G\big) 
\ssmapright{\delta} \cr 
&\ \ \ \ \ssmapright{\delta} H^2_s(G,A) 
\ssmapright{D} H_c^2(\g,\a) \ssmapright{P} \Hom\big(\pi_2(G),A\big)
\times \Hom\big(\pi_1(G), H^1_c(\g,\a)\big).  \cr}$$
Here the map $\delta$ assigns to a group homomorphism $\gamma \: \pi_1(G)
\to A^G$ the quotient of the semi-direct product 
$A \rtimes \tilde G$ by the graph $\{ (\gamma(d),d) \: d \in
\pi_1(G)\}$ of $\gamma$ which is a discrete central subgroup. 

\Proof. First we verify that $\tilde P$ vanishes on $B^2_c(\g,\a)$, so
that the map $P$ is well-defined. 
In Theorem VI.7 we have seen that $[\omega] \in \im(D)$ is equivalent
to $\tilde P(\omega) = 0$. If $[\omega] = 0$, then 
$\a \oplus_\omega \g \cong \a \rtimes \g$ and the semi-direct product 
$A \rtimes G$ is a corresponding extension of $G$ by $A$, so that 
Theorem~VI.7 leads to $\tilde P(\omega) =0$. As 
$\tilde P$ is a group homomorphism, it factors to a homomorphism 
$P$ on $H^2_c(\g,\a)$. 

The exactness of the sequence in $H^1_s(G,A)$, 
$H^1_s(\tilde G,A)$ and $\Hom(\pi_1(G), A^G)$ follows from
Example~D.11(b) and the exactness in $H^2_c(\g,\a)$ from Theorem~VI.7.
It therefore remains to verify the exactness in $H^2_s(G,A)$. 

First we need a more concrete interpretation of the map $\delta$ in
terms of abelian extensions. Let $\gamma \in \Hom(\pi_1(G), A^G)$ and 
$f \in C^1_s(\tilde G, A)$ as in Lemma~D.7 applied with $N = \pi_1(G)$ with 
$f(gd) = f(g) + \gamma(d)$ for $g \in \tilde G, d \in \pi_1(G)$. 
Then the arguments in Remark~D.10 show that the map 
$$ \Phi \: A \times_{d_{\tilde G} f} \tilde G \to A \rtimes \tilde G,
\quad (a,g)  \mapsto (a + f(g), g) $$
is a bijective group homomorphism. Since, in addition, $\Phi$ is
a local diffeomorphism, it also is an isomorphism of Lie groups, and
therefore the cocycle $\delta(f) := \oline{d_{\tilde G} f} \in
Z^2_s(G,A)$ satisfies 
$$ A \times_{\delta(f)} G 
\cong (A \times_{d_{\tilde G} f} \tilde G)/(\{\0\} \times \pi_1(G)) 
\cong (A \rtimes \tilde G)/\Phi(\{\0\} \times \pi_1(G)) 
\cong (A \rtimes \tilde G)/\{(d,\gamma(d)) \: d \in \pi_1(G)\}. $$
Now the inclusion $\im(\delta) \subeq \ker(D)$ follows from Proposition~VII.1 because 
for a cocycle $f \in Z^2_s(G,A)$ vanishing in an identity
neighborhood we clearly have $Df = 0$. 

Conversely, let $f \in Z^2_s(G,A)$ be a locally smooth group cocycle for which 
$\omega := Df$ is a coboundary and let 
$q \: \hat G = A \times_f G \to G$ be a corresponding 
Lie group extension (Proposition~II.6). 
Then the Lie algebra extension 
$\hat\g \cong \a \oplus_\omega \g \to \g$ splits, and there exists a continuous projection 
$p_\a \: \hat\g \to \a$
whose kernel is a closed subalgebra isomorphic to $\g$. Considering 
$p_\a$ as an element of 
$C^1_c(\hat\g,\a)$, we have 
$$  (d_\g p_\a)(x,y) = x.p_\a(y) - y.p_\a(x) - p_\a([x,y]) 
= p_\a([x-p_\a(x), p_\a(y)-y]) = 0, $$
for $x,y \in \hat\g$, so that $p_\a \in Z^1_c(\hat\g,\a)$. 
Let $q_{\hat G} \: G^\sharp \to \hat G$ denote the universal covering group of $\hat G$. 
Then the corresponding equivariant $1$-form $p_\a^{\rm eq}$ on
$G^\sharp$ is closed (Lemma B.5), 
so that we find a smooth function 
$$\phi \: G^\sharp \to \a \quad \hbox{ with } \quad \phi(\1) = 0 \quad \hbox{ and } \quad 
d\phi = p_\a^{\rm eq}, $$
and Lemma III.2 implies that $\phi \in Z^1_s(\hat G,\a)$ is a group cocycle. 

Using the local description of $\hat G$, resp., $G^\sharp$ by a $2$-cocycle, 
we see that the inclusion map $A_0 \into \hat G$ of the identity component of $A$ lifts to a Lie group morphism  
$\eta_\a \: \a \to G^\sharp$
whose differential is the inclusion $\a \into \hat\g$. Since $p_\a\res_\a = \id_\a$ and 
the image of $\eta_\a$ acts trivially on $\a$, the 
composition $\phi \circ \eta_\a \: \a\to \a$ is a morphism of Lie groups whose differential is $\id_\a$, 
which implies that $\phi \circ \eta_\a = \id_\a.$
Moreover, the cocycle condition implies that 
$$ \phi(ag) = \phi(a)+ \phi(g), \quad a \in \eta_\a(\a), g \in G^\sharp. \leqno(7.1) $$

Let $U \subeq G$ be a connected open identity neighborhood on which there exists a smooth 
section $\sigma \: U \to G^\sharp$ of the quotient map $q^\sharp := q
\circ q_{\hat G} \: G^\sharp \to G$. We then obtain another smooth map by   
$$ \sigma_1 \: U \to G^\sharp, \quad x \mapsto \eta_\a(\phi(\sigma(x))^{-1}) \sigma(x). $$
In view of (7.1), this map is also a section of $q^\sharp$. Moreover, 
$\im(\sigma_1) \subeq \phi^{-1}(0).$

From the description of $\hat G$ with the cocycle $f$ it follows that there exists an open 
$\1$-neighborhood in $G^\sharp$ of the form 
$$ U^\sharp := \eta_\a(U_\a) \sigma_1(U), $$
where $U_\a \subeq \a$ is an open $0$-neighborhood. Restricting $\phi$ to $U^\sharp$, we see that 
$\sigma_1(U) = \phi^{-1}(0) \cap U^\sharp.$
Since $\phi^{-1}(0)$ is a subgroup of $G^\sharp$, we have 
$$ (\sigma_1(U) \sigma_1(U)) \cap U^\sharp \subeq \sigma_1(U). $$
Let $V \subeq U$ be an open symmetric $\1$-neighborhood in $G$ 
such that there exists a smooth section 
$\sigma_V \: V \to \tilde G$ of the universal covering map $q_G \: \tilde G \to G$ 
and, in addition, $V V \subeq U$ and $\sigma_1(V)\sigma_1(V) \subeq U^\sharp$. 
For $x,y \in V$ we then have $xy \in U$, and 
$\sigma_1(x) \sigma_1(y) \in U^\sharp$ implies the existence of $z \in U$ with 
$ \sigma_1(z) = \sigma_1(x)\sigma_1(y)$. Applying $q^\sharp$ to both sides leads to 
$$z = q^\sharp\sigma_1(z) = q^\sharp(\sigma_1(x)\sigma_1(y)) = xy.$$
We therefore have 
$$ \sigma_1(xy) = \sigma_1(x) \sigma_1(y) \quad \hbox{ for } \quad x,y \in V. $$
Hence there exists a unique group homomorphism 
$ f \: \tilde G \to G^\sharp$
with $f \circ \sigma_V = \sigma_1$ ([HoMo98, Cor.\ A.2.26]). 
Composing $f$ with the covering map $q_{\hat G} \: G^\sharp \to \hat G$, we obtain a smooth 
homomorphism $\hat f \: \tilde G \to \hat G$ with $q \circ \hat f =
q_G$. In view of Proposition~VII.1, 
this implies that $\hat G$ is isomorphic to a group of the type 
$(A \rtimes \tilde G)/\Gamma(\gamma),$
where $\gamma \: \pi_1(G) \to A^G$ is a group homomorphism. 
\qed

Since the fundamental group $\pi_1(\tilde G)$ vanishes, we obtain in
particular: 
\Corollary VII.3. The map 
$\tilde D \: H^2_s(\tilde G,A) \to H^2_s(\g,\a)$
is injective. 
\qed

In view of Corollary VII.3, 
we may identify $H^2_s(\tilde G,A)$ with a subgroup of 
$H^2_c(\g,\a)$. The inflation map 
$$ I \: H^2_s(G,A) \to H^2_s(\tilde G, A) \quad \hbox{ satisfies }
\quad 
\tilde D \circ I = D \: H^2_s(G,A) \to H^2_c(\g,\a). $$

\Remark VII.4. At first sight, the following argument 
seems to be more natural to prove that 
$\ker D \subeq \im \delta$: If the group $\hat G$ is 
regular (cf.\ [Mil83]), then the Lie algebra morphism 
$\sigma \: \g \to \hat\g$ whose existence is 
guaranteed by $[Df] = 0$ can be integrated to a Lie group 
morphism $\tilde G \to \hat G$, and we can argue as above. Unfortunately this argument requires the 
regularity of the group $\hat G$, which is not needed for the argument given above. 
\qed

\sectionheadline{VIII. Abelian extensions with smooth global sections} 

\nin In this subsection we discuss the existence of a smooth
cross section for an abelian Lie group extension $A \into \hat G \onto G$ which is 
equivalent to the existence of a smooth global cocycle 
$f \: G \times G \to A$ with $\hat G \cong G \times_f A$. Moreover, we
will show that for simply connected groups, it is equivalent to the
exactness of the equivariant $2$-form $\omega^{\rm eq}$ on $G$,
where $\omega = Df$. 

The following lemma will be helpful in the proof of
Proposition~VIII.2. 

\Lemma VIII.1. Let $G$ be a connected Lie group, $A$ a smooth
$G$-module and $f \in Z^2_s(G,A)$ such that all functions 
$f_g \: G \to A, x \mapsto f(g,x)$ are smooth. Then 
$f \: G \times G \to A$ is a smooth function. 

\Proof. We write the cocycle condition as 
$$  f(xy,z) = f(x,yz) + \rho_A(x).f(y,z)- f(x,y), \quad x,y,z \in G. $$
For $x$ fixed, this function is smooth as a function of the pair
$(y,z)$ in a neighborhood of $(\1,\1)$. This implies that 
$f$ is smooth on a neighborhood of the points $(x,\1)$, $x \in G$. 
Fixing $x$ and $z$ shows that there exists a $\1$-neighborhood $V
\subeq G$ (independent of $x$) such that the functions $f(\cdot, z)$, $z \in V$, are smooth
in a neighborhood of $x$. Since $x \in G$ was arbitrary, we conclude
that the functions $f(\cdot, z)$, $z \in V$, are smooth.  
Now 
$$ f(\cdot, yz) = f(\cdot y,z) - \rho_A(\cdot).f(y,z) + f(\cdot,y)$$
shows that the same holds for the functions $f(\cdot, u)$, $u \in V^2$. 
Iterating this process, using $G = \bigcup_{n \in \N} V^n$, we
derive that all functions $f(\cdot, x)$, $x \in G$, are smooth. 
Finally we see that the function 
$$ (x,y) \mapsto  f(x,yz)  = f(xy,z) - \rho_A(x).f(y,z) + f(x,y) $$
is smooth in a neighborhood of each point $(x_0, \1)$, hence that $f$
is smooth in each point $(x_0, z_0)$, and this proves that $f$ is
smooth on $G \times G$. 
\qed

\Proposition VIII.2. Let $G$ be a 
connected Lie group, $\a$ a sequentially complete locally convex 
smooth $G$-module, $\omega \in Z^2_c(\g,\a)$ a
continuous $2$-cocycle, and $\omega^{\rm eq} \in \Omega^2(G,\a)$ 
the corresponding equivariant $2$-form on $G$ with $\omega^{\rm eq}_\1 = \omega$. 
We assume that 
\litem{(1)} $\omega^{\rm eq} = d\theta$ for some $\theta \in \Omega^1(G,\a)$
and 
\litem{(2)} for each $g \in G$ the closed $1$-form $\lambda_g^* \theta - \rho_\a(g) \circ \theta$ is exact. 

\par\nin Then the product manifold $\hat G := \a \times G$ carries a
Lie group structure which is given by a smooth $2$-cocycle $f \in
Z^2_s(G,\a)$ with $D[f] = [\omega]$ via  
$$ (a,g) (a',g') := (a + g.a' + f(g,g'),gg'). $$

\Proof. For each $g \in G$ the relation $\rho_\a(g) \circ \omega^{\rm eq} = \lambda_g^*\omega^{\rm eq}$ implies  
$$ d\big(\rho_\a(g) \circ \theta - \lambda_g^*\theta\big) = 
\rho_\a(g) \circ \omega^{\rm eq} - \lambda_g^*\omega^{\rm eq} = 0.$$
In view of (2), for each $g \in G$ there exists a smooth function $f_g \: G \to \a$
with $f_g(\1) = 0$ and 
$$df_g = \lambda_g^*\theta - \rho_\a(g) \circ \theta. $$
Observe that $f_\1 = 0$. For $g, h \in G$ this leads to 
$$ \eqalign{ df_{gh} 
&= \lambda_{gh}^*\theta - \rho_\a(gh) \circ\theta
= \lambda_h^*(\lambda_{g}^*\theta - \rho_\a(g)\circ\theta) + \lambda_h^*(\rho_\a(g)\circ\theta) - \rho_\a(gh) \circ\theta \cr 
&= \lambda_h^*df_g + \rho_\a(g) (\lambda_h^*\theta - \rho_\a(h) \circ\theta) 
= \lambda_h^*df_g + \rho_\a(g) \circ df_h 
= d(f_g \circ \lambda_h + \rho_\a(g) \circ f_h). \cr} $$
Comparing values of both functions in $\1$, we get 
$$ f_{gh} = f_g \circ \lambda_h + \rho_\a(g) \circ f_h - f_g(h). \leqno(8.1) $$
Now we define $f \: G \times G \to \a$ by $f(x,y) := f_x(y)$. Then (8.1)  means that 
$$ f(gh,u) = f(g,hu) + \rho_\a(g).f(h,u) - f(g,h), \quad g,h,u \in G, $$
i.e., $f$ is a group cocycle. 

Moreover, the concrete local formula for $f_x$ in the Poincar\'e Lemma ([Ne02, Lemma 3.3]) 
and the smooth dependence of the integral on $x$ imply that 
$f$ is smooth on a neighborhood of $(\1,\1)$, so that Lemma~VIII.1
implies that $f \: G \times G \to \a$ is a smooth function. 
We therefore obtain on the space 
$\hat G := \a \times G$
a Lie group structure with the multiplication given by 
$$ (a,g) (a',g') := (a + g.a' + f(g,g'), gg') $$
(Lemma II.1), 
and Lemma II.7 implies that the corresponding Lie bracket is given by 
$$ [(a,x), (a',x')] 
= \big(x.a' - x'.a + d^2 f(\1,\1)(x,x') - d^2 f(\1,\1)(x',x), [x,x']\big).$$

Now we relate this formula to the Lie algebra cocycle $\omega$. 
The relation $df_g = \lambda_g^* \theta - \rho_\a(g)\circ\theta$ leads to 
$$ df(g,\1)(0,y) =  df_g(\1)y = (\lambda_g^* \theta - \rho_\a(g)\circ \theta)_\1(y) 
= \la \theta, y_l \ra(g) - \rho_\a(g).\theta_\1(y), $$
where $y_l$ denotes the left invariant vector field with $y_l(\1) =
y$. Taking second derivatives, we further obtain for $x \in \g$: 
$$ \eqalign{ 
&\ \ \ \ d^2 f(\1,\1)(x, y)  \cr
&=  x_l(\la \theta, y_l \ra)(\1) - x.\theta_\1(y) 
=  (d\theta)(x_l, y_l)(\1) + y_l(\la \theta, x_l \ra)(\1) 
+ \theta([x_l,y_l])(\1)- x.\theta_\1(y) \cr
&=  \omega(x, y) + y_l(\la \theta, x_l \ra)(\1) + \theta_\1([x,y])- x.\theta_\1(y), \cr} $$
Subtracting 
$d^2 f(\1,\1)(y,x)  =  y_l(\la \theta, x_l \ra)(\1) - y.\theta_\1(x),$
leads to 
$$ (Df)(x,y) 
=  \omega(x, y) + \theta_\1([x,y]) - x.\theta_\1(y) + y.\theta_\1(x) 
=  \omega(x, y) - (d\theta_\1)(x,y). $$
Since this cocycle is equivalent to $\omega$, the assertion follows. 
\qed

\Corollary VIII.3. If $G$ is simply connected and 
$\omega^{\rm eq}$ is exact, then there exists a smooth cocycle $f \: G\times G
\to \a$ with $D[f] = [\omega]$, so that $\hat G := \a \times_f G$ is a Lie group with Lie
algebra $\hat\g = \a \oplus_\omega \g$. 

\Proof. Since $\pi_1(G)$ is trivial, condition (2) in
Proposition~VIII.2 is automatically satisfied. 
\qed

For central extensions of finite-dimensional groups, 
the construction described in Proposition VIII.2 is  
due to E.~Cartan, 
who used it to construct a central extension of a simply connected
finite-dimensional Lie group $G$ by the group $\a$. Since in this case 
$$H^2_{\rm dR}(G,\a) \cong \Hom(\pi_2(G),\a)= {\bf 0} 
\quad \hbox{ and } \quad 
H^1_{\rm dR}(G,\a) \cong \Hom(\pi_1(G),\a) = {\bf 0},
$$
(cf.\ [God71]), the requirements of the construction are satisfied for
every Lie algebra cocycle $\omega \in Z^2_c(\g,\a)$. 

\Proposition VIII.4. If $G$ is a connected Lie group which is smoothly
paracompact, then the conclusion of {\rm Proposition VIII.2} remains
valid under the assumptions: 
\litem{(1)} $\omega^{\rm eq}$ is an exact $2$-form, and 
\litem{(2)} $F_\omega = 0$. 

\Proof. In view of (1), we can apply Proposition VIII.2 to the universal covering 
group $q_G \: \tilde G \to G$ of $G$, which leads to an $\a$-extension  
$$ q^\sharp \: G^\sharp := \a\times_f \tilde G \to \tilde G, \quad
(a,g) \mapsto g, $$
where $f \in Z^2_s(\tilde G,\a)$ is a smooth cocycle with 
$D[f] = [\omega]$. 
In view of Corollary VI.5, the vanishing of $F_\omega$ implies
the existence of a homomorphism $\gamma \: \pi_1(G) \to Z(G^\sharp)$ with
$q^\sharp \circ \gamma = \id_{\pi_1(G)}$. Then 
$\im(\gamma)$ is a discrete central subgroup of $G^\sharp$, so that 
$\hat G :=  G^\sharp/\im(\gamma)$ is a Lie group, and 
we obtain an $\a$-extension of $G$ by 
$$ q \: \hat G \to G, \quad g \im(\gamma)
\mapsto q_G \circ q^\sharp(g). $$

As $\hat G$ is a principal $\a$-bundle over $G$, its fibers are 
affine spaces whose translation group is~$\a$. 
If $G$ is smoothly paracompact, we can therefore use a smooth
partition of unity subordinated to a trivializing open cover of the
$\a$-bundle $\hat G \to G$ to patch smooth local sections together to a global
smooth section $\sigma \: G \to \hat G$. Then the map 
$$ \a \times_{f_G} G \to \hat G, \quad (a,g) \mapsto a \sigma(g) $$
is an isomorphism of Lie groups, where 
$f_G \in Z^2_s(G,\a), (g,g') \mapsto \sigma(g)\sigma(g')
\sigma(gg')^{-1}$
is a globally smooth cocycle. 
\qed

\Remark VIII.5. Let $G$ be a connected Lie group and $A$ a smooth
$G$-module of the form $\a/\Gamma_A$. Let $Z^2_{gs}(G,A)$ denote the
group of smooth $2$-cocycles $G \times G \to A$ 
and $B^2_{gs}(G,A) \subeq Z^2_{gs}(G,A)$ the cocycles of the form $d_G
h$, where $h \in C^\infty(G,A)$ is a smooth function with $h(\1) =
0$. Then one can show that we have an injection 
$$ H^2_{gs}(G,A) := Z^2_{gs}(G,A)/B^2_{gs}(G,A) \into H^2_s(G,A), $$
the space $H^2_{gs}(G,A)$ classifies those $A$-extensions of $G$ with
a smooth global section, and we have an exact sequence 
$$ \Hom(\pi_1(G),\a^G) \sssmapright{\delta} H^2_{gs}(G,A) 
\sssmapright{D} H^2_c(\g,\a) \ssmapright{P} H^2_{\rm dR}(G,\a) 
\times \Hom\big(\pi_1(G), H^1_c(\g,\a)\big), $$
where $P([\omega]) = ([\omega^{\rm eq}], F_\omega).$ 
The proof is an easy adaptation from the corresponding arguments for
central extensions in Section~8 of [Ne02]. 
\qed

\sectionheadline{IX. Applications to diffeomorphism groups} 

In the present section we apply the general results of this paper to 
diffeomorphism groups of a compact manifold $M$. In this case the Lie algebra 
is the Fr\'echt--Lie algebra ${\cal V}(M)$ of smooth vector fields on $M$ 
and we obtain interesting Lie algebra $2$-cocycles with values in 
the space $C^\infty(M,V)$ of smooth $V$-valued functions from closed $V$-valued 
$2$-forms on $M$. In this case the period map and the flux cocycle can be made 
more concrete in geometric terms which makes it possible to evaluate the 
obstructions to the existence of abelian extensions in many concrete examples.

\Definition IX.1. Let $M$ be a compact manifold. 

\par\nin (a) We write $\Diff(M)$ for the group
of all diffeomorphisms of
$M$ and ${\cal V}(M)$ for the Lie algebra of smooth vector fields on $M$,
i.e., the set of all smooth maps $X \: M \to TM$ with 
$\pi_{TM} \circ X = \id_M$, where $\pi_{TM} \: TM \to M$ is the bundle
projection of the tangent bundle. 
We define the Lie algebra structure on ${\cal V}(M)$ in
such a way that $[X,Y].f = X.(Y.f) - Y.(X.f)$ holds for 
$X,Y \in {\cal V}(M)$ and $f \in C^\infty(M,\R)$. 

Then $\Diff(M)$ is a Lie group whose Lie algebra is ${\cal V}(M)^{\rm op}$
(the same space with the apposite bracket $(X,Y) \mapsto -[X,Y]$) and
we have a smooth exponential function 
$$ \exp \: {\cal V}(M) \to \Diff(M) $$
given by $\exp(X) = \Phi_X^1$, where $\Phi_X^t \in \Diff(M)$ is
the flow of the vector field $X$ at time $t$  
([KM97]). 

The tangent bundle of $\Diff(M)$ can be identified with the set 
$$ T(\Diff(M)) := \{ X \in C^\infty(M,TM) \: \pi_{TM} \circ X \in
\Diff(M)\},  $$
where the map 
$$ \pi \: T(\Diff(M)) \to \Diff(M), \quad X \mapsto \pi_{TM}
\circ X $$
is the bundle projection. Then 
$T_\phi(\Diff(M)) := \pi^{-1}(\phi)$
is the fiber over the diffeomorphism $\phi$. 

In view of the natural
action of $\Diff(M)$ on $TM$ given by 
$\psi.v := T(\psi).v$, we obtain natural left and right
actions of $\Diff(M)$ on $T(\Diff(M))$ by 
$$ (\phi.X)(x) = \phi(x).X(x), \quad 
X.\phi := X \circ \phi. $$
Then 
$$ \pi_{TM} \circ (\phi.X) = \phi \circ (\pi_{TM} \circ X)
\quad \hbox{ and } \quad 
\pi_{TM} \circ (X \circ \phi) = (\pi_{TM} \circ X) \circ \phi, $$
so that the left, resp., right action of $\Diff(M)$ on 
$T(\Diff(M))$ covers the left, resp., right multiplication action
of the group $\Diff(M)$ on itself.  In the following we shall mostly
consider the opposite group $\Diff(M)^{\rm op}$ whose Lie algebra is
${\cal V}(M)$. The adjoint action of this group is given by 
$$ \Ad \: \Diff(M)^{\rm op} \times {\cal V}(M) \to {\cal V}(M), \quad 
(\phi, X) \mapsto \phi^{-1}.(X \circ \phi) = \phi^{-1}.(X.\phi). $$

\par\nin (b) Let $J \subeq \R$ be an interval and $\phi \: J \to
\Diff(M)^{\rm op}$ be 
a smooth curve. Then for each $t \in J$ we obtain a vector field 
$$ \delta^r(\phi)(t) := \phi(t)^{-1}.\phi'(t) $$
called the {\it right logarithmic derivative of $\phi$ in
$t$}. We likewise define the {\it left logarithmic derivative} by 
$$ \delta^l(\phi)(t) := \phi'(t) \circ \phi(t)^{-1}. 
\qeddis

\Definition IX.2. Let $M$ be a compact smooth manifold and
$\g := {\cal V}(M)$ the Lie algebra of smooth vector fields on
$M$. If $V$ is Fr\'echet space and $\a := C^\infty(M,V)$ the
space of smooth $V$-valued functions on $M$, then 
$(X.f)(p) := df(p)X(p)$ turns $C^\infty(M,V)$ into a topological ${\cal V}(M)$-module. 
We observe that $C^\infty(M,V)$ 
and ${\cal V}(M)$ are Fr\'echet modules of the Fr\'echet 
algebra $R := C^\infty(M,\R)$. 

In the Lie algebra complex $(C^p_c(\g,\a),d_\g)_{p \in \N_0}$ formed by the
continuous alternating maps $\g^p \to \a$, we have the subcomplex
given by the subspaces 
$C^p_R(\g,\a) \subeq C^p_c(\g,\a)$ consisting of $R$-multilinear maps $\g^p
\to \a$. Using partitions of unity, it is easy to see that the
elements of $C^p_R(\g,\a)$ can be identified with smooth $V$-valued
$p$-forms, so that $C^p_R(\g,\a) \cong \Omega^p(M,V)$ ([Hel78]), and the de
Rham differential coincides with the Lie algebra differential $d_\g$
to $C^p_R(\g,\a)$. 

We thus obtain natural maps $Z^p_{\rm dR}(M,V) \to Z^p_c(\g,\a)$ and 
$j_p \: H^p_{\rm dR}(M,V) \to H^p_c(\g,\a).$
\qed

\Lemma IX.3. If $M$ is connected, then 
$V \cong C^\infty(M,V)^{{\cal V}(M)} = \a^\g$
consists of the constant functions $M \to V$. 
\qed

\Lemma IX.4. The map $j_1 \: H^1_{\rm dR}(M,V) \to H^1_c(\g,\a)$ is
injective. 

\Proof. Let $\alpha \in \Omega^1(M,V)$ be a closed $V$-valued
$1$-form on $M$. If $j_1([\alpha]) = 0$, then there exists an element 
$f \in \a = C^\infty(M,V)$ with $\alpha = d_\g f$, which means
that $\alpha = df$. Hence $\alpha$ is exact and therefore $j_1$ is
injective. 
\qed

Lemma VI.1 in [MN03] implies that we have a smooth action of the group
$G := \Diff(M)_0^{\rm op}$ on $\a$ by $\phi.f := f \circ \phi$. The derived
action of ${\cal V}(M)$ on this space is given by 
$$ (X.f)(p) = \derat0 (\exp(tX).f)(p) = \derat0 f(\exp(tX).p) 
= df(p)X(p) $$
which is compatible with Definition~IX.2. 
We view each smooth $V$-valued
$2$-form $\omega_M \in \Omega^2(M,V)$ as an element $\omega_\g \in
C^2_c(\g,\a)$. In the following we shall obtain some 
information on the period map and the flux homomorphism 
$$ \per_\omega \: \pi_2(\Diff(M)) \to \a^\g \cong V  
\quad \hbox{ and } \quad 
F_\omega \: \pi_1(\Diff(M)) \to H^1_c(\g,\a) $$
which makes it possible to verify the integrability criteria from 
Sections VI and VII in many special cases. 

\subheadline{More on the period group} 

The following proposition is very helpful in verifying the
discreteness of the image of the period map for the group 
$G := \Diff(M)_0^{\rm op}$. In the following we write 
$(m,g) \mapsto g(m)$ for the canonical right action of $G$ on
$M$. 

\Proposition IX.5. Let $\omega_M \in Z^2_{\rm dR}(M,V)$ be a closed
$V$-valued $2$-form on $M$, 
$\sigma \: \SS^2 \to G = \Diff(M)_0^{\rm op}$ smooth and $m \in M$. Then 
$$ \per_{\omega_\g}([\sigma])(m) = \int_{\eta_m \circ \sigma} \omega_M
\in V \cong C^\infty(M,V)^{{\cal V}(M)},
$$
where $\eta_m \: G \to M, g \mapsto g(m)$. 
In particular the period group $\Pi_{\omega_\g} =
\im(\per_{\omega_\g})$ is contained in the group 
$\int_{\pi_2(M)} \omega_M$ of spherical periods of $\omega_M$. 

\Proof. Since $\a^\g$ consists of constant functions $M \to V$, it
suffices to calculate the value of $\per_{\omega_\g}([\sigma]) \in
C^\infty(M,V)$ in the point~$m$. 

We claim that 
$$ \eta_m^*\omega_M = \ev_{m} \circ \omega_\g^{\rm eq}, \leqno(9.1) $$
where $\ev_m \: C^\infty(M,V) \to V$ is the evaluation in $m$. 
First we note that for $g \in G$ we have 
$\eta_m \circ \lambda_g = \eta_{g(m)}.$
Further 
$$ d\eta_m(\1)(X) = \derat0 \exp(tX).m = X(m) \quad \hbox{ for }
\quad X \in {\cal V}(M). $$
For $g \in G$ and vector fields $X, Y \in \g = {\cal V}(M)$ this leads to 
$$ \eqalign{
&\ \ \ \ (\eta_m^*\omega_M)(g.X, g.Y) \cr  
&= \omega_M(g(m))(d\eta_m(g)d\lambda_g(\1).X, d\eta_m(g)d\lambda_g(\1).Y)\cr 
&= \omega_M(g(m))(d(\eta_m \circ \lambda_g)(\1).X, d(\eta_m \circ \lambda_g)(\1).Y)\cr 
&= \omega_M(g(m))(d\eta_{g(m)}(\1).X, d\eta_{g(m)}(\1).Y)\cr 
&= \omega_M(g(m))(X(g(m)), Y(g(m))) 
= \Big(g.\big(\omega_\g(X,Y)\big)\Big)(m) 
= (\ev_m \circ \omega_\g^{\rm eq})(g.X, g.Y). \cr}$$
This proves (9.1). 
We now obtain 
$$ \eqalign{ \per_{\omega_\g}([\sigma])(m) 
&= \ev_m \int_\sigma \omega_\g^{\rm eq} 
= \int_\sigma \ev_m \circ \omega_\g^{\rm eq} 
= \int_\sigma \eta_m^*\omega_M 
= \int_{\eta_m \circ \sigma} \omega_M. \cr} 
\qeddis 

We immediately derive the following sufficient criterion for the
discreteness of $\im(\per_{\omega_\g})$. 
\Corollary IX.6. If the subgroup $\int_{\pi_2(M)} \omega_M := \{
\int_\sigma \omega_M \: \sigma \in C^\infty(\SS^2,M) \} \subeq
V$ of spherical periods of $\omega_M$ 
is discrete, then the image of $\per_{\omega_\g}$ is discrete. 
\qed

\Example IX.7. (1) The preceding corollary applies in particular to
all manifolds $M$ for which 
$\pi_2(M)/\tor(\pi_2(M))$ is a cyclic group. In fact, for each torsion
element $[\sigma] \in \pi_2(M)$ we have $\int_\sigma \omega_M = 0$, so
that $\int_{\pi_2(M)} \omega_M$ is the image of the cyclic group
$\pi_2(M)/\tor(\pi_2(M))$, hence cyclic and therefore discrete. 

Examples of such manifolds are spheres and tori:  
$$ \pi_2(\SS^d) \cong \cases{ 
\{0\} & for $d \not= 2$ \cr 
\Z & for $d = 2$ \cr } 
\quad \hbox{ and } \quad \pi_2(\T^d) \cong \pi_2(\R^d) = \{0\}, \quad  d\in \N. $$

The only compact connected manifolds $M$ with $\dim M \leq 2$ and $\pi_2(M)$
non-trivial are the $2$-sphere $\SS^2$ and the real projective
plane $\P_2(\R)$. This follows from $\pi_2(M) \cong \pi_2(\tilde M)$
for the universal covering $\tilde M \to M$ and the fact that a simply
connected $2$-dimensional manifold is diffeomorphic to $\SS^2$
or $\R^2$. Further all orientable $3$-manifolds which are irreducible in the
sense of Kneser have trivial $\pi_2$. 
In particular the complement of a knot $K \subeq \SS^3$ has trivial
$\pi_2$ (cf.\ [Mil03, p.1228]). 

\par\nin (2) For $M = \SS^2$ we have 
$$ \pi_2(\Diff(M)) \cong \pi_2(\SO_3(\R)) = \{\1\} 
\quad \hbox{ and } \quad 
\pi_2(\SS^2) \cong \Z. $$
If $\omega_M \in Z^2_{\rm dR}(M,\R)$ is the closed $2$-form with 
$\int_M \omega_M = 1$, we have 
$\int_{\pi_2(M)} \omega_M = \Z$
which is larger than $\Pi_{\omega_\g} = \im(\per_{\omega_\g}) = \{0\}$. 
\qed 

\Problem IX. Find an example of a closed $2$-form $\omega$ for which
the group $\Pi_{\omega_\g} = \im(\per_{\omega_\g})$ 
is discrete and $\int_{\pi_2(M)} \omega_M$ is not. 
\qed

\subheadline{The flux cocycle} 

We continue with the setting where $M$ is a compact manifold and 
$G = \Diff(M)_0^{\rm op}$ is the identity component of its diffeomorphism
group endowed with the opposite multiplication. 
For any Fr\'echet space $V$ the space 
$\Omega^1(M,V)$ is a smooth $G$-module with respect to 
$(\phi,\beta) \mapsto \phi^*\beta$. To verify the
smoothness of this action, we can think
of $\Omega^1(M,V)$ as a closed subspace of $C^\infty(TM,V)$ and
observe that $\Diff(M)$ acts smoothly on $TM$, so that Lemma VI.1 in 
[MN03] applies. The corresponding derived module of $\g = {\cal V}(M)$
is given by $(X,\beta) \mapsto {\cal L}_X.\beta$, where ${\cal L}_X =
d \circ i_X + i_X \circ d$
denotes the Lie derivative. The subspace $dC^\infty(M,V)$ 
of exact $1$-forms is a closed subspace  because 
$$ dC^\infty(M,V) = \Big\{ \beta \in \Omega^1(M,V) \: 
(\forall \gamma \in C^\infty(\SS^1,M)) \int_\gamma \beta = 0\Big\}
\leqno(9.2) $$
and the linear maps $\Omega^1(M,V) \to V, \beta \mapsto \int_\gamma \beta$ are continuous. 
We can therefore form the quotient module 
$$ \hat H^1_{\rm dR}(M,V) := \Omega^1(M,V)/dC^\infty(M,V) $$ 
containing $H^1_{\rm dR}(M,V) = Z^1_{\rm dR}(M,V)/dC^\infty(M,V)$ as a closed subspace. 

\Lemma IX.8. For each closed $V$-valued $2$-form $\omega \in
\Omega^2(M,V)$ the continuous linear map 
$$ f_{\omega} \: {\cal V}(M) \to  \hat H^1_{\rm dR}(M,V), \quad X
\mapsto [i_X \omega] $$ 
is a Lie algebra $1$-cocycle. 

\Proof. For $X, Y \in {\cal V}(M)$ we use the formulas 
$i_{[X,Y]} = [{\cal L}_X, i_Y]$ and ${\cal L}_X = i_X \circ d + d
\circ i_X$ to
obtain 
$$\eqalign{  i_{[X,Y]}\omega 
&= {\cal L}_X i_Y \omega - i_Y {\cal L}_X \omega 
= d i_X i_Y \omega + i_X d(i_Y \omega) - i_Y (d i_X \omega + i_X d  \omega)  \cr
&= d i_X i_Y \omega + i_X d(i_Y \omega) - i_Y (d i_X \omega). \cr} $$  
In view of 
$[{\cal L}_X i_Y \omega] 
= [d i_X  i_Y \omega + i_X d i_Y \omega] 
= [i_X d i_Y \omega]$
in $\hat H^1_{\rm dR}(M,V)$, 
this means that 
$$ f_\omega([X,Y]) = X.f_\omega(Y) - Y.f_\omega(X), $$
i.e., $f_\omega$ is a cocycle. 
\qed

\Definition IX.9. Let $q_G \: \tilde G \to G$ denote the universal
covering morphism of $G
= \Diff(M)_0^{\rm op}$ and define the $\tilde G$-action on $C^\infty(M,V)$,
$\Omega^1(M,V)$, $\hat H^1_{\rm dR}(M,V)$ etc.\ by pulling it back
with $q_G$ to $\tilde G$. Then Proposition~III.4 implies that there exists a
smooth $1$-cocycle 
$$ F_\omega \: \tilde G \to \hat H^1_{\rm dR}(M,V) =\Omega^1(M,V)/dC^\infty(M,V) \quad \hbox{
with } \quad d F_\omega(\1) = f_\omega. $$ 
This cocycle is called the {\it flux cocycle corresponding to $\omega$}. 
Its differential $d F_\omega$ coincides with the equivariant $1$-form
$f_\omega^{\rm eq}$. 
\qed

\Remark IX.10. (a) If $g \in \tilde G$ and $\tilde \gamma\: [0,1] \to
\tilde G$ is a piecewise smooth curve with $\tilde\gamma(0) = \1$ and 
$\tilde\gamma(1) = g$, then $\tilde\gamma$ is the unique lift of 
$\gamma := q_G \circ \tilde\gamma \: [0,1] \to G$. The value of the
flux cocycle in $g$ is determined by 
$$ \eqalign{ F_\omega(g) 
&= \int_0^1 dF_\omega(\tilde\gamma(t))(\tilde\gamma'(t))\, dt  
= \int_0^1 (f_\omega^{\rm eq})(\tilde\gamma(t))(\tilde\gamma'(t))\, dt  \cr
&= \int_0^1 \gamma(t).f_\omega(\tilde\gamma(t)^{-1}.\tilde\gamma'(t))\, dt  
= \int_0^1 \gamma(t).f_\omega(\gamma(t)^{-1}.\gamma'(t))\, dt\cr 
&= \int_0^1 \gamma(t).f_\omega(\delta^l(\gamma)(t))\, dt
= \int_0^1 [\gamma(t)^*.i_{\delta^l(\gamma)(t)}\omega]\, dt \cr 
&= \int_0^1 [i_{\delta^r(\gamma)(t)}(\gamma(t)^*\omega)]\, dt 
\in \hat H^1_{\rm dR}(M,V).\cr}$$ 
Here we have used the relation 
$\phi^*(i_X \omega) = i_{\Ad(\phi).X} (\phi^*\omega)$ for 
$\phi \in \Diff(M)^{\rm op}$. 

\par\nin (b) For the special case when the curve $\gamma \: [0,1] \to \Diff(M)$ has
values in the subgroup 
$$\Sp(M,\omega) := \{ \phi \in \Diff(M) \: \phi^*\omega = \omega \}, $$
all vector fields $\delta^l(\gamma)(t)$ are contained in the Lie
algebra 
$$\sp(M,\omega) := \{ X \in {\cal V}(M) \: {\cal L}_X.\omega = 0\} $$
([NV03, Lemma I.4]). For ${\cal L}_X \omega = 0$ we have 
$d(i_X\omega) = {\cal L}_X \omega = 0$, so that all $1$-forms 
$i_X\omega$ are closed. This in turn implies that for 
each $\phi \in \Diff(M)_0$ the $1$-form $\phi^* i_X \omega - i_X
\omega$ is exact ([NV03, Lemma 1.3]). For the flux cocycle this leads
to the simpler formula 
$$ \eqalign{ F_\omega(g) 
= \int_0^1 [i_{\delta^l(\gamma)(t)}\omega]\, dt. \cr}$$
Hence $F_\omega(g)$ is the flux associated to the curve 
$\gamma \: [0,1] \to \Sp(M,\omega)$ in the context of symplectic
geometry [MDS98]. 

\par\nin (c) If the closed form $\omega$ is exact, 
$\omega = d\theta$, then 
$$ f_\omega(X)= [i_X\omega] = [i_X d\theta] = [{\cal L}_X\theta] = X.[\theta]$$
in $\hat H^1_{\rm dR}(M,V)$ implies that $f_\omega$ is a coboundary. Hence 
it integrates to a group cocycle given by 
$$ F_\omega \: \Diff(M)^{\rm op} \to \hat H^1_{\rm dR}(M,V), 
\quad \phi \mapsto 
[\phi^*\theta - \theta]. 
\qeddis 

On the space $\hat H^1_{\rm dR}(M,V)$ the integration maps 
$\hat H^1_{\rm dR}(M,V) \to V, [\beta] \mapsto \int_\alpha \beta$ for  
$\alpha \in C^\infty(\SS^1,M)$ separate points (cf.\ (9.2)), 
so that the element $F_\omega(g) \in \hat H^1_{\rm dR}(M,V)$ is determined
by the integrals $\int_\alpha F_\omega(g)$ which are evaluated in the
proposition below. 

\Proposition IX.11. For $\alpha \in C^\infty(\SS^1,M)$ and a smooth curve 
$\gamma \: [0,1] \to G = \Diff(M)_0^{\rm op}$ 
with $\gamma(0) = \id_M$ we consider the smooth map  
$$ H \: \SS^1 \times [0,1] \to M, \quad (t,s) \mapsto
\gamma(t)(\alpha(s)). $$
Let $\tilde\gamma \: [0,1] \to \tilde G$ be the smooth lift with
$\tilde\gamma(0) = \1$. Then the value of the flux cocycle in
$\tilde\gamma(1)$ is determined by the integrals 
$$ \int_\alpha F_\omega(\tilde\gamma(1)) 
= \int_H \omega. $$

\Proof. First we note that 
$$ {\partial H \over \partial t}(t,s) 
= \gamma'(t)(\alpha(s)) 
= \gamma'(t) \circ \gamma(t)^{-1} \circ \gamma(t)(\alpha(s)) 
= \delta^l(\gamma)(t)(H(t,s)) $$ 
and ${\partial H \over \partial s}(t,s) 
= \gamma(t).\alpha'(s).$
We therefore obtain with Remark IX.10(a) the formula 
$$ \eqalign{ 
\int_\alpha F_\omega(\tilde\gamma(1)) 
&= \int_\alpha \int_0^1
[\gamma(t)^*.i_{\delta^l(\gamma)(t)}\omega]\, dt \cr
&= \int_0^1 \int_0^1 \omega_{\gamma(t).\alpha(s)}(\delta^l(\gamma)(t)(\gamma(t).\alpha(s)),
\gamma(t).\alpha'(s))\, dt\, ds \cr
&= \int_0^1 \int_0^1 \omega_{H(t,s)}\Big({\partial H(t,s) \over \partial t}(t,s), 
{\partial H(t,s) \over \partial s}(t,s)\Big)\, dt\, ds \cr
&=  \int_{[0,1]^2} H^*\omega =  \int_H \omega. \cr} $$
\qed

The preceding proposition justifies the term `flux cocycle' because it says
that $\int_\alpha F_\omega(\tilde\gamma(1))$ measures the
`$\omega$-surface area' of the surface obtained by moving the loop
$\alpha$ by the curve $\gamma$ in $\Diff(M)$. 

\Corollary IX.12. If $\gamma(1) = \gamma(0) = \id_M$, then 
$F_\omega(\tilde\gamma(1)) \in H^1_{\rm dR}(M,V)$, 
and we obtain a homomorphism 
$$ F_\omega\res_{\pi_1(\Diff(M))} \: \pi_1(\Diff(M)) \to H^1_{\rm dR}(M,V). $$

\Proof. We keep the notation from Proposition IX.11. If 
the curve $\gamma$ in $\Diff(M)$ is closed and 
$\tilde \gamma$ is the corresponding map $\SS^1 \to \Diff(M)$, then 
$H$ induces a continuous map $\tilde H \: \T^2 \to M, (t,s) \mapsto
\tilde\gamma(t).\alpha(s)$ and 
$$ \int_\alpha F_\omega(\tilde\gamma(1)) 
= \int_H \omega = \int_{\tilde H} \omega = \tilde H^*[\omega] \in H^2(\T^2, V) \cong V. $$
As homotopic curves $\alpha_1$ and $\alpha_2$ lead to homotopic maps 
$\tilde H_1, \tilde H_2 \: \T^2 \to M$, we obtain 
$$ \int_{\alpha_1} F_\omega(\tilde\gamma(1)) 
= \int_{\alpha_2} F_\omega(\tilde\gamma(1)) $$
whenever $\alpha_1$ and $\alpha_2$ are homotopic, and this implies
that $F_\omega(\tilde\gamma(1)) \in H^1_{\rm dR}(M,V)$. 

That the restriction of $F_\omega$ to $\pi_1(\Diff(M))$ is a
homomorphism follows from the cocycle property of $F_\omega$ and the fact that
$\pi_1(\Diff(M)) = \ker q_G$ acts trivially on $\hat H^1_{\rm dR}(M,V)$. 
\qed

Let $\omega_M \in
\Omega^2(M,V)$ be a closed $2$-form and identify it with a Lie algebra $2$-cocycle 
$\omega_\g \in Z^2_c(\g,\a)$ for $\g = {\cal V}(M)$ and $\a =
C^\infty(M,V)$. Next we show that the flux cocycle 
$$ F_{\omega_\g} \: \tilde G \to \hat H^1_c(\g,\a) $$
coincides with flux cocycle $F_{\omega_M}$ from Definition~IX.9. For that we 
recall  from 
Lemma~IX.4 that we can view $H^1_{\rm dR}(M,V)$ as a subspace of
$H^1_c(\g,\a)$ because $B^1_c(\g,\a) = dC^\infty(M,V)$, which leads to an embedding 
$$ \hat H^1_{\rm dR}(M,V) \into \hat H^1_c(\g,\a) :=
C^1_c(\g,\a)/B^1_c(\g,\a). $$

\Lemma IX.13. For a closed $2$-form $\omega_M \in \Omega^2(M,V)$ we have
$$ F_{\omega_\g}  = F_{\omega_M} \: \tilde G \to 
\hat H^1_{\rm dR}(M,V)\subeq  \hat H^1_c(\g,\a). $$

\Proof. We parametrize $\SS^1 \cong \R/\Z$ by the unit interval $[0,1]$. Then we
have 
for any smooth curve $\gamma \: [0,1] \to G = \Diff(M)_0^{\rm op}$ starting in
$\1$ and $X \in \g = {\cal V}(M)$: 
$$ \eqalign{ 
I_\gamma(X) 
&:= \int_\gamma i_{X_r}\omega_\g^{\rm eq} 
= \int_0^1 \omega_\g^{\rm eq}(X \gamma(t), \gamma'(t))\, dt \cr
&= \int_0^1 \gamma(t).\omega_\g(\Ad(\gamma(t))^{-1}.X,
\gamma(t)^{-1}\gamma'(t))\, dt \cr
&= \int_0^1 \gamma(t).\omega_M(\gamma(t).(X \circ \gamma(t)^{-1}),
\delta^l(\gamma)(t))\, dt \cr
&= \int_0^1 \omega_M\big(\gamma(t).(X \circ \gamma(t)^{-1}), \delta^l(\gamma)(t)\big)
\circ \gamma(t)\, dt. \cr} $$
From this formula it is easy to see that 
$I_\gamma \in \Lin(\g,\a)$ 
defines a $1$-form on $M$ whose value in $v \in T_p(M)$ is given by 
$$ \eqalign{I_\gamma(v) 
&= \int_0^1 (\omega_M)_{\gamma(t).p}(\gamma(t).v,
\delta^l(\gamma)(t)(\gamma(t).p))\, dt.  \cr} $$
This means that 
$$ I_\gamma = -\int_0^1 \gamma(t)^* \big(i_{\delta^l(\gamma)(t)} \omega_M\big)\, dt, $$
which, in view of Remark~IX.10, implies that 
$$ F_{\omega_\g}(\tilde\gamma(1)) = [-I_\gamma] =
F_{\omega_M}(\tilde\gamma(1)) \in \hat H^1_{\rm dR}(M,V)
\subeq \hat H^1_c(\g,\a).$$
The remaining assertions now follow from Corollary~IX.12. 
\qed

\Corollary IX.14. $F_\omega(\pi_1(G))$ vanishes if and only if for each
smooth loop $\alpha \: \SS^1 \to M$ and each smooth loop $\gamma \:
\SS^1 \to \Diff(M)$ we have $\int_H \omega = 0$ for the map 
$H \: \T^2 \to M, H(t,s) = \gamma(t).\alpha(s).$
\qed

The condition in the preceding corollary is in particular satisfied if
the set of homotopy classes of based maps $\T^2 \to M$ or at least the
corresponding homology classes in $H_2(M)$ are trivial. 

\Remark IX.15. It is interesting to observe that the discreteness of the period
map for $\omega \in \Omega^2(M,V)$ leads to a condition on the group
of spherical cycles, i.e., the image of $\pi_2(M)$ in $H_2(M)$, and
the vanishing of $F_\omega(\pi_1(G))$ leads to a condition on the larger subgroup of
$H_2(M)$ generated by the cycles coming from maps $\T^2 \to M$.
That the latter group contains the former follows from the existence 
of a map $\T^2 \to \SS^2$ inducing an isomorphism $H_2(\T^2) \to
H_2(\SS^2)$. 
\qed


\subheadline{Examples} 

\Example IX.16. Let $\z$ be a Fr\'echet space, $\Gamma_Z \subeq \z$ a discrete
subgroup,  $Z := \z/\Gamma_Z$ and $q_Z \: \z\to Z$ the quotient map,
which can also be considered as the exponential map of the Lie group $Z$. 

Further let $q \: P \to M$ be a
smooth $Z$-principal bundle over the compact manifold $M$, 
$\theta \in \Omega^1(P,\z)$ a principal connection $1$-form and
$\omega \in \Omega^2(M,\z)$ the corresponding curvature, i.e., 
$q^*\omega = -d\theta$. We call a vector field $X \in {\cal V}(P)$
{\it horizontal} if $\theta(X) = 0$. Write ${\cal V}(P)^Z$ for the Lie algebra 
of $Z$-invariant vector fields on $P$. Then we have an isomorphism 
$$ \sigma \:  {\cal V}(M) \to {\cal V}(P)_{\rm hor}^Z := \{ X \in
{\cal V}(P)^Z \: \theta(X) = 0 \} $$
which is uniquely determined by $q_* \sigma(X) = X$ for $X \in {\cal
V}(M)$. For two horizontal vector fields $\tilde X, \tilde Y$ on $P$ we then have 
$$ (q^*\omega)(\tilde X,\tilde Y) 
= -d\theta(\tilde X,\tilde Y) = \tilde Y.\theta(\tilde X) - \tilde
X.\theta(\tilde Y) -\theta([\tilde Y,\tilde X]) = \theta([\tilde
X,\tilde Y]). $$
This means that 
$$\omega(X,Y) = (q^*\omega)(\sigma(X), \sigma(Y)) = 
\theta([\sigma(X), \sigma(Y)]) = \theta([\sigma(X),
\sigma(Y)]-\sigma([X,Y])) \leqno(9.3) $$ 
can be viewed as the cocycle of the abelian extension 
$$ \a := \gau(P) \cong C^\infty(M,\z) \into \hat\g := {\cal V}(P)^Z 
\onto \g = {\cal V}(M) $$
with respect to the section $\sigma \: \g \to \hat\g$. 

On the group level we find that the inverse image $\hat G$ of $G =
\Diff(M)_0^{\rm op}$ in $\Aut(P)^{\rm op}$ is an
extension of $G$ by the abelian gauge group $A := \Gau(P) \cong C^\infty(M,Z)$ 
and we have already seen above that its Lie algebra is 
$\hat\g \cong \a \oplus_\omega \g$.

The exponential function of the abelian Lie group $A$ is given by
$$ \exp_A \: \a = C^\infty(M,\z) \to C^\infty(M,Z), \quad \xi \mapsto
q_Z \circ \xi. $$
Its image is the identity component $A_0$ of $A$. 
The characteristic map 
$$ \oline\theta_A \: \pi_0(A) \to H^1_c(\g,\a), \quad [f] \mapsto [D(d_G f)] $$
considered in Proposition VI.4 can be made more explicit by observing
that 
$$ (d_G f)(g) = g.f - f = f \circ g - f, $$
so that 
$$ D(d_G f)(X) = X.f = \la df, X \ra $$
(cf.\ Definition A.2). This means that $D(d_G f)$ can be identified with the $1$-form $df \in
H^1_{\rm dR}(M,\z) \subeq H^1_c(\g,\a)$. Therefore the homomorphism 
$\oline\theta_A \: \pi_0(A) \to H^1_c(\g,\a)$ from Proposition VI.4 is obtained
by factorization of the map 
$$ A = C^\infty(M,Z) \to H^1_{\rm dR}(M,\z), \quad f \mapsto [df] $$
whose kernel is the identity component 
$A_0 = q_Z \circ C^\infty(M,\z)$ of $A$ to the injective homomorphism 
$$ \pi_0(A) \cong C^\infty(M,Z)/q_Z \circ C^\infty(M,\z) 
\to H^1_{\rm dR}(M,\z), \quad [f] \mapsto [df]. $$
According to [Ne02, Prop.~3.9], its image consists of the subspace 
$$ H^1_{\rm dR}(M,\Gamma_Z) := \Big\{[\alpha] \in H^1_{\rm dR}(M,\z) \: (\forall \gamma \in
C^\infty(\SS^1,M))\, \int_\gamma \alpha \in \Gamma_Z\Big\}, $$
so that 
$$ \oline\theta_A \: \pi_0(A) \to H^1_{\rm dR}(M,\Gamma_Z), \quad [f] \mapsto [df] $$
is an isomorphism. 

In view of Proposition VI.3, the flux homomorphism satisfies 
$F_\omega = -\oline\theta_A \circ 
\delta$, where $\delta \: \pi_1(G) \to \pi_0(A)$ is the connecting
homomorphism corresponding to the long exact homotopy sequence of the
$A$-bundle $\hat G  \to G$. As $\oline\theta_A$ is an isomorphism, 
$F_\omega$ is essentially the same as $\delta$, and we can view it
as a homomorphism 
$$ F_\omega \: \pi_1(G) \to H^1_{\rm dR}(M,\Gamma_Z) \subeq H^1_{\rm dR}(M,\z). $$

Note that we cannot expect $F_\omega(\pi_1(G))$ to vanish because 
the abelian extension $A \into \hat G \to G$ is not an extension by a
connected group. 
\qed

\Example IX.17. (a) We consider the special case where the manifold $M$ is
a torus: $M = T = \t/\Gamma_T$, where $\t$ is a finite-dimensional
vector space and $\Gamma_T \subeq \t$ is a discrete subgroup for which
$\t/\Gamma_T$ is compact. 

Then the group $T$ acts by multiplication maps on itself, and we
obtain a homomorphism $T \into G = \Diff(M)_0^{\rm op}$ which induces a
homomorphism
$$ \eta_T \: \pi_1(T) \to \pi_1(G). $$

Let $\omega_T \in \Omega^2(T,\z)$ be an invariant $\z$-valued $2$-form
on $T$ and $\omega = \Omega_\1 \in Z^2_c(\t,\z)$. 
Then $\omega_T$ is closed because $T$ is abelian. 
If $e_1,\ldots, e_n$ is an integral basis of $\Gamma_T$, then 
the maps 
$$ \T^2 \to T, \quad (t,s) \mapsto t e_i + s e_j + \Gamma_T, \quad i <
j$$
lead to an integral basis of $H_2(T) \cong \Z^{{\dim T \choose 2}}$,
so that the period group of $\omega_T$ is 
$$ \Gamma_\omega := \span_\Z \omega(e_i, e_j) 
= \span_\Z \omega(\Gamma_T, \Gamma_T) \subeq \z. $$
We assume that $\Gamma_Z \subeq \z$ is a discrete subgroup
with 
$$ \omega(\Gamma_T, \Gamma_T) \subeq \Gamma_Z $$
and put $Z := \z/\Gamma_Z$. 

In view of 
$\pi_2(T) = \{0\}$, we have $\per_{\omega} = 0$ by
Proposition~IX.5. Next we are making the map 
$$ F_\omega \circ \eta_T \: \pi_1(T) = \Gamma_T \to 
H^1_{\rm dR}(T,\Gamma_Z) 
\cong \Hom(\Gamma_T, \Gamma_Z) $$
more explicit. 
For $x, y \in \Gamma_T$ and the corresponding loops 
$\gamma_x(t) = tx + \Gamma_T$ and $\gamma_y(t) = ty + \Gamma_T$ in $T$
we have for 
$$H \: \T^2 \to T, (t,s) \mapsto \gamma_x(t) + \gamma_y(s) 
= [tx + s y] $$ 
the formula 
$$ \int_{\gamma_y} F_\omega([\gamma_x])
= \int_H \omega = \omega(x,y) $$
(Proposition IX.11, Lemma IX.13). 
This means that $F_\omega \circ \eta_T \: \pi_1(T) \to
\Hom(\pi_1(T),\Gamma_Z)$ can be identified with the map 
$x \mapsto i_x\omega$. 

If $\omega \not=0$, then $F_\omega(\pi_1(G))\not=\{0\}$, which means that there
is no abelian extension $A \into \hat T \onto T$ with a connected
abelian group $A$ of the form $\a/\Gamma_A$ for $\a = C^\infty(T,\z)$. 
Another reason for this is that any such extension would be central, 
but all central extensions of tori by connected Lie groups are flat in
the sense that their Lie algebra cocycle vanishes (cf.\ [Ne02]). 

On the other hand, the existence of a $Z$-bundle over $T$ with
curvature $\omega$ implies the existence of an abelian extension 
$$ A := C^\infty(T,Z) \into \hat T \onto T, $$
 where $T$ acts on $A$ by $(t.f)(x) = f(x + t)$ (cf.\ Example IX.16). 
The corresponding Lie
 algebra cocycle 
$\omega \in Z^2_c(\t,C^\infty(T,\z))$ is given by 
$(x,y) \mapsto \omega(x,y) \in \z$ whose values lie in 
$\z \cong \a^T$. 

\par\nin (b) Let $\t$ be a locally convex space, $\Gamma_T \subeq \t$
 a discrete subgroup and consider the connected abelian Lie group $T
 := \t/\Gamma_T$. 
Let further $\z$ be a sequentially complete locally convex space, 
$\Gamma_Z \subeq \z$ be a discrete subgroup and $Z := \z/\Gamma_Z$,
 considered as a trivial $T$-module. 
We fix an alternating continuous map $\omega \in Z^2_c(\t,\z)$ and define 
$\omega_Z \in Z^2_s(\t,Z)$ by $f_Z := q_Z \circ {1\over 2}\omega$, where 
$q_Z \: \z \to Z$ is the quotient map. 

Let $H := Z \times_{f_Z} \t$ denote the corresponding central
extension of $\t$ by $Z$. 
Then $Z^\sharp := Z \times_{f_Z} \Gamma_T$ is a normal subgroup of $H$ because 
all commutators lie in $Z$. Since $H/Z^\sharp  \cong
\t/\Gamma_T = T$, we can think of $H$ as an extension 
$$ Z^\sharp \into H \onto T. $$

Since $Z$ is divisible and $\Gamma_T$ discrete, the 
central extension $Z \into Z^\sharp \onto \Gamma_T$ is trivial if and
only if it is an abelian group, which means that its commutator map 
$\Gamma_T \times \Gamma_T \to Z$
vanishes. The commutator map is given by 
$$ \eqalign{ (z,t)(z',t') (z,t)^{-1}(z',t')^{-1} 
&= (f_Z(t,t'),t+t')(f_Z(t',t),t+t')^{-1} \cr
&= (f_Z(t,t')-f_Z(t',t),0) 
= (2f_Z(t,t'),0) = ( q_Z(\omega(t,t')),0).\cr} $$
Therefore $Z^\sharp$ is a trivial extension of $\Gamma_T$ if and only
if 
$$ \omega(\Gamma_T, \Gamma_T) \subeq \Gamma_Z. \leqno(9.4) $$

The condition for the extistence of a $Z$-bundle $P \to T$ with
curvature $\omega_T$ is also given by (9.4). 
The necessity of this condition in the infinite-dimensional case can
be seen by restricting to two-dimensional subtori. 
If (9.4) is satisfied, then we can view $\Gamma_T$ as a subgroup of
$Z^\sharp$ because there exists a homomorphism $\sigma \: \Gamma_T \to Z^\sharp$
splitting the extension $Z^\sharp \onto \Gamma_T$. 
Now we form the homogeneous space 
$P := H/\sigma(\Gamma_T)$ which defines a $Z$-bundle 
$$ Z \into P = H/\sigma(\Gamma_T) \onto T \cong H/Z^\sharp. $$

As $Z$ is central in $H$, the left action of $H$ on $P$ induces a
homomorphism 
$$ H \to \Aut(P) = \Diff(P)^Z $$
restricting to a homomorphism 
$$ j_Z \: Z^\sharp \cong Z \times_{f_Z} \Gamma_T\to \Gau(P) \cong C^\infty(T,Z), $$
where the elements of $Z$ correspond to constant functions. The group
$\Gamma_T$ acts on $P$ by 
$$ x.(q_Z(z),y) 
= (q_Z(z + {\textstyle{1\over 2}}\omega(x,y)), y) 
= (q_Z(z),y).f_Z(x,y), $$
so that 
$$ j_Z(z,x)(y + \Gamma_T) = z + f_Z(x,y). $$

If $\omega(\Gamma_T,\t) \not=\{0\}$, then the map 
$$F_\omega \: \pi_1(T) \cong\Gamma_T\to H^1_c(\t,\z) = \Lin(\t,\z), 
\quad F_\omega(x)(y) = \omega(y,x) $$
does not vanish, but if $\omega(\Gamma_T, \Gamma_T) \subeq \Gamma_Z$,
then the extension $Z \into Z^\sharp \onto \Gamma_T$ is
trivial. Therefore the natural sequence 
$$ H^2_s(T,Z) \to H^2_s(\tilde T,Z) \cong H^2_c(\t,\z) \to
H^2(\Gamma_T,Z) \leqno(9.5) $$
is {\sl not} exact in $H^2_s(\tilde T, Z)$ (cf.\ Theorem VII.2). 

Identifying 
$$ H^1_{\rm dR}(T,\Gamma_Z) \cong d C^\infty(T,Z)/dC^\infty(T,\z)
\subeq H^1_{\rm dR}(T,\z)$$ 
with a subspace of 
$H^1_c(\t,\a)$ (cf.\ Lemma IX.4), we can view 
$F_\omega$ as a map 
$$ \pi_1(T) \to H^1_{\rm dR}(T,\Gamma_Z) \into \Hom(\pi_1(T),
\Gamma_Z). 
\qeddis 

\sectionheadline{X. The diffeomorphism group of the circle} 

\nin In this section we apply the general results from Sections VI and VII to the 
group of orientation preserving diffeomorphisms of the circle $\SS^1$ and 
the modules ${\cal F}_\lambda$ of $\lambda$-densities on $\SS^1$ whose cohomology for the group $\Diff(\SS^1)_0$ has been determined in [OR98]. 
We shall also point out how the picture changes if $\Diff(\SS^1)_0$ 
is replaced by its universal covering group. 

Let $G := \Diff(\SS^1)_0^{\rm op}$ be the group of orientation
preserving diffeomorphisms of the circle $\SS^1\cong\R/\Z$. Then 
its universal covering group $\tilde G$ can be
identified with the group 
$$ \tilde G := \{ f \in \Diff(\R)^{\rm op} \: (\forall x\in \R)\, f(x + 1) =
f(x) + 1\}, $$
and the covering homomorphism $q_G \: \tilde G \to G$ is given by 
$q(f)([x]) = [f(x)],$
where $[x] = x + \Z \in \SS^1 \cong \R/\Z$. The kernel of $q_G$ 
consists of all translations $\tau_a$, $a \in \Z$, and since 
$\tilde G$ is an open  convex subset of a closed subspace of 
$C^\infty(\R,\R)$, it is a contractible manifold. In particular, we obtain 
$$ \pi_1(G) \cong \Z \quad \hbox{ and } \quad 
\pi_k(G) = \{\1\}, \quad k \not= 1. $$

The group $G$ has an import series of representation 
${\cal F}_\lambda$, $\lambda \in \R$, where 
${\cal F}_\lambda$ is the space of $\lambda$-densities on the circle $\SS^1$. 
As the tangent bundle $T\SS^1$ is trivial, we may identify the space 
${\cal F}_\lambda$ with the space  
$C^\infty(\SS^1,\R)$ of $1$-periodic functions on $\R$ with the representation 
$$ \rho_\lambda(\phi).\xi 
= (\phi')^{\lambda} \cdot (\xi \circ \phi) $$
which corresponds symbolically to 
$\phi^*(\xi (dx)^\lambda) = (\xi \circ \phi)\cdot (\phi')^\lambda
\cdot (dx)^\lambda.$
Note that ${\cal F}_0 = C^\infty(\SS^1,\R)$ is a Fr\'echet algebra and that, 
as $G$-modules, 
$$ {\cal F}_1 \cong \Omega^1(\SS^1,\R) \quad \hbox{ and } \quad 
{\cal F}_{-1} \cong {\cal V}(\SS^1) = \g. $$

For the Lie algebra $\g = {\cal V}(\SS^1)$ of $G$ the derived 
representation is given on $X = \xi {d \over dx}$ by 
$$ \rho_\lambda(\xi).f  = \xi f' + \lambda f \xi'. \leqno(10.1) $$
This follows directly from $\rho_\lambda(g).f = 
(g')^\lambda \cdot (f \circ g)$ and the product rule. In the following we shall 
identify $\g$ with $C^\infty(\SS^1,\R)$ and denote elements of $\g$ by $\xi$, $\eta$ etc.  

\Lemma X.1. On the Fr\'echet--Lie group $A :=
C^\infty(\SS^1,\R^\times) = {\cal F}_0^\times$ we have a smooth $G$-action by 
$g.f := f \circ g$ and the derivative 
$\eta \: G \to A, f \mapsto f'$ 
is a smooth $1$-cocycle. 

\Proof. For $g,h \in G$ we have 
$\eta(gh) = (gh)' = (h \circ g)' = (h' \circ g) \cdot g' 
= (g.\eta(h)) \cdot \eta(g).$ 
\qed 

\Remark X.2. The representation on ${\cal F}_\lambda$ has the form 
$\rho_\lambda(g).f = \eta(g)^\lambda \cdot (f \circ g)$ 
and the fact that $\eta^\lambda \: G \to A$ is a cocycle
implies that $\rho_\lambda \: G \to \GL({\cal F}_\lambda)$ 
is a group homomorphism.
\qed

\subheadline{The cohomology on the Lie algebra level} 

\Proposition X.3. The cohomology in degrees $0,1,2$ of the $\g$-module  
${\cal F}_\lambda$ has the following structure: 
$$H^0_c(\g,{\cal F}_\lambda) = {\cal F}_\lambda^\g 
= \cases{ \{0\} & for $\lambda \not= 0$ \cr 
\R 1 & for $\lambda = 0$. \cr}.  $$
For $n \in \N_0$ let $\alpha_n(\xi) = \xi^{(n)}$ 
denote the $n$-fold derivative. Then 
$$ H^1_c(\g,{\cal F}_0) = \span \{[\alpha_0], [\alpha_1]\}, \quad 
H^1_c(\g,{\cal F}_1) = \R [\alpha_2], \quad 
H^1_c(\g,{\cal F}_2) = \R [\alpha_3] $$
and $H^1_c(\g, {\cal F}_\lambda)$ vanishes for $\lambda\not=0,1,2$. 
In degree $2$ we have 
$$H^2_c(\g,{\cal F}_\lambda) 
\cong \cases{ \R^2 & for $\lambda = 0,1,2$ \cr 
\R  & for $\lambda = 5,7$ \cr
\{0\}  & otherwise. \cr} $$
For $\lambda = 0,1,2$ the cohomology classes of the following elements
form a basis of $H^2_c(\g, {\cal F}_\lambda)$: 
$$ \oline\omega_0(\xi,\eta) := \left|\matrix{ \xi & \eta \cr \xi' & \eta'\cr}\right|,
\quad \omega_0(\xi,\eta) := \int_0^1 \left|\matrix{ \xi' & \eta' \cr \xi'' &
\eta''\cr}\right| \quad \hbox{ for } \quad \lambda =0, $$
$$ \oline\omega_1(\xi,\eta) := \left|\matrix{ \xi & \eta \cr \xi'' &
\eta''\cr}\right|, \quad 
\omega_1(\xi,\eta) := \left|\matrix{ \xi' & \eta' \cr \xi'' & \eta''\cr}\right|
\quad \hbox{ for } \quad \lambda =1, $$
and 
$$ \oline\omega_2(\xi,\eta) := \left|\matrix{ \xi & \eta \cr \xi''' &
\eta'''\cr}\right|, \quad 
\omega_2(\xi,\eta) := \left|\matrix{ \xi' & \eta' \cr \xi''' & \eta'''\cr}\right|  
\quad \hbox{ for } \lambda = 2. $$

\Proof. (cf.\ [OR98]) We have 
$$ {\cal F}_\lambda^\g 
\cong \{ f \in C^\infty(\SS^1,\R) \: (\forall \xi \in C^\infty(\SS^1,\R)) 
\xi f' + \lambda \xi' f = 0\}. $$
For constant functions $\xi$ the differential equation 
from above reduces to $f'\xi = 0$, so that $f$ is constant, and now 
$\lambda \xi' f= 0$ for each $\xi$ implies 
$\lambda f = 0$. This proves the assertion about $H_c^0(\g,{\cal F}_\lambda)$. 

According to [Fu86, p.176], we have 
$$ H^q_c(\g, {\cal F}_\lambda) = 0 
\quad \hbox{ for } \quad 
\lambda \not\in \Big\{ {3 r^2 \pm r \over 2} \: r \in \N_0\Big\} 
= \{ 0,1,2,5,7,12,15, \ldots \}. $$
If $r \in \N_0$ and $\lambda = {3 r^2 \pm r \over 2}$, then 
$$ H^q_c(\g, {\cal F}_\lambda) \cong 
\cases{ 
H_{\rm sing}^{q-r}(Y(\SS^1), \R) & for $q \geq r$ \cr 
\{0\} & for $q < r$, \cr } $$
where $Y(\SS^1) = \T^2\times \Omega \SS^3$ and 
$\Omega \SS^3$ is the loop space of $\SS^3$. 
The cohomology algebra 
$$H^*_c(\g, {\cal F}_0) \cong H^*_{\rm sing}(Y(\SS^1), \R) 
\cong H^*_{\rm sing}(\SS^1, \R) \otimes H^*_{\rm sing}(\SS^1, \R) \otimes 
H^*_{\rm sing}(\Omega \SS^3, \R) $$
is a free anticommutative real algebra with 
generators $a,b,c$ satisfying 
$$ \deg(a) = \deg(b) = 1, \quad \deg(c) =2, \quad 
a^2 = b^2 = 0. $$
It follows in particular that 
$$ H^0_c(\g, {\cal F}_0) = \R,\quad 
 H^1_c(\g, {\cal F}_0) = \R a + \R b \cong \R^2, \quad 
 H^2_c(\g, {\cal F}_0) = \R c + \R ab \cong \R^2. $$ 
The structure of $H^*_c(\g, {\cal F}_\lambda)$ is now determined by 
the fact that it is a free module of the algebra 
$H^*(Y(\SS^1),\R)\cong H^*_c(\g,{\cal F}_0)$ with one generator in 
degree $r$. Here the algebra structure on $H^*_c(\g,{\cal F}_0)$ is obtained 
from the multiplication on ${\cal F}_0$ as in Appendix F, and the multiplication 
${\cal F}_0 \times {\cal F}_\lambda \to {\cal F}_\lambda$ yields the 
$H^*_c(\g, {\cal F}_0)$-module structure 
$([\alpha], [\beta]) \mapsto [\alpha \wedge \beta]$ on $H^*_c(\g, {\cal F}_\lambda)$. 

From [Fu86, Th.~2.4.12] we see that generators of 
$H^*_c(\g, {\cal F}_0)$ are given by the classes of 
$\alpha_0, \alpha_1$ and $\omega_0$. Therefore a second basis element of 
$H^2_c(\g, {\cal F}_0)$ is represented by 
$$ (\alpha_0 \wedge \alpha_1)(\xi,\eta) = \alpha_0(\xi)\alpha_1(\eta) 
- \alpha_0(\eta)\alpha_1(\xi) = \xi \eta' - \xi' \eta = \oline\omega_0(\xi,\eta). $$

The space $H^1(\g, {\cal F}_\lambda)$ is non-zero for 
$r = 0,1$ which corresponds to $\lambda \in \{ 0,1,2\}$. 
For $r = 0$ it is two-dimensional and for $r = 1$ it is one-dimensional. 
For $\lambda = 1$ a generator is given by 
$[\alpha_2]$ ([Fu86, Th.~2.4.12]; there is a misprint in the formula!). 
From the $H^*_c(\g, {\cal F}_0)$-module structure of 
$H^*_c(\g, {\cal F}_1)$ we obtain the generators of $H^2_c(\g, {\cal F}_1)$:
$$ (\alpha_0 \wedge \alpha_2)(\xi, \eta) = \xi \eta'' - \eta \xi'' = \oline{\omega}_1, 
\quad (\alpha_1 \wedge \alpha_2)(\xi, \eta) = \xi' \eta'' - \eta' \xi'' = {\omega}_1. $$

Averaging over the rotation group, we see that every cocycle is equivalent to a 
rotation invariant one. From that it is easy to verify that for 
$\lambda=2$ a generator of $H^1_c(\g, {\cal F}_2)$ 
is given by $[\alpha_3]$, and we obtain for the basis elements of 
$H^2_c(\g, {\cal F}_2)$: 
$$ (\alpha_0 \wedge \alpha_3)(\xi, \eta) = \xi \eta''' - \eta \xi''' = \oline{\omega}_2, 
\quad (\alpha_1 \wedge \alpha_3)(\xi, \eta) = \xi' \eta''' - \eta' \xi''' = {\omega}_2. 
\qeddis 

For an explicit description of a basis of $H^2_c(\g, {\cal F}_\lambda)$ 
for $\lambda = 5,7$ we refer to [OR98]. 

\subheadline{Integrating Lie algebra cocycles to group cocycles} 

Now we translate the information on the Lie algebra cohomology 
$H^p_c(\g, {\cal F}_\lambda)$ for $p=0,1,2$ (Proposition X.3) to the group $G$. 
Since the group $G$ is connected, we have 
$$ H^0_s(G, {\cal F}_\lambda) = {\cal F}_\lambda^G = {\cal F}_\lambda^\g = \cases{ 
\{0\} & for $\lambda \not= 0$ \cr 
\R 1 & for $\lambda = 0$. \cr} $$

\msk 

In degree $1$, we can use Proposition III.4 to see that we have an exact sequence 
$$ \0 \to H^1_s(G, {\cal F}_\lambda) \sssmapright{D} 
H^1_c(\g, {\cal F}_\lambda) \sssmapright{P} 
{\cal F}_\lambda^\g. $$
For $\lambda \not=0$ this implies that 
$D \: H^1_s(G, {\cal F}_\lambda) \to H^1_c(\g, {\cal F}_\lambda)$
is an isomorphism.  For $\lambda = 0$ we have to calculate the period map 
$P$. 
Let $\t := \R \1 \cong \R {d\over dx} \subeq \g$ 
be the one-dimensional subalgebra corresponding to the rotations of the 
circle $\SS^1$ and $T \cong \T \subeq G$ the corresponding subgroup. 
Then the inclusion 
$T \into G$ induces an isomorphism $\pi_1(T) \to \pi_1(G)$, 
so that we can calculate $P$ by  restricting to $T$. 
Since $\t$ corresponds to constant functions, the cocycle $\alpha_1$ vanishes on 
$\t$, and the cocycle $\alpha_0$ is non-trivial on $\t$. Hence 
$$ H^1_s(G, {\cal F}_0) \cong \ker P = \R [\alpha_1]. $$
The group cocycle corresponding to $\alpha_1(\xi) = \xi'$ is 
$\theta(\phi) = \log\phi'$ (cf.\ Lemma X.1) 
because for $\phi = \id_\R + \xi$ we have 
$$ \theta(\id + \xi) = \log(1 + \xi') \sim \xi' + \ldots, $$
which implies $D \theta = \alpha_1$.
Since the map $d \: {\cal F}_0 \cong C^\infty(\SS^1, \R) \to {\cal F}_1 \cong 
\Omega^1(\SS^1, \R)$ is equivariant, we obtain a group cocycle 
$$ d\circ \theta \in Z^1_s(G, {\cal F}_1), 
\quad (d\circ \theta)(f) := {\log(f')'} = {f'' \over f'}, $$
and for $\phi = \id + \xi$ the relation 
$(d \circ \theta)(\id + \xi) = {\xi'' \over 1 + \xi'}$
directly leads to $D(d \circ \theta) = \alpha_2.$
The Schwarzian derivative 
$$ S \in Z^1_s(G, {\cal F}_2), \quad 
S(\phi) := \Big({\phi'''\over \phi'} 
- {3\over 2}\Big({\phi'' \over \phi'}\Big)^2\Big) $$
satisfies $D S = \alpha_3$. We thus have 
$$ H^1_s(G, {\cal F}_\lambda) = \cases{ \{0 \} & for $\lambda \not= 0,1,2$ \cr 
\R [\theta]  & for $\lambda = 0$ \cr 
\R [d \circ \theta]  & for $\lambda = 1$ \cr 
\R [S] & for $\lambda = 2$. \cr } $$

On the simply connected covering group $q_G \: \tilde G \to G$ we have 
$H^1_s(\tilde G, {\cal F}_\lambda) \cong H^1_c(\g, {\cal F}_\lambda)$ 
(Proposition~III.4), so that we need 
an additional $1$-cocycle for $\lambda = 0$, which is given by 
$$ L(\phi) := \phi - \id_\R. $$
In fact, 
$L(\psi\phi) =  L(\phi\circ \psi) := \phi \circ \psi - \psi + \psi - \id_\R 
= \psi^*L(\phi) + L(\psi).$
Since $DL = \alpha_0$, we get 
$$ H^1_s(\tilde G,{\cal F}_0) = \R [L] + \R [\theta], $$
where $\theta(\phi) = \log \phi'$.

\msk 
Now we turn to the group cohomology in degree $2$: 
In view of $\pi_1(G) \cong \Z$ and Theorem~VII.2, we have a map 
$$ \delta \: \Hom(\pi_1(G), {\cal F}_\lambda^G) \cong {\cal F}_\lambda^G 
\to H^2_s(G,{\cal F}_\lambda),
\quad \delta(\gamma) = ({\cal F}_\lambda \rtimes \tilde G)/\Gamma(\gamma). $$
The kernel of this map coincides with the image of the restriction map 
$$ R \: H^1_s(\tilde G, {\cal F}_\lambda) \cong H^1_c(\g, {\cal F}_\lambda) \to 
\Hom(\pi_1(G), {\cal F}_\lambda^G) \cong {\cal F}_\lambda^G  $$
and the image of $D$ coincides with the kernel of the map 
$$ P \: H^2_c(\g, {\cal F}_\lambda) \to 
\Hom(\pi_1(G), H^1_c(\g, {\cal F}_\lambda)) 
\cong H^1_c(\g, {\cal F}_\lambda). $$

The following proposition clarifies the relation between second Lie algebra 
and Lie group cohomology for the modules ${\cal F}_\lambda$. 
We refer to Appendix F for the definition of the $\cap$-product of Lie group 
cocycles. 

\Proposition X.4. For each $\lambda \in \R$ the map 
$D \: H^2_s(G,{\cal F}_\lambda) \to H^2_c(\g,{\cal F}_\lambda)$
is injective.  It is bijective for $\lambda \not\in \{0,1,2\}$. 
For $\lambda \in \{0,1,2\}$ we have
$$ H^2_s(G, {\cal F}_0) = \R [B_0], \quad 
 H^2_s(G, {\cal F}_1) = \R [B_1], \quad 
 H^2_s(G, {\cal F}_2) = \R [B_2] $$
for 
$$ B_0(\phi, \psi) := -\int_0^1 \log((\psi \circ \phi)')d(\log \phi'), \quad 
B_1  := \theta \cap (d \circ \theta) \quad \hbox{ and } \quad 
B_2  := \theta \cap S. $$

\Proof. (cf.\ [OR98]) First we show that $D$ is injective for each $\lambda$. 
As above, let $T \cong \T \subeq G$ be the subgroup corresponding to 
$\t = \R 1$ in $\g$. 
Since the inclusion 
$T \into G$ induces an isomorphism $\pi_1(T) \to \pi_1(G)$, 
we can calculate $R$ by  
using the factorization 
$$ H^1_s(\g, {\cal F}_0) \to H^1_s(\t, {\cal F}_0) 
\to \Hom(\pi_1(T), {\cal F}_0^G) \cong {\cal F}_0^G 
\cong \Hom(\pi_1(G), {\cal F}_0^G). $$
It is clear that the cocycle $\alpha_1$ vanishes on $\t$, but 
$\alpha_0$ satisfies 
$\per_{\alpha_0}([\id_T]) = 1 \in {\cal F}_0^G.$
Therefore 
the restriction map $R$ is surjective for $\lambda = 0$, which implies $\delta = 0$. 
For all other values of 
$\lambda$ the map $\delta$ vanishes because ${\cal F}_\lambda^G$ is
trivial. Therefore $D$ is injective for each $\lambda$. 

For $\lambda \not\in \{0,1,2\}$ the space $H^1_c(\g, {\cal F}_\lambda)$ vanishes, 
so that $P  = 0$ and $\im(D) = \ker(P)$ imply that $D$ is surjective. 

For $\lambda = 0,1,2$ the space $H^2_c(\g, {\cal F}_\lambda)$ is two-dimensional 
(Proposition~X.3). To calculate $P$ in these cases, let 
$$ \gamma \: [0,1] \to T \subeq G, \quad t \mapsto (x \mapsto x + t + \Z)  $$
be the generator of $\pi_1(G)$. We have 
$$ I_\gamma(x) 
= \int_0^1 (i_{x_r}.\omega^{\rm eq})(\gamma'(t))\, dt 
= \int_0^1 \gamma(t).\omega\big(\Ad(\gamma(t))^{-1}.x,1\big)\, dt. $$
This means that $I_\gamma$ is the $T$-equivariant part of the linear map 
$-i_1 \omega \: \g \to {\cal F}_\lambda.$

For the cocycle 
$\oline\omega_\lambda(\xi,\eta) := \xi \eta^{(\lambda+1)} 
- \eta \xi^{(\lambda+1)}$ we have 
$$\big(i_1\oline\omega_\lambda\big)(\eta) = 
\oline\omega_\lambda(1,\eta) = \eta^{(\lambda+1)}. $$
As $1$ acts on each ${\cal F}_\lambda$ by $\xi \mapsto \xi'$ the linear map 
$\oline\omega_\lambda(1,\cdot)$ is $T$-equivariant, hence equal to $I_\gamma$, 
and we obtain 
$$ F_{\oline\omega_\lambda}(1) = -[I_\gamma], \quad I_\gamma(\eta) = -\eta^{(\lambda+1)}, 
\quad \hbox{ for } \quad \lambda = 0,1,2. $$
For 
$\omega_0(\xi,\eta) := \int_{\SS^1} \xi' \eta'' - \xi'' \eta'$ 
we have $\omega_0(1,\eta) = 0$, so that $F_{\omega_0} = 0$, and likewise 
$\omega_\lambda(1,\eta)=0$ for $\lambda = 1,2$ leads to $F_{\omega_\lambda} = 0$ for 
$\lambda = 1,2$. 
 
We conclude that for $\lambda = 0,1,2$ the kernel of $P$ is one-dimensional, and that 
$$ \im(D) = \ker(P) = \R [\omega_\lambda]. $$
For $\lambda = 0$ the Thurston--Bott cocycle (for $\Diff(\SS^1)^{\rm op}$) 
$$ B_0 \in Z^2_s(G, \R) \subeq Z^2_s(G, {\cal F}_0), 
\quad B_0(\phi, \psi) = -\int_0^1 \log((\psi \circ \phi)')d(\log \phi') $$
satisfies $DB_0 = \omega_0$ (cf.\ [GF68]). 
For $\lambda = 1,2$ we reall that $\omega_\lambda = \alpha_1 \wedge \alpha_{\lambda+1}$, 
so that Lemma F.3 implies that  the cocycles 
$$ B_1  := \theta \cap (d \circ \theta) \quad \hbox{ and } \quad 
B_2  := \theta \cap S  $$
satisfy $D B_\lambda = \omega_\lambda$. This completes the proof. 
\qed

\Proposition X.5. For the simply connected covering group $\tilde G$ of $G$ we have 
$$ H^2_s(\tilde  G, {\cal F}_\lambda) = 
\R [B_\lambda] \oplus \R [\oline B_\lambda] \cong \R^2 \quad \hbox{ for } \quad 
\lambda = 0,1,2,$$ 
where 
$$ \oline B_0 := L \cap \theta, \quad 
 \oline B_1 := L \cap (d \circ \theta) \quad \hbox{ and } \quad 
 \oline B_2 := L \cap S $$
and $B_\lambda$ is the pull-back of the corresponding cocycle on $G$. 

\Proof. Since the 
simply connected covering group $\tilde G$ is contractible, the derivation map 
$$ D \: H^2_s(\tilde G,{\cal F}_\lambda) \to H^2_c(\g, {\cal F}_\lambda) $$ 
is bijective, so that we obtain larger cohomology spaces of $\tilde G$ than for 
$G$. For $\lambda = 0,1,2$ we have 
$\oline\omega_\lambda = \alpha_0 \wedge \alpha_{\lambda + 1},$
so that the cocycles $\oline B_j$, $j =0,1,2$, 
satisfy $D \oline B_\lambda = \oline \omega_\lambda$ (Lemma F.3). 
Combining this with the pull-backs of the cocycles $B_\lambda$ from $G$, 
the assertion follows. 
\qed

\subheadline{A non-trivial abelian extension of $\SL_2(\R)$} 

We consider the right action of $\SL_2(\R)$ on the projective line 
$\P_1(\R) = \R \cup \{\infty\}$ by 
$$ x.\pmatrix{ a & b \cr c & d \cr} 
:= \pmatrix{ a & b \cr c & d \cr}^{-1}.x := {dx - b \over -c x + a}. $$
In particular the action of the rotation group $\SO_2(\R)$ is given by 
$$ \pmatrix{ \cos \pi t & -\sin \pi t \cr \sin \pi t & \cos \pi t \cr}.x 
= {\cos \pi t \cdot x - \sin  \pi t
\over \sin \pi t \cdot x + \cos \pi t},  $$
so that 
$$ \pmatrix{ \cos \pi t & -\sin \pi t \cr \sin \pi t & \cos \pi t \cr}.0 = - \tan \pi t $$
and the map 
$t \mapsto \tan \pi t$ induces a diffeomorphism $\R/\Z \to \P_1(\R)$.  
We use this diffeomorphism to identify $\SS^1 = \R/\Z$ with $\P_1(\R)$ and 
to obtain a smooth right action of $\SL_2(\R)$ on $\SS^1$. 
Then $\sL_2(\R)$ is isomorphic to a $3$-dimensional subalgebra of ${\cal V}(\SS^1)$ 
and $\so_2(\R)$ corresponds to $\R 1 = \t$. We put 
$$ U := \pmatrix{ 0 & -1 \cr 1 & 0 \cr} $$
and observe that this element corresponds to the constant function 
${1\over \pi}$. From $\ad U((\ad U)^2 + 4) = 0$ on $\sL_2(\R)$ and the formula for 
commutators in ${\cal V}(\SS^1)$ we therefore derive 
$$ \sL_2(\R) = \span \{ 1, \cos(2\pi t), \sin(2\pi t)\} $$
as a subalgebra of ${\cal V}(\SS^1) \cong C^\infty(\SS^1)$. 
We may therefore pick $H,P \in \sL_2(\R)$ with $[U,H] = - 2P$ and $[U,P] = 2 H$ such that 
$H$ corresponds to the function $\cos(2\pi t)$ and $P$ to the function 
$\sin(2\pi t)$. 

The corresponding group homomorphism 
$$ \sigma \: \SL_2(\R) \to \Diff(\SS^1)^{\rm op}_0 $$
is homotopy equivalent to the twofold covering of $T \cong \SS^1$, hence induces 
an injection 
$$ \pi_1(\sigma) \: \pi_1(\SL_2(\R))\cong \Z \to \pi_1(\Diff(\SS^1)) \cong \Z $$
onto a subgroup of index $2$. 

From the action of $\SL_2(\R)$ on $\SS^1$, we obtain a smooth action 
on the Fr\'echet spaces 
$$ {\cal F}_\lambda := C^\infty(\SS^1, \R), \quad (g.f)(x) := \big(\sigma(g)'\big)^\lambda 
f(x.g). $$
By restriction to the subalgebra $\sL_2(\R) \subeq {\cal V}(\SS^1)$, 
we obtain the $2$-cocycle 
$\omega(\xi,\eta) = \xi' \eta'' - \xi'' \eta'$ in 
$Z^2_c(\sL_2(\R), {\cal F}_1)$. 
Let $\gamma \: I \to \SL_2(\R), t \mapsto \exp(2\pi t U)$ be the canonical generator of 
$\pi_1(\SL_2(\R))$. 
As in the proof of Proposition X.4, it then follows that 
$$ F_{\omega} \: \pi_1(\SL_2(\R)) \to H^1_c(\sL_2(\R), {\cal F}_1) $$
is given by $F_{\omega}([\gamma]) = -[I_\gamma]$, where 
$I_\gamma$ is the $\t$-invariant part of $-2 i_1 \omega = 0$, 
hence $F_\omega = 0$. 

Next we show that $[\omega] \not=0$ in $H^2_c(\sL_2(\R), {\cal F}_1)$. 
If this is not the case, then there exists a linear map 
$\alpha \: \sL_2(\R) \to {\cal F}_1$ with $\omega = d \alpha$. 
Since $\omega$ is $T$-equivariant, we may assume, 
after averaging over the compact group $T$, that $\alpha$ is also $T$-invariant, 
i.e., 
$$\alpha([U,x]) = U.\alpha(x), \quad x \in \sL_2(\R). $$
Now 
$$ 0 = i_U \omega = i_U d \alpha = {\cal L}_U.\alpha - d i_U \alpha = - d i_U \alpha $$
implies 
$$i_U \alpha = \alpha(U) \in Z^0(\sL_2(\R),{\cal F}_1) = {\cal F}_1^{\sL_2(\R)} = \{0\}.$$ 
We now derive from $[H,P] \in \R U$: 
$$ \omega(H,P) = d\alpha(H,P) 
= H.\alpha(P) - P.\alpha(H) - \alpha([H,P]) 
= H.\alpha(P) - P.\alpha(H). $$
Further the equivariance of $\alpha$ implies the existence of 
$a,b \in \R$ with 
$$ \alpha(P) = a \cos(2\pi t) + b \sin(2\pi t) 
\quad \hbox{ and } \quad 
\alpha(H) = {1\over 2}\alpha([U,P]) = {1\over 2} U.\alpha(P) 
= - a \sin(2\pi t) + b \cos(2\pi t). $$
We further have 
$$ H.\alpha(P) 
= \cos(2\pi t).(a \cos(2\pi t) + b \sin(2\pi t)) 
=  (a\cos^2(2\pi t) + b \sin(2\pi t)\cos(2\pi t))' $$
and 
$$ P.\alpha(H) 
= \sin(2\pi t).(-a \sin(2\pi t) + b \cos(2\pi t)) 
=  (- a\sin^2(2\pi t) + b \sin(2\pi t)\cos(2\pi t))', $$
so that 
$$ \omega(H,P) = H.\alpha(P) - P.\alpha(H) = a(\cos^2(2\pi t) + \sin^2(2\pi t))' = a 1' = 0, $$
contradicting 
$$ \omega(H,P)  =\cos(2\pi t)' \sin(2\pi t)''  -\cos(2\pi t)'' \sin(2\pi t)' 
= 8\pi^3 (\sin^3(2\pi t) + \cos^3(2\pi t)) \not= 0. $$
Therefore $[\omega] \not= 0$. Since $F_\omega$ and $\per_\omega$ vanish, and 
$$ H^2_{\rm dR}(\SL_2(\R), {\cal F}_1) \cong H^2_{\rm dR}(\SS^1, {\cal F}_1) = \{0\}, $$
there exists a smooth $2$-cocycle 
$f \in Z^2_s(\SL_2(\R), {\cal F}_1)$ with $Df = \omega$ (Proposition VII.4).
Then the group 
$$ {\cal F}_1 \times_f \SL_2(\R) $$
is a non-trivial abelian extension of $\SL_2(\R)$. 
It is diffeomorphic to the direct product 
vector space $C^\infty(\SS^1,\R) \times \R^3$, 
hence contractible. 

If $V$ is a trivial $\sL_2(\R)$-module, then  
the range of each $2$-cocycle lies in a $3$-dimensional subspace, 
hence is a coboundary, because the corresponding assertion holds for finite-dimensional modules. Therefore all central extensions 
of $\SL_2(\R)$ by abelian Lie groups 
of the form $A = \a/\Gamma_A$ are trivial (Theorem VII.2). 
The preceding example shows that $H^2_c(\sL_2(\R), {\cal F}_1) \not= \{0\}$, 
which provides the non-trivial extension of $\SL_2(\R)$. 

The choice of the cocycle $\omega$ above is most natural because one can show that the 
cohomology of the $\sL_2(\R)$-modules 
${\cal F}_\lambda$ satisfies 
$$ \dim H^2_c(\sL_2(\R), {\cal F}_\lambda) = \cases{ 
0 & for $\lambda \not=0,1$ \cr 
1 & for $\lambda =0$ \cr 
2 & for $\lambda =1$, \cr} 
\quad 
\dim H^1_c(\sL_2(\R), {\cal F}_\lambda) = \cases{ 
0 & for $\lambda \not=0,1$ \cr 
2 & for $\lambda =0$ \cr 
1 & for $\lambda =1$. \cr} $$
For $\lambda = 0$ the flux homomorphism yields an injective map 
$$ H^2_c(\sL_2(\R), {\cal F}_\lambda) 
\to \Hom(\pi_1(\SL_2(\R)), H^1_c(\sL_2(\R), {\cal F}_\lambda)
\cong  H^1_c(\sL_2(\R), {\cal F}_\lambda), \leqno(10.2) $$ 
so that we only obtain non-trivial abelian extensions of the universal covering group 
$\tilde\SL_2(\R)$. For $\lambda = 1$ the kernel of (10.2) is one-dimensional and 
spanned by $[\omega]$, so that $[\omega]$ is, up to scalar multiples, 
the only non-trivial $2$-cohomology class associated to the modules 
${\cal F}_\lambda$ which integrates to a group cocycle on $\SL_2(\R)$. 

\sectionheadline{XI. Central extensions of groups of volume preserving 
diffeomorphisms} 

In the present section we discuss certain central extensions of the group 
$\Diff(M,\mu)$ of diffeomorphisms of a compact connected 
orientable manifold $M$ preserving a volume form $\mu$, resp., 
its identity component $D(M,\mu)$. 
Each closed $\z$-valued $2$-form $\omega$ 
on $M$ defines a central extension of the 
corresponding Lie algebra ${\cal V}(M,\mu)$ of $\mu$-divergence free 
vector fields because composing integration over $M$ with respect to $\mu$ 
with the $C^\infty(M,\z)$-valued cocycle defined by the $2$-form 
(cf.\ Sectio~IX) leads to a 
$\z$-valued $2$-cocycle, the so-called {\it Lichnerowicz cocycle} 
(cf.\ [Vi02], [Li74]). 
We shall see that if $\pi_2(M)$ vanishes, then the only obstruction to the 
integrability of the corresponding central extension is given by 
the flux homomorphism $\pi_1(D(M,\mu)) \to H^1_{\rm dR}(M,\z)$. 
If $M = G$ is a compact Lie group, 
we show that the flux becomes trivial on the covering group 
$\tilde D(G,\mu)$ of 
$D(G,\mu)$ acting on the universal covering manifold $\tilde G$ of $G$, 
which leads to central Lie group extensions of this group.

\subheadline{Some facts on the flux homomorphism for volume forms} 

In this  short subsection we collect some facts on the flux homomorphism 
of a volume form on a compact connected manifold. These results will be used 
to show that each closed $2$-form on a compact Lie group $G$  
defines a central extension of the covering $\tilde D(G,\mu)$ 
of identity component $D(G,\mu)$ of the 
group of volume preserving diffeomorphisms of $G$ which acts faithfully on 
the universal covering group $\tilde G$. 

Let $M$ be a smooth compact 
manifold, $\z$ a sequentially complete locally convex space 
and $\omega \in \Omega^p(M,\z)$ a closed $\z$-valued $p$-form. 
For a piecewise smooth curve $\alpha \: I \to \Diff(M)$ we define the 
{\it flux form} 
$$ \tilde F_\omega(\alpha) 
:= \int_0^1 \alpha(t)^*\big(i_{\delta^l(\alpha)(t)}\omega\big)\, dt 
= \int_0^1 i_{\alpha(t)^{-1}.\alpha'(t)}(\alpha(t)^*\omega)\, dt \in \Omega^{p-1}(M,\z). $$

Let $\alpha \: I \to \Diff(M)$ be a piecewise smooth path and 
$\sigma \: \Delta_{p-1} \to M$ a smooth singular simplex. Further define 
$$ \alpha.\sigma \: I \times \Delta_{p-1} \to M, \quad (t,x) \mapsto \alpha(t).\sigma(x). $$
Then 
$$ \eqalign{ 
&\ \ \ \ ((\alpha.\sigma)^*\omega)(t,x)\Big({\partial \over \partial t}, v_1, \ldots, v_{p-1}\Big) \cr
 &= \omega(\alpha(t).\sigma(x))\big(\alpha'(t)(\sigma(x)), \alpha(t).d\sigma(x)v_1, 
\ldots, \alpha(t).d\sigma(x)v_{p-1}) \cr
 &= (\alpha(t)^*\omega)(\sigma(x))\big(\alpha(t)^{-1}.\alpha'(t)(\sigma(x)), d\sigma(x)v_1, 
\ldots, d\sigma(x)v_{p-1}) \cr
 &= \big(i_{\alpha(t)^{-1}.\alpha'(t)}\big(\alpha(t)^*\omega\big)\big)(\sigma(x))
\big(d\sigma(x)v_1, \ldots, d\sigma(x)v_{p-1}) \cr} $$
(cf.\ [NV03, Lemma 1.7]) implies 
$$ \int_{\alpha.\sigma} \omega 
= \int_{I \times \Delta_{p-1}} (\alpha.\sigma)^*\omega = \int_\sigma \tilde F_\omega(\alpha). $$
We thus obtain 
$$ \int_{\alpha.\Sigma} \omega = \int_\Sigma \tilde F_\omega(\alpha) $$
for each singular chain $\Sigma$ if we extend the map 
$\sigma \mapsto \alpha.\sigma$ additively to the group of piecewise smooth singular 
chains. If $\Sigma$ is a boundary and $\alpha$ is closed, then 
$\alpha.\Sigma$ is a boundary, so that the integral vanishes by Stoke's Theorem, 
and therefore $\int_\Sigma \tilde F_\omega(\alpha)$ vanishes. We conclude that 
$\tilde F_\omega(\alpha)$ is a closed $(p-1)$-form, so that we obtain a group homomorphism 
$$ F_\omega \: \pi_1(\Diff(M)) \to H^{p-1}_{\rm dR}(M,\z), \quad 
[\alpha] \mapsto [\tilde F_\omega(\alpha)]. $$

\Lemma XI.1. If $M$ is an oriented compact 
manifold of dimension $n$, $m_0 \in M$,  
and $\mu$ a volume form on $M$ with $\int_M \mu = 1$, then the corresponding 
flux homomorphism 
$$ F_\mu \: \pi_1(\Diff(M)) \to H^{n-1}_{\rm dR}(M,\R), \quad 
[\alpha] \mapsto [\tilde F_\mu(\alpha)] $$
factors through the kernel of the map 
$$ \pi_1(\ev_{m_0}) \: \pi_1(\Diff(M)) \to \pi_1(M,m_0). $$

\Proof. (We are grateful to Stephan Haller for communicating the idea of the 
following proof.)  To each smooth loop $\alpha \: \SS^1 \to \Diff(M)$ with $\alpha(1) = \id_M$ 
we associate 
a locally trivial fiber bundle $q_\alpha \: P_\alpha \to \SS^2$ 
whose underlying topological space is obtained as follows. 
We think of $\SS^2$ as a union of two closed discs $B_1$ and $B_2$ 
with $B_1 \cap B_2 = \SS^1$. Then we put 
$$ P_\alpha := \Big((B_1 \times M) \dot\cup (B_2 \times M)\Big)/\sim, $$
where 
$$ (x,m) \sim (x', m') \: \Leftrightarrow \cases{ 
x = x' \not\in \partial B_1 \cup \partial B_2, \ m = m' \cr 
x = x' \in \partial B_1, m' = \alpha(x)(m). \cr} $$
Then $q_\alpha([x,m]) := x$ defines the structure of a locally trivial fiber bundle 
with fiber $M$ over $\SS^2$. 

A section of $P_\alpha$ is a pair of two continuous maps 
$\tilde\sigma_j \: B_j \to M$, $j = 1,2$, such that the 
restrictions $\sigma_j := \tilde\sigma_j \res_{\partial B_j}$ satisfy 
$\sigma_2(x) = \alpha(x)(\sigma_1(x))$ for all 
$x \in \partial B_j$. This means that $\sigma_1$ and $\sigma_2$ are contractible 
loops in $M$ with $\alpha.\sigma_1 = \sigma_2$. Conversely, every pair 
of contractible loops $\sigma_1$ and $\sigma_2$ in $M$ satisfying 
$\alpha.\sigma_1 = \sigma_2$ can be extended to continuous maps 
$B_j \to M$ and thus to a section of $P_\alpha$. 

If $\sigma_1$ is a contractible loop based in $m_0$, then 
$\alpha.\sigma_1$ is a loop based in $m_0$ homotopic to the loop 
$x \mapsto \alpha(x)(m_0)$. Therefore the existence of a continuous section of 
$P_\alpha$ is equivalent to  $[\alpha] \in \ker \pi_1(\ev_{m_0})$.

Suppose that $[\alpha] \in \ker \pi_1(\ev_{m_0})$ 
and that $\sigma \: \SS^2 \to P_\alpha$ 
is a corresponding section. It follows easily from the construction of $P_\alpha$ 
that the manifold $P_\alpha$ is orientable if $M$ is orientable. Hence 
the $2$-cycle $[\sigma]$ has a Poincar\'e dual 
$[\beta] \in H^n_{\rm sing}(P_\alpha,\Z)$ whose restriction to a fiber $M$ 
is the Poincar\'e dual of the intersection of $\im(\sigma)$ with a fiber, 
hence the fundamental class $[\mu] \in H^n_{\rm sing}(M,\Z)$ ([Bre93, p.372]). 
Therefore the fundamental class of $M$ extends to an $n$-dimensional cohomology class 
in $P$. 

On the other hand we obtain from [Sp66, p.455] 
the exact Wang cohomology sequence associated 
to $P_\alpha$: 
$$  \ldots \to H^n_{\rm sing}(P_\alpha,\Z) \to H^n_{\rm sing}(M,\Z) 
\sssmapright{\partial_\alpha} 
H^{n-1}_{\rm sing}(M,\Z) \to H^{n+1}_{\rm sing}(P,\Z) \to \ldots, $$
where $\partial_\alpha$ satisfies 
$$\la \partial_\alpha[\beta], [\Sigma]\ra = \la [\beta], [\alpha.\Sigma]\ra $$
for each $(n-1)$-cycle $\Sigma$ in $M$, 
and the kernel of $\partial_\alpha$ 
consists of those cohomology classes extending to $P_\alpha$. 
As this is the case for the fundamental class of $M$, it follows that 
$[\alpha.\Sigma] = 0$ holds for all $(n-1)$-cycles $\Sigma$ on $M$. 
We conclude 
that $\tilde F_\mu(\alpha)$ is an exact 
$(n-1)$-form if $[\alpha] \in \ker \pi_1(\ev_{m_0})$. 
\qed

\Remark XI.2. Suppose that $G$ is a compact Lie group of dimension $d$. 
Then $G$ is orientable 
and we can identify $G$ with the group $\lambda(G)$ of left translations in $\Diff(G)$. 
Then 
$$\Diff(G) = \Diff(G)_\1 \lambda(G) \cong \Diff(G)_\1 \times G $$
as smooth manifolds. In particular we have 
$$ \pi_1(\Diff(G)) \cong \pi_1(\Diff(G)_\1) \times \pi_1(G). $$
If $\mu$ is a normalized biinvariant volume form on $G$, then Lemma IX.1 implies 
that the corresponding flux homomorphism 
$$ F_\mu \: \pi_1(\Diff(G)) \to H^{d-1}_{\rm dR}(G,\R)  $$
factors through a homomorphism 
$$ F_\mu^\sharp \: \pi_1(G) \to H^{d-1}_{\rm dR}(G,\R).  $$

Let $q_G \: \tilde G \to G$ denote the universal covering homomorphism 
and 
$$ \tilde\Diff(G) := \{ \tilde\phi \in \Diff(\tilde G) \: 
(\exists \phi \in \Diff(G))\ \phi \circ q_G = q_G \circ \tilde \phi\}. $$
Then we have a canonical homomorphism 
$$ Q_G \: \tilde\Diff(G) \to \Diff(G), \quad \tilde \phi \mapsto \phi $$
whose kernel coincides with the group of deck transformations that is 
isomorphic to $\pi_1(G)$. We endow $\tilde\Diff(G)$ with the Lie group structure 
turning $Q_G$ into a covering map. We then have 
$$ \tilde\Diff(G) 
= \tilde\Diff(G)_\1 \tilde G \cong \tilde\Diff(G)_\1 \rtimes \tilde G
\cong \Diff(G)_\1 \rtimes \tilde G $$
as smooth manifolds, so that 
$$ \pi_1(\tilde\Diff(G)) \cong \pi_1(\Diff(G)_\1). $$
The identity component $\tilde\Diff(G)_0$ is a covering of 
$\Diff(G)_0$ and since the flux homomorphism vanishes on its fundamental group 
(Lemma IX.1), the flux cocycle 
$$ f_\mu \: {\cal V}(G) \to \hat H^{d-1}_{\rm dR}(G,\R), \quad X \mapsto [i_X\mu] $$
integrates to a group cocycle 
$$ F_\mu \: \tilde\Diff(G)_0 \to \hat H^{d-1}_{\rm dR}(G,\R) 
= \Omega^{d-1}(G,\R)/ d\Omega^{d-2}(G,\R)$$
with $D F_\mu = f_\mu$. 

\subheadline{Application to central extensions} 

In this subsection we apply the tools developed in 
the present paper to central extensions of 
groups of volume preserving diffeomorphisms of compact manifolds. 

Let $M$ denote an orientable connected  
compact manifold and $\mu$ a volume form on $M$, normalized by $\int_M \mu = 1$. 
We write 
$$D(M,\mu) := \{ \phi \in \Diff(M)^{\rm op} \: \phi^*\mu = \mu\}_0 $$
for the identity component of the 
{\it group of volume preserving diffeomorphisms of $(M,\mu)$} and 
$$ \g_\mu := {\cal V}(M,\mu) := \{ X \in {\cal V}(M) \: {\cal L}_X \mu = 0\} $$
for its Lie algebra.  Further let $\tilde D(M,\mu) \subeq \Diff(\tilde M)$ 
denote the identity component of the inverse image of 
$D(M,\mu)$ in $\tilde\Diff(M)$. Then we have a covering map 
$\tilde D(M,\mu) \to D(M,\mu)$
which need not be universal. We write $\tilde{D(M,\mu)}$ for the universal 
covering group of $D(M,\mu)$ which also is a covering group of 
$\tilde D(M,\mu)$. 

Let $\z$ be a Fr\'echet space. 
On the space $C^\infty(M,\z)$ of smooth $\z$-valued functions on $M$ we then have the 
integration map 
$$ I \: C^\infty(M,\z) \to \z, \quad f \mapsto \int_M f \mu. $$
Then $I$ is equivariant for the natural action of $D(M,\mu)$ on 
$C^\infty(M,\z)$, where we consider $\z$ as a trivial module. 
On the infinitesimal level this means that 
$$ \int_M (X.f)\, \mu = 0 \quad \hbox{ for}  \quad f \in C^\infty(M,\R), 
\quad X \in {\cal V}(M,\mu). $$

Each closed $\z$-valued $p$-form $\omega \in \Omega^p(M,\z)$ defines a 
$C^\infty(M,\z)$-valued $p$-cochain for the action of the Lie algebra 
$\g_\mu$ on $C^\infty(M,\z)$ and since $I$ is $\g_\mu$-equivariant, 
we obtain continuous linear maps 
$$ \Phi \: \Omega^p(M,\z) \to C^p_c(\g_\mu, \z), 
\quad \Phi(\omega)(X_1,\ldots, X_p) := I(\omega(X_1, \ldots, X_p)) 
= \int_M \omega(X_1, \ldots, X_p)\, \mu. $$
The equivariance of $I$ implies that $\Phi(d\omega) = d_{\g_\mu} \Phi(\omega)$, 
so that $\Phi$ induces maps 
$$ \Phi \: H^p_{\rm dR}(M,\z) \to H^p_c(\g_\mu, \z). $$

\Remark XI.3. If $\pi_2(M) = \{0\}$ and $\tilde{D(M,\mu)}$ denotes the 
simply connected covering group of $D(M,\mu)$, then for each closed $2$-form 
$\omega \in Z^2_{\rm dR}(M,\z)$ the period map of the corresponding 
Lie algebra cocycle vanishes (Proposition~IX.5), 
so that, in view of Theorem VII.2, $\Phi$ induces a map 
$$ \Phi \: H^2_{\rm dR}(M,\z) \to H^2_s(\tilde{D(M,\mu)}, \z). $$

If, more generally, $\Gamma_Z \subeq \z$ is a discrete subgroup with 
$\int_{\pi_2(M)} \omega \subeq \Gamma_Z$ and $Z := \z/\Gamma_Z$, then 
Theorem VII.2 implies that 
the Lie algebra cocycle $\omega$ integrates to a central extension 
$$ Z \into \hat D(M,\mu) \onto \tilde{D(M,\mu).}$$

Let 
$${\cal V}(M,\mu)_{\rm ex} := \{ X \in {\cal V}(M,\mu) \: i_X \mu \in d \Omega^{p-2}(M,\R)\} $$
denote the Lie algebra of exact divergence free vector fields. 
It can be shown that this is the commutator algebra of ${\cal V}(M,\mu)$ 
(cf.\ [Li74]), hence a perfect Lie algebra. It follows in particular that 
$$H^1_c({\cal V}(M,\mu)_{\rm ex},\z) = \Hom_{\rm Lie\ alg}({\cal V}(M,\mu)_{\rm ex},\z) = \{0\} $$
vanishes for each trivial module $\z$. Therefore restricing the cocycles 
from above to ${\cal V}(M,\mu)_{\rm ex}$, resp.\ the corresponding connected subgroup 
$D(M,\mu)_{\rm ex}$ of exact volume preserving diffeomorphisms leads to a trivial 
flux homomorphism. Hence $\int_{\pi_2(M)} \omega \subeq \Gamma_Z$ 
implies the existence of a central $Z$-extension of $D(M,\mu)_{\rm ex}$. 
We refer to Ismagilov ([Is96]) and Haller-Vizman ([HV04]) for geometric constructions 
of these central extensions (for the case $\z = \R, Z = \T = \R/\Z$). 
\qed

\Proposition XI.4. Let $G$ be a compact connected Lie group and $\mu$ an invariant 
normalized volume form on $G$. 
Then the flux cocycle restricts to a surjective Lie algebra homomorphism 
$$ f_\mu \: {\cal V}(G,\mu) \to H^{d-1}_{\rm dR}(G,\R)  $$
whose kernel is the commutator algebra and whose restriction to $\z(\g) \subeq 
\g \subeq {\cal V}(G,\mu)$ is bijective. This Lie algebra homomorphism 
integrates to a homomorphism of connected Lie groups 
$$ F_\mu^G \: \tilde D(G,\mu) \to H^{d-1}_{\rm dR}(G,\R) $$
whose restriction to $Z(\tilde G)_0 \subeq \tilde G \subeq \tilde D(G,\mu)$ 
is an isomorphism. Moreover, each Lie algebra homomorphism 
$\phi_\g \: {\cal V}(G,\mu) \to \a$ to an abelian Lie algebra 
integrates to a group homomorphism 
$\phi_G \: \tilde D(G,\mu) \to \a$  which factors through $F_\mu^G$. 

\Proof. Since $f_\mu$ defines a Lie algebra homomorphism 
${\cal V}(G,\mu) \to H^{d-1}_{\rm dR}(G,\R)$, 
the restriction of the flux cocycle $F_\mu \: 
\tilde\Diff(G)_0 \to \hat H^{d-1}_{\rm dR}(G,\R)$ 
to the subgroup $\tilde D(G,\mu)$ is a group 
homomorphism 
$$ F_\mu \: \tilde D(G,\mu) \to H^{d-1}_{\rm dR}(G,\R) \cong H^{d-1}(\g,\R) $$
which on the subgroup $\tilde G$ of $\tilde D(G,\mu)$ is the Lie group 
homomorphism obtained by integrating the Lie algebra quotient homomorphism 
$$ \g \to H^{d-1}(\g,\R), \quad x \mapsto [i_x \mu_\g], $$
where $\mu_\g := \mu(\1) \in C^d(\g,\R)$. 
Note that Poincar\'e Duality implies that 
$$ H^{d-1}_{\rm dR}(G,\R) \cong H^{1}_{\rm dR}(G,\R) \cong 
\Hom(\g,\R) \cong \z(\g)^* $$
so that 
$H^{d-1}_{\rm dR}(G,\R) \cong Z(\tilde G)_0 \cong \z(\g)$
and we can think of the flux homomorphism as a group homomorphism 
$$ F_\mu^G \: \tilde D(G,\mu) \to \z(\g). $$

On the Lie algebra level we have $\g \subeq {\cal V}(G,\mu)$, 
$[{\cal V}(G,\mu), 
{\cal V}(G,\mu)] \subeq \ker f_\mu,$ and 
$f_\mu$ maps $\z(\g)$ isomorphically onto $H^{d-1}_{\rm dR}(G,\R)$. 
This leads to 
$$ {\cal V}(G,\mu) = [{\cal V}(G,\mu), {\cal V}(G,\mu)] \rtimes \z(\g) $$
with $H_1({\cal V}(G,\mu)) \cong \z(\g)$ and we conclude that  
the flux homomorphism $F_\mu^G \: \tilde D(G,\mu) \to \z(\g)$ 
is universal in the sense that each Lie algebra homomorphism 
${\cal V}(G,\mu) \to \a$, where $\a$ is an abelian Lie algebra, 
integrates to a Lie group homomorphism 
$\tilde D(G,\mu) \to \a$. 
\qed

\Theorem XI.5. Let $G$ be a connected compact Lie group, $\mu$ an invariant 
normalized volume form, 
$\z$ a sequentially complete locally convex space and $\omega \in \Omega^2(G,\z)$ a 
closed $2$-form. Then the Lichnerowicz cocycle on ${\cal V}(G,\mu)$ given by 
$$ (X,Y) \mapsto \int_G \omega(X,Y)\cdot \mu $$ 
integrates to a central Lie group extension 
$$ \z \to \hat D(G,\mu)\to \tilde D(G,\mu). $$

\Proof. First we recall that $\pi_2(G) = \{0\}$ ([Ca52]), so that Remark XI.3 implies 
that the period map of $\tilde D(G,\mu)$  
vanishes for each closed $2$-form $\omega \in \Omega^2(G,\z)$ on $G$. 
Moreover, the flux cocycle is a Lie 
algebra homomorphism 
$$ f_\omega \: \g_\mu = {\cal V}(G,\mu) \to H^1_c(\g_\mu,\z) 
\cong \Hom(\g_\mu,\z) \cong \Hom(\z(\g),\z) $$
so that Proposition XI.4 implies that the corresponding flux homomorphism vanishes on 
the fundamental group 
$\pi_1(\tilde D(G,\mu))$, so that Theorem VII.2 implies that 
$\omega$ defines a Lie algebra cocycle in $Z^2_c({\cal V}(G,\mu),\z)$ corresponding to a 
global central extension as required.  
\qed

\Remark XI.6. In view of 
$$ H^2_{\rm dR}(G,\z) \cong H^2_c(\g,\z) = H^2_c(\z(\g),\z) 
= \Alt^2(\z(\g),\z) = \Lin(\Lambda^2(\z(\g)), \z), $$ 
we obtain a universal Lichnerowicz cocycle with values in the space 
$\z := \Lambda^2(\z(\g))$. 
\qed

The preceding remark applies in particular to the $d$-dimensional torus $G = \T^d 
:= \R^d/\Z^d$, 
which we consider as the quotient of $\R^d$ modulo the integral lattice. 
We write $x_1, \ldots, x_d$ for the canonical coordinate functions on $\R^d$ 
and observe that their differential $d x_j$ can also be viewed as $1$-forms on 
$\T^d$. In this sense we have 
$$ H^2_{\rm dR}(\T^d, \R) \cong \bigoplus_{i < j} \R [d x_i \wedge d x_j]  
\cong \R^{{d \choose 2}}. $$
Therefore the central extensions of $\tilde D(\T^d,\mu)$ described above 
correspond to the central extensions of the corresponding Cartan type 
algebras discussed in [Dz92]. We conclude in particular that these cocycle 
do not integrate to central extensions of 
$D(\T^d, \mu)$, but that they integrate to central extensions of the covering 
group $\tilde D(\T^d, \mu)$ which we can considered as a group of 
diffeomorphisms of $\R^d$.

\sectionheadline{Appendix A. Differential forms and Alexander--Spanier
cohomology} 

\nin In this appendix we discuss a smooth version of Alexander--Spanier
cohomology for smooth manifolds and define a homomorphism of chain
complexes from the smooth Alexander--Spanier complex $(C^n_{AS,s}(M,A),
d_{AS}), n \geq 1,$ with values in
an abelian Lie group $A$ with Lie algebra $\a$ to the $\a$-valued 
de Rham complex $(\Omega^n(M,\a),d)$. In Appendix B this map is used to relate Lie
group cohomology to Lie algebra cohomology.  
The main point is Proposition~A.6 which provides an explicit map from 
smooth Alexander--Spanier cohomology to de Rham cohomology. 

\Definition A.1. {\rm(1)} Let $M$ be a smooth manifold and $A$ an abelian Lie
group. For $n \in \N_0$ let $C^n_{AS,s}(M,A)$ denote the set of 
germs of smooth $A$-valued functions on the diagonal in $M^{n+1}$. For
$n = 0$ this is the space $C^0_{AS,s}(M,A) \cong C^\infty(M,A)$ of
smooth $A$-valued functions on $M$.  
An element $[F]$ of this space is represented by a smooth function 
$F \: U \to A$, where $U$ is an open neighborhood of the diagonal in
$M^{n+1}$, and two functions $F_i \: U_i \to A$, $i =1,2$, define the
same germ if and only of their difference vanishes on a neighborhood
of the diagonal. The elements of the space $C^n_{AS,s}(M,A)$ are
called {\it smooth $A$-valued Alexander--Spanier $n$-cochains on $M$}. 

We have a differential 
$$d_{AS} \: C^n_{AS,s}(M,A) \to C^{n+1}_{AS,s}(M,A) $$
given by 
$$ (d_{AS}F)(m_0, \ldots, m_{n+1}) := 
 \sum_{j=0}^{n+1} (-1)^j F(m_0, \ldots, \hat{m_j}, \ldots, m_{n+1}),
 $$ 
where $\hat{m_j}$ indicates omission of the argument $m_j$. 
To see that $d_{AS}F$ 
defines a smooth function on an open neighborhood of the diagonal in
$M^{n+2}$, consider for $i =0,\ldots, n+1$ 
the projections $p_i \: M^{n+2} \to M^{n+1}$ 
obtained by omitting the $i$-th component. Then for each open subset
 $U \subeq M^{n+1}$  containing the diagonal the subset 
$\bigcap_{i =0}^{n+1} p_i^{-1}(U)$ 
is an open neighborhood of the diagonal in $M^{n+2}$ on which $d_{AS} F$ is defined. 
It is easy to see that $d_{AS}$ is well-defined on germs and that we
 thus obtain a differential complex $(C^*_{AS,s}(M,A), d_{AS})$. 
Its cohomology groups are denoted $H^n_{AS,s}(M,A)$. 

\par\nin (2) If $M$ is a smooth manifold, then an atlas for the tangent bundle $TM$ is obtained
directly from an atlas of $M$, but we do not consider the cotangent
bundle as a manifold because this requires to choose a topology on the
dual spaces, for which there are many possibilities. 
Nevertheless, there is a natural concept of a smooth $p$-form on $M$. 
If $V$ is a locally convex space, then a
{\it $V$-valued $p$-form} $\omega$ on $M$ is a function 
$\omega$ which associates to each $x \in M$ a $k$-linear 
alternating map $T_x(M)^p \to V$ such that in local coordinates the
map 
$(x,v_1, \ldots, v_p) \mapsto \omega(x)(v_1, \ldots, v_p)$
is smooth. We write $\Omega^p(M,V)$ for the space of smooth $p$-forms
on $M$ with values in $V$. 

The {\it de Rham differential} 
$d \: \Omega^p(M,V) \to \Omega^{p+1}(M,V)$ is defined by 
$$ \eqalign{ (d\omega)(x)(v_0, \ldots, v_p) 
&:= \sum_{i = 0}^p (-1)^i \big(X_i.\omega(X_0, \ldots, \hat X_i, \ldots,
X_p)\big)(x) \cr
&\ \ \ \ + \sum_{i < j} (-1)^{i+j} \omega([X_i, X_j], X_0, \ldots, \hat X_i,
\ldots, \hat X_j, \ldots, X_p)(x) \cr}$$
for $v_0,\ldots, v_p \in T_x(M)$, where $X_0,\ldots, X_p$ are smooth
vector fields on a neighborhood of $x$ with $X_i(x) = v_i$. 

To see that $d$ defines indeed a map $\Omega^p(M,V) \to
\Omega^{p+1}(M,V)$ one has to verify that the right hand side of the
above expression does not depend on the choice of the vector fields
$X_i$ with $X_i(x) = v_i$ and that it defines an element of
$\Omega^{p+1}(M,V)$, i.e., in local coordinates the map  
$$ (x,v_0, \ldots, v_p) \mapsto (d\omega)(x)(v_0, \ldots, v_p) $$
is smooth, multilinear and alternating in $v_0, \ldots, v_p$. 
For the proof we refer to (cf.\ [KM97]). 

Extending $d$ to a linear map on $\Omega(M,V) := \bigoplus_{p \in
\N_0} \Omega^p(M,V)$, we have the relation $d^2 = 0$. The space 
$$ Z^p_{\rm dR}(M,V) := \ker(d\res_{\Omega^p(M,V)}) $$
of {\it closed forms} therefore contains the space 
$B^p_{\rm dR}(M,V) := d(\Omega^{p-1}(M,V))$
of exact forms, and 
$$ H^p_{\rm dR}(M,V) := Z^p_{\rm dR}(M,V) / B^p_{\rm dR}(M,V) $$
is the {\it $V$-valued de Rham cohomology space of $M$}. 
\qed

\Definition A.2. If $M$ is a smooth manifold, $A$ an abelian Lie
group, $\a$ its Lie algebra, 
$f \: M \to A$ a smooth function and $Tf \: TM\to TA$ its tangent map, 
then we define the {\it logarithmic derivative of $f$} as the $\a$-valued $1$-form 
$$ df \: TM \to \a, \quad v \mapsto f(m)^{-1}.Tf(v),
\quad \hbox{ for } \quad v \in T_m(M). $$

In terms of the canonical trivialization $\theta \: TA \to A \times
\a, v \mapsto a^{-1}.v$ (for $v \in T_a(A)$) of the tangent bundle of $A$, this means that 
$$ d f = \pr_2 \circ \theta \circ Tf \: TM \to \a. 
\qeddis 

\def\vbar{\hbox{\vrule width0.5pt 
                height 5mm depth 3mm${{}\atop{{}\atop{\scriptstyle t_i=0}}}$}}

\Definition A.3. Let $M_1, \ldots, M_n$ be smooth manifolds, $A$ an
abelian Lie group, and 
$$ f \: M_1 \times \ldots \times M_n \to A $$
be a smooth function. 
For $n \in \N$ we define a function 
$$ d^n f \: TM_1 \times \ldots \times TM_n \to \a $$
as follows. Let $q \: TM \to M$ be the canonical projection. 
For $v_1, \ldots, v_n \in TM$ with $q(v_i) = m_i$ we consider smooth curves 
$\gamma_i \: ]-1,1[ \to M$ with 
$\gamma_i(0) = m_i$ and 
$\gamma_i'(0) = v_i$ and define 
$$ (d^n f)(m_1,\ldots, m_n)(v_1,\ldots, v_n) := {\partial^n \over \partial t_1 \cdots
\partial t_n} \vbar f(\gamma_1(t_1), \ldots,
\gamma_n(t_n)),  $$
where for $n \geq 2$ the iterated higher derivatives are derivatives of $\a$-valued
functions in the sense of Definition A.2. 
One readily verifies that the right hand side does not depend on the
choice of the curves $\gamma_i$ and that it defines for each tuple 
$(m_1, \ldots, m_n) \in M_1\times \ldots \times M_n$ a continuous $n$-linear map 
$$ (d^nf)(m_1,\ldots, m_n) \: 
T_{m_1}(M_1) \times \ldots \times T_{m_n}(M_n) \to \a. $$

If $X$ is a smooth vector field on $M_i$, then we also define a smooth
function 
$$ \partial_i(X)f \: M_1 \times \ldots \times M_n \to \a, \quad 
(m_1, \ldots, m_n) \mapsto 
df(m_1,\ldots, m_n)(0,\ldots,0, X(m_i), 0,\ldots, 0) $$
by the partial derivative of $f$ in the direction of the vector field
$X$. 
For vector fields $X_i$ on $M_i$ we then obtain by iteration of this
process 
$$ \big(\partial_1(X_1)\cdots \partial_n(X_n) f\big)(m_1,\ldots, m_n) 
= (d^nf)(m_1,\ldots, m_n)(X_1(m_1),\ldots, X_n(m_n)) $$
and 
$$ \partial_1(X_1)\cdots \partial_n(X_n) f \: M_1 \times \ldots \times
M_n \to \a $$ 
is a smooth function. 
\qed

\Definition A.4. Let $M$ be a smooth manifold and $A$ an abelian Lie
group. We write $\Delta_n \: M \to M^{n+1}, m \mapsto (m,\ldots, m)$ for
the diagonal map. 

For $[F] \in C^n_{AS,s}(M,A)$, $p \in M$ and $v_1,\ldots, v_n \in
T_p(M)$ we define 
$$ \tau(F)(p)(v_1,\ldots, v_n) := 
\sum_{\sigma \in S_n} \sgn(\sigma)\cdot (d^n F)(p,\ldots,
p)(0,v_{\sigma(1)},\ldots, 
v_{\sigma(n)}) $$
and observe that 
$\tau(F)$ defines a smooth $\a$-valued $n$-form on $M$ depending only
on the germ $[F]$ of $F$. We thus obtain for $n \geq 1$ a group homomorphism 
$$ \tau \: C^n_{AS,s}(M,A) \to \Omega^n(M,\a). $$
If $A = \a$, then we also define $\tau$ for $n = 0$ as the identical
map 
$$ \tau \: C^0_{AS,s}(M,A) \cong C^\infty(M,A) \to \Omega^0(M,\a)
\cong C^\infty(M,\a). $$

If $X_1, \ldots, X_n$ are smooth vector fields on an open subset $V
\subeq M$, we have on $V$ the relation 
$$ \tau(F)(X_1,\ldots, X_n) 
= \sum_{\sigma \in S_n} \sgn(\sigma)
\cdot \big(\partial_1(X_{\sigma(1)})\cdots
\partial_n(X_{\sigma(n)}).F\big) \circ \Delta_n. $$
As the operators $\partial_i(X)$ and $\partial_j(Y)$ commute for
$i \not=j$ and vector fields $X$ and $Y$ on $M$, this
can also be written as 
$$ \tau(F)(X_1,\ldots, X_n) 
= \sum_{\sigma \in S_n} \sgn(\sigma)
\cdot \big(\partial_{\sigma(1)}(X_1)\cdots
\partial_{\sigma(n)}(X_n).F\big) \circ \Delta_n. $$

\nin For small $n$ we have in particular the formulas
\par\nin $n =0$: $\tau(F) = F$ (if $A = \a$). 
\par\nin $n =1$: $\tau(F)(X) = \partial_1(X).F$.
\par\nin $n =2$: $\tau(F)(X,Y) = \partial_1(X)\partial_2(Y).F- 
\partial_1(X)\partial_2(Y).F$.
\qed

The following proposition builds on a construction one finds in 
the appendix of [EK64]. First we need a combinatorial lemma. 

\Lemma A.5. Let $\sigma \in S_{n+1}$ be a permutation with 
$k := \sigma(1) < \ell := \sigma(i+1)$
and such that the restriction of $\sigma$ defines an increasing map 
$\{1,\ldots, n\} \setminus \{1,i+1\} \to 
\{1,\ldots, n\} \setminus \{k,\ell\}.$
Then $\sgn(\sigma) = (-1)^{i+k+l}.$

\Proof. Replacing $\sigma$ by $\sigma_1 := \sigma \circ \alpha$, where 
$\alpha = (i+1\ \ i\ \  i-1\ \ldots \ 3\ 2)$
is a cycle of length $i$, we obtain a permutation $\sigma_1$ 
that restricts to an increasing map 
$$ \{3,4,\ldots, n\}  \to \{1,\ldots, n\} \setminus \{k,\ell\}.$$
Next we put $\sigma_2 := \beta \circ \sigma_1$, where 
$\beta = (1\ 2\ 3\ \ldots \ k-1\ k)$
is a cycle of length $k$ to obtain an increasing map 
$$ \{3,4,\ldots, n\}  \to \{2,\ldots, n\} \setminus \{\ell\}.$$
Eventually we put $\sigma_3 := \gamma \circ \sigma_2$, where 
$\gamma = (2\ 3\ \ldots \ \ell-1\ \ell)$
is a cycle of length $\ell-1$ to obtain an increasing map 
$$ \{3,4,\ldots, n\}  \to \{3,4,\ldots, n\},$$
which implies that $\sigma_3$ fixes all these elements. 
Further 
$$\sigma_3(1) 
= \gamma\beta\sigma\alpha(1)
= \gamma\beta\sigma(1) 
= \gamma\beta(k)
= \gamma(1) = 1 $$
implies that $\sigma_3 = \id$. This implies that 
$$ \sgn(\sigma) = \sgn(\alpha)\sgn(\beta)\sgn(\gamma) 
= (-1)^{i-1} (-1)^{k-1} (-1)^{\ell} = (-1)^{i+k+\ell}. 
\qeddis

The following proposition generalizes an observation of van Est and
Korthagen in the Appendix of [EK64]: 

\Proposition A.6. {\rm(van Est--Korthagen)} If $M$ is 
smooth manifold, then the map 
$$ \tau \: C^n_{AS,s}(M,A) \to 
\cases{ C^\infty(M,A) & for $n = 0$ \cr 
\Omega^n(M,\a) & for $n \geq 1$\cr} $$
intertwines the Alexander--Spanier differential with the de Rham
differential, hence induces a map 
$$ \tau \: H^n_{AS,s}(M,A) \to H^n_{\rm dR}(M,\a). $$

\Proof. We have to show that 
$\tau(d_{AS} F) = d \tau(F)$ holds for $F \in C^\infty(U,A)$, where
$U$ is an open neighborhood of the diagonal in $M^{n+1}$.  

From the chain rule we obtain for a vector field $Y$ on $M$ the
relation 
$$ \leqalignno{ 
Y.\big(\big(\partial_1(X_1)\cdots \partial_n(X_n).F\big) &\circ \Delta_n\big)
= \big(\partial_0(Y) \partial_1(X_1)\cdots \partial_n(X_n).F\big)
\circ \Delta_n \cr
&\ \ \ \ + \sum_{i=1}^n 
\big(\partial_1(X_1)\cdots \partial_i(Y)\partial_i(X_i)\cdots 
\partial_n(X_n).F\big) \circ \Delta_n. &(A.1) \cr} $$
Now let 
$$ F_i(x_0, \ldots, x_{n+1}) 
:=F(x_0, \ldots, \hat{x_i}, \ldots, x_{n+1}). $$
Then 
$$ F_i \circ  \Delta_{n+1} = F \circ \Delta_n \leqno(A.2) $$
and $d_{AS} F = \sum_{i=0}^{n+1} (-1)^i F_i$. 
Since the function $F_i$ is independent of $x_i$, we have 
$$ \partial_1(X_1)\cdots \partial_{n+1}(X_{n+1}).F_i = 0, \quad i \geq
1. \leqno(A.3) $$
Therefore 
$$ \partial_1(X_1)\cdots \partial_{n+1}(X_{n+1}).(d_{AS} F) 
= \partial_1(X_1)\cdots \partial_{n+1}(X_{n+1})(F_0) 
= \big(\partial_0(X_1)\cdots \partial_n(X_{n+1}) F\big)_0. $$
In view of (A.2) and (A.1), this leads to 
$$ \eqalign{ 
&\ \ \ \ \big(\partial_1(X_1)\cdots \partial_{n+1}(X_{n+1}).(d_{AS}
F)\big)\circ \Delta_{n+1} = \big(\partial_0(X_1)\cdots \partial_n(X_{n+1}).F\big) \circ
\Delta_n \cr 
&= X_1.\Big(\big(\big(\partial_1(X_2)\cdots
\partial_n(X_{n+1}).F\big)\Big) \circ \Delta_n \cr
&\ \ \ \ - 
\sum_{i=1}^n \big(\partial_1(X_2)\cdots \partial_i(X_1)\partial_i(X_{i+1})\cdots 
\partial_n(X_{n+1}).F\big) \circ \Delta_n.\cr} $$
Alternating the first summand, we get an expression of the form 
$$ \eqalign{ 
&\ \ \ \ \sum_{\sigma \in S_{n+1}}\sgn(\sigma) 
X_{\sigma(1)}.\Big(\partial_1(X_{\sigma(2)})\cdots
\partial_n(X_{\sigma(n+1)}).F\Big) \circ \Delta_n \cr
&= \sum_{i=1}^{n+1} \sum_{\sigma(1) = i} \sgn(\sigma) 
X_i.\Big(\partial_1(X_{\sigma(2)})\cdots
\partial_n(X_{\sigma(n+1)}).F\Big) \circ \Delta_n \cr} $$
We write any permutation $\sigma \in S_{n+1}$ with $\sigma(1) = i$ as 
$\sigma = \alpha_i \beta$, where $\beta(1) = 1$ and $\alpha_i(1) = i$ 
and $\alpha_i$ is the cycle 
$$ \alpha_i = (i\ i-1\ i-2\ \ldots\ 2\ 1). $$
We further identify $S_n$ with the stabilizer of $1$ in $S_{n+1}$. 
Then the above sum turns into 
$$ \eqalign{ 
&= \sum_{i=1}^{n+1} \sgn(\alpha_i) \sum_{\beta \in S_n} \sgn(\beta) 
X_i.\Big(\partial_1(X_{\alpha_i\beta(2)})\cdots
\partial_n(X_{\alpha_i\beta(n+1)}).F\Big) \circ\Delta_n \cr
&= \sum_{i=1}^{n+1} (-1)^{i-1} 
X_i.\tau(F)(X_{\alpha_i(2)}, \ldots, X_{\alpha_i(n+1)}) \cr
&= \sum_{i=1}^{n+1} (-1)^{i-1} 
X_i.\tau(F)(X_1, \ldots, \hat{X_i}, \ldots, X_{n+1}).\cr} $$

In view of 
$$ \eqalign{ d(\tau(F))(X_1,\ldots, X_{n+1}) 
&= \sum_{i=1}^{n+1} (-1)^{i-1} X_i.\tau(F)(X_1,\ldots, \hat X_i,
\ldots, X_{n+1}) \cr
&\ \ \ \  + 
\sum_{k<\ell} (-1)^{k+\ell} \tau(F)([X_k, X_\ell], X_1,\ldots, \hat X_k,
\ldots, \hat X_\ell, \ldots, X_{n+1}), \cr} $$
and 
$$ \leqalignno{ 
&\ \ \ \ -\sum_{k< \ell} (-1)^{k+\ell} \tau(F)([X_k, X_\ell], X_1,\ldots, \hat X_k,
\ldots, \ldots, \hat X_\ell, \ldots X_{n+1}) \cr 
&= \sum_{k< \ell} (-1)^{k+\ell+1} \sum_{\beta \in S_n} \sgn(\beta)
\big(\partial_{\beta(1)}([X_k, X_\ell])\partial_{\beta(2)}(X_1) \cdots \hat \partial(X_k) 
\cdots \hat \partial(X_\ell) \ldots \partial_{\beta(n)}(X_{n+1}).F\big) \circ
\Delta_n,\cr} $$
it remains to show that, as operators on functions on $M^{n+1}$, 
alternation of 
$$ \sum_{i=1}^n \partial_1(X_2)\cdots \partial_i(X_1)\partial_i(X_{i+1})\cdots 
\partial_n(X_{n+1}) \leqno(A.3) $$
leads to 
$$ \leqalignno{ 
&\ \ \ \ \sum_{k< \ell} (-1)^{k+\ell+1} \sum_{\beta \in S_n} \sgn(\beta)  
\partial_{\beta(1)}([X_k, X_\ell])\partial_{\beta(2)}(X_1) \cdots \hat \partial(X_k) 
\cdots \hat \partial(X_\ell) \ldots
\partial_{\beta(n)}(X_{n+1})  \cr 
&= \sum_{k < \ell} (-1)^{k + \ell+1} \la \partial_1 \wedge \ldots \wedge
\partial_n, [X_k, X_\ell] \wedge X_1 \wedge \cdots \wedge \hat X_k \wedge
\cdots \wedge \hat X_\ell \wedge \cdots \wedge X_{n+1}\ra\cr
&= \sum_{k < \ell} (-1)^{k + \ell+1} 
\sum_{i=1}^n (-1)^{i+1} \partial_i([X_k,X_\ell]) \circ 
\la \partial_1 \wedge \ldots \wedge \hat \partial_i\wedge \ldots \wedge \partial_n, 
X_1 \wedge \cdots \wedge \hat X_k \wedge 
\cdots \wedge \hat X_\ell \wedge \cdots \wedge X_{n+1}\ra.\cr} $$ 

Alternating (A.3) leads to the expression 
$$ \eqalign{ 
&\ \ \ \ \ \sum_{\sigma \in S_{n+1}} \sgn(\sigma) 
\sum_{i=1}^n \partial_1(X_{\sigma(2)})\cdots \partial_i(X_{\sigma(1)})
\partial_i(X_{\sigma(i+1)})\cdots \partial_n(X_{\sigma(n+1)})  \cr
&= \sum_{i = 1}^n \sum_{\sigma(1) < \sigma(i+1)} \sgn(\sigma) 
\sum_{i=1}^n \partial_1(X_{\sigma(2)})\cdots
\partial_i([X_{\sigma(1)}, X_{\sigma(i+1)}])
\cdots \partial_n(X_{\sigma(n+1)})  \cr
&= \sum_{i = 1}^n \sum_{k < \ell} 
\sum_{{\sigma(1)=k \atop\sigma(i+1)= \ell}} 
\sgn(\sigma) 
\sum_{i=1}^n \partial_1(X_{\sigma(2)})\cdots
\partial_i([X_k, X_\ell]) \cdots \partial_n(X_{\sigma(n+1)}).  \cr}$$
We can write each permutation $\sigma\in S_{n+1}$ as 
$\sigma = \sigma_0 \beta$, where $\beta$ fixes $1$ and $i+1$, so that
we can identify it with an element of $S_{n-1}$, and 
$$ \sigma_0 \: \{2,\ldots, n+1\} \setminus \{i+1\} \to \{1,\ldots,
n+1\} \setminus \{k,\ell\} $$ 
is increasing. 
In view of Lemma A.5, we then have 
$\sgn(\sigma_0) = (-1)^{i+k+\ell}$ for $k = \sigma(1)$ and $\ell = \sigma(i+1)$.
Therefore alternating (A.3) gives 
$$ \eqalign{ 
&= \sum_{i = 1}^n \sum_{k < \ell} (-1)^{i+k+\ell} 
\sum_{\beta \in S_{n-1}} \sgn(\beta) 
\sum_{i=1}^n \partial_1(X_{\sigma_0\beta(2)})\cdots
\partial_i([X_k, X_\ell]) \cdots \partial_n(X_{\sigma_0\beta(n+1)}) \cr 
&= \sum_{i = 1}^n \sum_{k < \ell} (-1)^{i+k+\ell} 
\partial_i([X_k, X_{\ell}]) \circ 
\la \partial_1 \wedge \cdots \wedge \hat\partial_i \wedge \cdots \wedge
\partial_n, X_{\sigma_0(2)} \wedge \cdots \wedge X_{\sigma_0(n+1)})\ra\cr 
&= \sum_{k < \ell} \sum_{i = 1}^n (-1)^{i+k+\ell} 
\partial_i([X_k, X_{\ell}]) \la \partial_1 \wedge \cdots \hat\partial_i \cdots
\wedge\partial_n, X_2 \wedge \cdots \wedge \hat X_k \wedge \cdots
\wedge \hat X_\ell 
\wedge \cdots \wedge X_{n+1})\ra.\cr 
}$$
This completes the proof of Proposition~A.6. 
\qed

\sectionheadline{Appendix B. Cohomology of Lie groups and Lie algebras} 

In this appendix we show that for $n \geq 2$ there is a natural ``derivation map''
$$ D_n \: H^n_s(G,A) \to H^n_c(\g,\a) $$
from locally smooth Lie group cohomology to continuous Lie algebra
cohomology. For $n = 1$ we have a map 
$D_1 \: Z^1_s(G,A) \to Z^1_c(\g,\a)$, and if, in addition, 
$A \cong \a/\Gamma_A$ holds for a discrete subgroup $\Gamma_A$ of
$\a$, then this map induces a map between the cohomology groups. 

\Definition B.1. Let $V$ be a {\it topological module of the
topological Lie algebra $\g$}. 
For $p \in \N_0$, let $C^p_c(\g,V)$ denote the space of continuous 
alternating maps $\g^p \to V$, i.e., 
the {\it Lie algebra $p$-cochains with values in the module $V$}. 
Note that $C^1_c(\g,V) = \Lin(\g,V)$ is the space of continuous linear
maps $\g \to V$. We use the convention $C^0_c(\g,V) = V$. 
We then obtain a chain complex with the differential 
$$ d_\g \: C^p_c(\g,V) \to C^{p+1}_c(\g,V) $$
given on $f \in C^p_c(\g,V)$ by 
$$ \eqalign{ (d_\g f)(x_0, \ldots, x_p) 
&:= \sum_{j = 0}^p (-1)^j x_j.f(x_0, \ldots, \hat x_j, \ldots, 
x_p) \cr
& + \sum_{i < j} (-1)^{i + j} f([x_i, x_j], x_0, \ldots, \hat
x_i, \ldots, \hat x_j, \ldots, x_p), \cr } $$
where $\hat x_j$ indicates omission of $x_j$. Note that the continuity of the bracket 
on $\g$ and the action on $V$ imply that $d_\g f$ is continuous. 

We thus obtain a subcomplex of the algebraic Lie algebra complex associated to 
$\g$ and~$V$. Hence $d_\g^2 = 0$, and the space 
$Z^p_c(\g,V) := \ker(d_\g\res_{C^p_c(\g,V)})$
of {\it $p$-cocycles} contains the space 
$B^p_c(\g,V) :=  d_\g(C^{p-1}_c(\g,V))$
of {\it $p$-coboundaries} (cf.\ [We95, Cor.~7.7.3]). The quotient 
$$ H^p_c(\g,V) := Z^p_c(\g,V)/B^p_c(\g,V) $$
is the {\it $p$-th continuous cohomology space of $\g$ with values in
the $\g$-module $V$}. We write $[f] := f + B^p_c(\g,V)$ for the 
cohomology class $[f]$ of the cocycle $f$. 
\qed

\Definition B.2. Let $G$ be a Lie group and $A$ an abelian Lie group. 
We call $A$ a {\it smooth $G$-module} if it is endowed with a
$G$-module structure defined by a smooth action map $G \times A \to
A$. 

Let $A$ be a smooth $G$-module. 
Then we define 
$\tilde C^n_s(G,A)$ to be the space of all functions 
$F \: G^{n+1} \to A$ which are smooth in a neighborhood of the
diagonal, equivariant with respect to the action of $G$ on
$G^{n+1}$ given by 
$$ g.(g_0, \ldots, g_n) := (g g_0, \ldots, g g_n), $$
and vanish on all tuples of the form 
$(g_0, \ldots, g,g, \ldots, g_n).$
As the $G$-action preserves the diagonal, it preserves the space 
$\tilde C^n_s(G,A)$. Moreover, the Alexander--Spanier differential
$d_{AS}$ defines a group homomorphism 
$$ d_{AS} \: \tilde C^n_s(G,A) \to \tilde C^{n+1}_s(G,A), $$
and we thus obtain a differential complex 
$(\tilde C^*_s(G,A), d_{AS}).$

Let $C^n_s(G,A)$ denote the space of all function $f \: G^n \to A$
which are smooth in an identity neighborhood and normalized in the
sense that $f(g_1,\ldots, g_n)$ vanishes if $g_j = \1$ holds for some
$j$. We call these functions {\it normalized locally smooth group
cochains}. Then the map 
$$ \Phi_n \: C^n_s(G,A) \to \tilde C^n_s(G,A), \quad 
\Phi_n(f)(g_0, \ldots, g_n) 
:= g_0.f(g_0^{-1}g_1, g_1^{-1}g_2, \ldots,g_{n-1}^{-1} g_n) $$
is a linear bijection whose inverse is given by 
$$ \Phi_n^{-1}(F)(g_1,\ldots, g_n) := F(\1,g_1, g_1 g_2, \ldots, g_1\cdots g_n).$$
By 
$$ d_G := \Phi_{n+1}^{-1} \circ d_{AS} \circ \Phi_n
\: C^n_s(G,A) \to C^{n+1}_s(G,A) $$
we obtain the differential 
$d_G \: C^n_s(G,A) \to C^{n+1}_s(G,A)$ turning 
$(C^*_s(G,A), d_G)$ into a differential complex. 
We write $Z^n_s(G,A)$ for the corresponding group of cocycles, 
$B^n_s(G,A)$ for the subgroup of coboundaries and 
$$ H^n_s(G,A) := Z^n_s(G,A)/B^n_s(G,A) $$
is called the {\it $n$-th Lie cohomology group with values in the smooth
module $A$}. 
\qed

\Lemma B.3. The group differential $d_G \: C^n_s(G,A) \to
C^{n+1}_s(G,A)$ is given by 
$$ \eqalign{ (d_G f)(g_0, \ldots, g_n) 
&= g_0.f(g_1,\ldots, g_n) \cr
&\ \ \ \ + \sum_{j=1}^n (-1)^j f(g_0, \ldots, g_{j-1} g_j,\ldots, g_n) 
+ (-1)^{n+1} f(g_0, \ldots, g_{n-1}). \cr}$$

\Proof. In fact, $d_{AS} F = \sum_{i=0}^{n+1} (-1)^i F_i$ leads with $F =
\Phi_n(f)$ to 
$d_G f = \sum_{i=0}^{n+1} (-1)^i \Phi_{n+1}^{-1}(F_i)$ and hence to 
$$ \eqalign{ 
&\ \ \ \ (d_G f)(g_0, \ldots, g_{n}) \cr
&= \sum_{i=0}^{n+1} (-1)^i F_i(\1,g_0, g_0 g_1, \ldots, g_0 \cdots
g_{n}) \cr 
&= \sum_{i=0}^{n+1} (-1)^i F(\1,g_0, g_0 g_1, \ldots, 
g_0\cdots g_{i-1}, g_0\cdots g_{i+1}, \ldots, g_0 \cdots g_{n}) \cr 
&= g_0.f(g_1,\ldots, g_{n}) 
+ \sum_{i=1}^n (-1)^i f(g_0, g_1, \ldots, g_i g_{i+1}, \ldots,
g_{n}) + (-1)^{n+1} f(g_0,\ldots, g_{n-1}). \cr} $$ 
\qed

For $n = 0$ we have in particular 
$$ (d_G f)(g_0) = g_0.f - f, $$
and for $n = 1$: 
$$ (d_G f)(g_0,g_1) = g_0.f(g_1) - f(g_0g_1) + f(g_0). $$

\Definition B.4. Let $G$ be a Lie group and 
$\a$ a {\it smooth locally convex $G$-module}, i.e., $\a$ is a locally
convex space and the action map 
$\rho_\a \: G \times \a \to \a, (g,a) \mapsto g.a$ is smooth. We write 
$\rho_\a(g) \: \a \to \a, a \mapsto g.a$ 
for the corresponding continuous linear automorphisms 
of $\a$. 

We call a $p$-form $\Omega \in \Omega^p(G,\a)$ {\it equivariant} if we have 
for all $g \in G$ the relation 
$$ \lambda_g^*\Omega = \rho_\a(g) \circ \Omega. $$
The complex of equivariant differential forms has been introduced in the 
finite-dimensional setting by Chevalley and Eilenberg in [CE48]. 

If $\a$ is a trivial module, then an equivariant $p$-form is 
a left invariant $\a$-valued $p$-form on $G$. 
An equivariant $p$-form is uniquely determined by the corresponding element 
$\Omega_\1 \in C^p_c(\g,\a)$: 
$$ \Omega_g(g.x_1, \ldots, g.x_p) 
= \rho_\a(g) \circ \Omega_\1(x_1, \ldots, x_p), \quad \hbox{ for } \quad
g \in G, x_i \in \g \cong T_\1(G), \leqno(B.1) $$
where $G \times T(G) \to T(G), (g,x) \mapsto g.x$ denotes the 
natural action of $G$ on its tangent bundle $T(G)$ obtained by restricting the 
tangent map of the group multiplication. 

Conversely, (B.1) provides for each $\omega \in
C^p_c(\g,\a)$ a unique 
equivariant $p$-form $\omega^{\rm eq}$ on $G$ with $\omega^{\rm eq}_\1 = \omega$.
\qed

\Lemma B.5. For each $\omega \in C^p_c(\g,\a)$ 
we have $d(\omega^{\rm eq}) = (d_\g\omega)^{\rm eq}$. 
In particular the evaluation map 
$$ \ev_\1 \: \Omega^p(G,\a)^{\rm eq} \to C^p_c(\g,\a), \quad 
\omega \mapsto \omega_\1 $$
defines an isomorphism from the chain complex of equivariant
$\a$-valued differential forms on $G$ to the continuous $\a$-valued Lie algebra
cohomology. 

\Proof. (cf.\ [CE48, Th.~10.1]) For $g \in G$ we have 
$$ \lambda_g^* d\omega^{\rm eq}
= d\lambda_g^* \omega^{\rm eq}
= d(\rho_\a(g) \circ \omega^{\rm eq})
= \rho_\a(g) \circ (d\omega^{\rm eq}), $$
showing that $d\omega^{\rm eq}$ is also equivariant. 

For $x \in \g$ we write $x_l$ for the corresponding
left invariant vector field on $G$, i.e., $x_l(g) =g.x$. 
It suffices to calculate the value of $d\omega^{\rm eq}$ on $(p+1)$-tuples of 
left invariant vector fields in the identity element. 

In view of 
$$ \omega^{\rm eq}(x_{1,l},\ldots, x_{p,l})(g) 
= \rho_\a(g).\omega(x_1, \ldots, x_p), $$
we obtain 
$$ \big(x_{0,l}.\omega^{\rm eq}(x_{1,l},\ldots, x_{p,l})\big)(\1) 
= x_0.\omega(x_1, \ldots, x_p), $$
and therefore 
$$ \eqalign{
\Big(d\omega^{\rm eq}(x_{0,l}, \ldots, x_{p,l}\big)\Big)(\1) 
&=\sum_{i=0}^{p} (-1)^i x_{i,l}.\omega^{\rm eq}(x_{0,l}, \ldots, \hat{x_{i,l}}, 
\ldots, x_{p,l})(\1) \cr
&\ \ \ \ + \sum_{i < j}(-1)^{i+j} \omega^{\rm eq}([x_{i,l}, x_{j,l}], x_{0,l}, 
\ldots,\hat{x_{i,l}}, \ldots, \hat{x_{j,l}}, \ldots, x_{p,l})(\1)\cr           
&=\sum_{i=0}^{p} (-1)^i x_i.\omega(x_0, \ldots, \hat{x_i}, 
\ldots, x_p) \cr
&\ \ \ \ + \sum_{i < j}(-1)^{i+j} \omega([x_i, x_j], x_0, 
\ldots,\hat{x_i}, \ldots, \hat{x_j}, \ldots, x_p) \cr
&=(d_\g\omega)(x_0, \ldots, x_p). \cr} $$
This proves our assertion. 
\qed

\Theorem B.6. The maps 
$$ D_n \: \ev_\1 \circ \tau \circ \Phi_n \: C^n_s(G,A) \to
C^n_c(\g,\a), \quad n \geq 1, $$
induce a morphism of chain complexes 
$$ D \: (C^n_s(G,A), d_G)_{n \geq 1} \to (C^n_c(\g,\a), d_\g)_{n \geq 1} $$
and in particular homomorphisms
$$D_n \: H^n_s(G,A) \to H^n_c(\g,\a), \quad n \geq 2.$$

For $A = \a$ these assertions hold for all $n \in \N_0$ and if 
$A \cong \a/\Gamma_A$ for a discrete subgroup $\Gamma_A$ of
$\a$, then $D_1$ also induces a homomorphism 
$$D_1 \: H^1_s(G,A) \to H^1_c(\g,\a), \quad [f] \mapsto [df(\1)]. $$

\Proof. In view of Proposition A.6 and the definition of the group
differential $d_G$, the composition 
$$ \tau \circ \Phi_n \: C^n_s(G,A) \to \tilde C^n_s(G,A) \subeq
C^n_{AS,s}(G,A) \to \Omega^n(G,\a), \quad n \geq 1,  $$
defines a homomorphism of chain complexes. For $A = \a$ this relation 
also holds for $n = 0$. 

For $f \in C^n_s(G,A)$ the function $F := \Phi_n(f) \: G^{n+1} \to A$ 
is $G$-equivariant with respect to the diagonal action. 
For $g \in G$ let 
$$\mu_g \: G^{n+1} \to G^{n+1}, \quad (g_0,\ldots, g_n) \mapsto 
(gg_0, \ldots, g g_n) $$ 
and write $\rho_A(g)(a) := g.a$ for $a \in A$. 
Then the equivariance of $F$ means that 
$\mu_g^* F = F \circ \mu_g = \rho_A(g) \circ F$ which implies that 
$$ \rho_A(g) \circ \tau(F) 
= \tau(\rho_A(g) \circ F) = \tau(\mu_g^*F) = \lambda_g^* \tau(F). $$
This shows that the image of $\tau \circ \Phi_n$ consists of 
equivariant $\a$-valued $n$-forms on $G$. According to Lemma B.5,
evaluating an equivariant $n$-form in the identity intertwines the
de Rham differential on $\Omega^p(G,\a)$ with the Lie algebra
differential $d_\g$. This implies 
$$ d_\g \circ D_n = D_{n+1} \circ d_G $$
for each $n \in \N$, i.e., the $D_n$ define a morphism of
chain complexes (truncated to $n \geq 1$). 
For $A = \a$ it also holds for $n = 0$. 

If $A \cong \a/\Gamma_A$ and $n = 1$, then $D_1(B^1_s(G,A)) =
B^1_c(\g,\a)$ implies that $D_1$ induces a map $H^1_s(G,A) \to
H^1_c(\g,\a)$. If $A$ is not of this form, then we cannot conclude
that $D_1$ maps $B^1_s(G,A)$ into $B^1_c(\g,\a)$. 
\qed

To make $D_n$, $n \geq 2$, better accessible to calculations, 
we need a more concrete formula for the Lie
algebra cochain $D_n f$ for $f \in C^n_s(G,A)$. 
As $f$ vanishes on all tuples of the form 
$(g_1,\ldots, \1, \ldots, g_n)$, its $(n-1)$-jet in $\1$
vanishes and the term of order $n$ is the $n$-linear map 
$$ (d^n f)(\1,\ldots, \1) \: \g^n = T_\1(G)^n \to \a $$
(cf.\ Definition A.3). In fact, in local coordinates the $n$-th order
term of the Taylor expansion of $f$ in $(\1,\ldots, \1)$ is given by 
a symmetric $n$-linear map  
$$ (d^{[n]}f)(\1,\ldots, \1) \: (\g^n)^n \to \a $$
as 
$$ {1\over n!} (d^{[n]}f)(\1,\ldots, \1)(x,\ldots, x), \quad 
x = (x_1,\ldots, x_n) \in \g^n. $$
The normalization condition on $f$ implies that $(d^{[n]}f)(\1,\ldots,
\1)$ vanishes on all elements \break $(x^1,\ldots, x^n)$, $x^i = (x^i_l)\in \g^n$, 
for which the $j$-th component (in $\g$) vanishes for some $j$, i.e., 
$x^i_j = 0$ for all $i$. This implies that 
$$ (d^{[n]}f)(\1,\ldots, \1)(x,\ldots, x) $$
is a sum of $n!$ terms of the form 
$$ (d^{[n]}f)(\1,\ldots, \1)((0,\ldots, x_{\sigma(1)}, \ldots, 0),
(0, \ldots, x_{\sigma(2)}, \ldots,0), \ldots, (0, \ldots,
x_{\sigma(n)}, \ldots, 0)), $$
since all these terms are equal, we find 
$$ \eqalign{ {1\over n!} (d^{[n]}f)(\1,\ldots, \1)(x,\ldots, x)
&= (d^{[n]}f)(\1,\ldots, \1)((x_1,0,\ldots,0),\ldots, (0,\ldots, 0,x_n))\cr
&= (d^n f)(\1,\ldots, \1)(x_1,\ldots, x_n). \cr}$$

\Lemma B.7. For $f \in C^n_s(G,A)$ and $x_1,\ldots, x_n \in \g$ we
have 
$$ (D_nf)(x_1,\ldots, x_n) 
= \sum_{\sigma \in S_n} \sgn(\sigma)
(d^nf)(\1,\ldots,\1)(x_{\sigma(1)},\ldots, 
x_{\sigma(n)}). $$

\Proof. Recall that on an $n$-tuple $(x_1,\ldots, x_n) \in \g^n$ the
map $d^nf $ can be
calculated by choosing smooth vector fields $X_n$ on an open identity
neighborhood of $G$ with $X_i(\1) = x_i$ via 
$$ (d^n f)(\1,\ldots, \1)(x_1,\ldots, x_n) := (\partial_1(X_1) \cdots
\partial_n(X_n).f)(\1,\ldots, \1). $$
For $F = \Phi_n(f)$ we now get 
$$ \eqalign{ (D_n f)(x_1,\ldots, x_n) 
&= \tau(F)(x_1,\ldots, x_n) 
= \sum_{\sigma \in S_n} \sgn(\sigma) (d^n F)(\1,\ldots,
\1)(0,x_{\sigma(1)}, \ldots, x_{\sigma(n)}). \cr} $$
In view of 
$$ F(\1,g_1,\ldots, g_n) = f(g_1, g_1^{-1} g_2, \ldots,
g_{n-1}^{-1}g_n) $$
and $f(g_1,\1,\ldots) = 0$, we have 
$$ (\partial_1(X_1) F)(\1,\1, g_2,\ldots, g_n) 
= (\partial_1(X_1) f)(\1,g_2,g_2^{-1} g_3, \ldots, g_{n-1}^{-1} g_n),
$$
and inductively we obtain 
$$ \eqalign{ (\partial_1(X_1) \cdots \partial_n(X_n)F)(\1,\1, \ldots, \1) 
&= (\partial_1(X_1) \ldots \partial_n(X_n)f)(\1,\ldots, \1) \cr
&= (d^n f)(\1,\ldots, \1) (x_1,\ldots, x_n). \cr} $$ 
This implies the assertion. 
\qed

For $n = 1$ we obtain 
$(D_1 f)(x) = df(\1).x,$
and for $n = 2$ we have 
$$ \eqalign{ (D_2 f)(x,y) 
&= (d^2 f)(\1,\1)(x,y) - (d^2 f)(\1,\1)(y,x). \cr} $$
If $(d^{[n]} f)(\1,\1)$ denotes the symmetric $n$-linear map 
$(\g^n)^n \to \a$ representing the $n$-jet of $f$, this expression
equals 
$$ (d^{[2]} f)(\1,\1)((x,0)(0,y)) - (d^{[2]} f)(\1,\1)((y,0),(0,x)).  $$

\sectionheadline{Appendix C. Split Lie subgroups} 

\nin In this appendix we collect some general material on Lie group structures on
groups, (normal) Lie subgroups and quotient groups. 
In particular Theorem C.2 provides
a tool to construct Lie group structures on groups for which a subset
containing the identity is an open $0$-neighborhood of a locally convex space
such that the group operations are locally smooth in these
coordinates. We also give a condition on a normal subgroup
$N \trile G$ for the quotient group $G/N$ being a
manifold such that the quotient map $q \: G \to G/N$ defines on $G$
the structure of a smooth $N$-principal bundle. 

\Lemma C.1. Let $G$ be a group and ${\cal F}$ a filter basis of
subsets with $\bigcap {\cal F} = \{\1\}$ satisfying: 
\litem{(U1)} $(\forall U \in {\cal F})(\exists V \in {\cal F}) VV \subeq U. $
\litem{(U2)} $(\forall U \in {\cal F})(\exists V \in {\cal F}) V^{-1} \subeq U. $
\litem{(U3)} $(\forall U \in {\cal F})(\forall g \in G) (\exists V \in
{\cal F}) gVg^{-1} \subeq U. $

\par\nin Then there exists a unique group topology on $G$ such that
${\cal F}$ is a basis of $\1$-neighborhoods in $G$. This topology is
given by 
$\{ U \subeq G \: (\forall g \in U) (\exists V \in {\cal F}) gV
\subeq U\}.$

\Proof. [Bou88, Ch.\ III, \S 1.2, Prop.\ 1]
\qed

\Theorem C.2. Let $G$ be a group and 
$U = U^{-1}$ a symmetric subset. 
We further assume that $U$ is a smooth manifold such that 
\litem{(L1)} there exists an open $\1$-neighborhood $V \subeq
U$ with $V^2 = V \cdot V \subeq U$ such that the group multiplication 
$\mu_V \: V \times V \to U$ is smooth, 
\litem{(L2)} the inversion map $\eta_U \: U \to U, u \mapsto u^{-1}$ is
smooth, and 
\litem{(L3)} for each $g \in G$ there exists an open $\1$-neighborhood $U_g \subeq
U$ with $c_g(U_g) \subeq U$ and such that the conjugation map 
$$ c_g \: U_g \to U, \quad x \mapsto gxg^{-1} $$
is smooth. 

Then there exists a unique structure of a
Lie group on $G$ for which there exists an open $\1$-neighborhood 
$U_1 \subeq U$ such that the inclusion map $U_1 \to G$ induces a 
diffeomorphism onto an open subset of $G$. 

\Proof. (cf.\ [Ch46, \S 14, Prop.~2] or [Ti83, p.14] for the
finite-dimensional case) 
First we consider the filter basis 
$$ {\cal F} := \{ W \subeq G \: W \in {\cal U}_U(\1) \} $$
of all those subsets of $U$ which are $\1$-neighborhoods in $U$. 
Then (L1) implies (U1), (L2) implies (U2), and (L3) implies (U3). 
Moreover, the assumption that $U$ is Hausdorff implies that 
$\bigcap {\cal F} = \{\1\}$. Therefore Lemma C.1 implies that $G$
carries a unique structure of a (Hausdorff) topological group for which 
${\cal F}$ is a basis of $\1$-neighborhoods. 

After shrinking $V$ and $U$, we may assume that there exists a
diffeomorphism $\phi \: U \to \phi(U) \subeq E$, where $E$ is a topological 
$\K$-vector space, $\phi(U)$ an open subset,  
that $V$ satisfies $V = V^{-1}$, $V^4 \subeq U$, and
that $m \: V^2 \times V^2 \to U$ is smooth. For $g \in G$ we consider the maps 
$$ \phi_g \: gV \to E, \quad \phi_g(x) = \phi(g^{-1}x) $$
which are homeomorphisms of $gV$ onto $\phi(V)$. We claim that
$(\phi_g, gV)_{g \in G}$ is an atlas of $G$. 

Let $g_1, g_2 \in G$ and put $W := g_1 V \cap g_2 V$. If $W \not=
\eset$, then $g_2^{-1} g_1 \in V V^{-1} = V^2$. The smoothness of the 
map 
$$ \psi := \phi_{g_2} \circ \phi_{g_1}^{-1}\res_{\phi_{g_1}(W)}
\: \phi_{g_1}(W) \to \phi_{g_2}(W) $$
given by 
$$ \psi(x) 
= \phi_{g_2}(\phi_{g_1}^{-1}(x))
= \phi_{g_2}(g_1 \phi^{-1}(x))
= \phi(g_2^{-1} g_1 \phi^{-1}(x)) $$
follows from the smoothness of the
multiplication $V^2 \times V^2 \to U$. This proves that the charts 
$(\phi_g, gU)_{g \in G}$ form an atlas of $G$. Moreover, the
construction implies that all left translations of $G$ are smooth
maps. 
 
The construction also shows that for each $g \in G$ the conjugation $c_g \: G
\to G$ is smooth in a neighborhood of $\1$. Since all left
translations are smooth, and 
$$ c_g \circ \lambda_x = \lambda_{c_g(x)} \circ c_g, $$
the smoothness of $c_g$ in a neighborhood of  $x \in G$ follows. 
Therefore all conjugations and hence also all
right multiplications are smooth. The smoothness of the inversion
follows from its smoothness on $V$ and the fact that left and right
multiplications are smooth. Finally the smoothness of the
multiplication follows from the smoothness in $\1 \times \1$ because of
$$ \mu_G(g_1 x, g_2 y)= g_1 x g_2 y 
= g_1 g_2 c_{g_2^{-1}}(x) y = g_1 g_2 \mu_G(c_{g_2^{-1}}(x), y). $$ 
The uniqueness of the Lie group structure is clear because each locally
diffeomorphic bijective homomorphism between Lie groups is a
diffeomorphism. 
\qed

\Remark C.3. Suppose that the group $G$ in Theorem C.2 is generated by
each $\1$-neighborhood $U \subeq U$. Then condition (L3) can be
omitted. Indeed, the construction of the Lie group structure 
shows that for each $g \in V$ the conjugation $c_g \: G
\to G$ is smooth in a neighborhood of $\1$. Since the set of all these
$g$ is a submonoid of $G$ containing $V$, it contains $V^n$ for each
$n \in \N$, hence all of $G$ because $G$ is generated by $V$. 
Therefore all conjugations are smooth, and one can proceed as in the
proof of Theorem C.2. 
\qed

\Definition C.4. {\rm(a)} (Split Lie subgroups) Let $G$ be a Lie group. A subgroup $H$ 
is called a {\it split Lie subgroup} if it carries a Lie group
structure for which the canonical right action of 
$H$ on  $G$ defined by restricting the multiplication map of $G$ 
to $G \times H \to G$ defines a smooth principal bundle, i.e., the coset space 
$G/H$ is a smooth manifold and the quotient map $\pi \: G \to G/H$ has smooth local 
sections. 

\par\nin (b) If $G$ is a Banach--Lie group and $\exp \: \g \to G$ its
exponential function, then a closed subgroup $H \subeq G$ is called a
{\it Lie subgroup} if there exists an open $0$-neighborhood 
$U \subeq \g$ such that $\exp\res_U \: U \to \exp(U)$ is a
diffeomorphism onto an open subset of $G$ and the Lie algebra 
$$ \h := \{ x \in \g \: \exp(\R x) \subeq H \} $$
of $H$ satisfies 
$$ H \cap \exp(U) = \exp(U \cap \h). 
\qeddis 

Since the Lie algebra $\h$ of a Lie subgroup $H$ of a Banach Lie group
$G$ need not have a closed complement in 
$\g$, not every Lie subgroup is split. A simple example is the
subgroup $H := c_0(\N,\R)$ in $G := \ell^\infty(\N,\C)$. 

\Lemma C.5. If $H$ is a split Lie subgroup of $G$
or a Lie subgroup of the Banach--Lie group $G$, then for any smooth
manifold $X$ each smooth map $f \:X  \to G$ with $f(X) \subeq H$ is
also smooth as a map $X \to H$. 
If $H$ is a normal split 
Lie subgroup, then the conjugation action of $G$ on $H$ is smooth. 

\Proof. The condition that $H$ is a split Lie subgroup implies that 
there exists an open subset $U$ of some locally convex space $V$ 
and a smooth map $\sigma \: U \to G$ such that the map 
$$ U \times H \to G, \quad (x,h) \mapsto \sigma(x) h $$
is a diffeomorphism onto an open subset of $G$. Let 
$p \: \sigma(U) H \to U$ denote the smooth map given by 
$p(\sigma(x)h) = x$. If $X$ is a manifold and 
$f \: X \to G$ is a smooth map with values in $H$, then $f$ is smooth
as a map to $\sigma(U)H \cong U \times H$, hence smooth as a map $X \to
H$. 

If $H$ is a Lie subgroup of a Banach--Lie group and $f \: X \to G$ is
a smooth map with $f(X) \subeq H$, then we have to see that $f$ is
smooth as a map $X \to H$. To verify smoothness in a neighborhood of
some $x_0 \in X$, it suffices to consider the map $x \mapsto
f(x)f(x_0)^{-1}$, so that we may w.l.o.g.\ assume that 
$f(x_0) = \1$. 
Then we can use the natural chart of $H$ in $\1$ given by the exponential function to
see that $f$ is smooth in a neighborhood of $x_0$ because 
any smooth map $X \to \g$ with values in $\h$ is smooth as a map $X
\to \h$. 

Now suppose that $H \trile G$ is normal. Then the 
conjugation map $G \times H \to G, (g,h) \mapsto ghg^{-1}$, is smooth
with values in $H$, hence smooth as a map $G \times H \to H$. 
\qed 

\Theorem C.6. Let $G$ be a Lie group and $N \trile G$ a split normal 
subgroup. Then the quotient group $G/N$ has a natural Lie group 
structure such that the quotient map 
$q \: G \to G/N$ defines on $G$ the structure of a principal 
$N$-bundle. 

\Proof. There exists an open subset $U$ of a locally convex space $V$ 
and a smooth map $\sigma \: U \to G$ such that the map 
$$ U \times N \to G, \quad (u,n) \mapsto \sigma(u) n $$
is a diffeomorphism onto an open subset $W = \sigma(U)N$ of $G$. 
As $N$ is in particular closed, the quotient group $G/N$ has a natural
(Hausdorff) group topology.

Let $q \: G \to G/N$ denote the quotient map. Then 
$q(W) = q\circ \sigma(U)$ is an open subset of $G/N$ and 
$q(W) \cong W/N \cong (U \times N)/N \cong U$. Therefore the map 
$\phi :=  q \circ \sigma \: U \to q(W)$ is a homeomorphism. 

Let $K = K^{-1} \subeq q(W)$ be a symmetric open subset, 
and $U_K := \phi^{-1}(K)$, and 
endow $K$ with the manifold structure
obtained from the homeomorphism $\phi \: U_K \to K$. 

(L1): Let $V \subeq K$ be an open $\1$-neighborhood with 
$V^2 \subeq K$. We
identify $V$ with the corresponding open subset $U_V \subeq U$. Then
the group multiplication 
$\mu_V \: V \times V \to K$ is given by 
$$ \phi(x)\phi(y) = \sigma(x) N \cdot \sigma(y)N
= \sigma(x)\sigma(y)N= \phi(\phi^{-1}(\sigma(x)\sigma(y)N)), $$
and since the map 
$p \: W \to U, \sigma(u) n \to u$
is smooth, the map 
$$ (x,y) \mapsto \phi^{-1}(\sigma(x)\sigma(y)N)
= p(\sigma(x)\sigma(y)) $$
is smooth. 

(L2): We likewise see that the inversion map $K \to K$ 
corresponds to the smooth map 
$$ x \mapsto \phi^{-1}(\phi(x)^{-1}) 
= \phi^{-1}(\sigma(x)^{-1}N) = p(\sigma(x)^{-1}). $$

(L3): For each $g \in G$ we find an open $\1$-neighborhood $K_g \subeq
K$ with $c_g(K_g) \subeq K$. Then the conjugation map 
$$ c_g \: K_g \to K, \quad x \mapsto gxg^{-1} $$
is written in $\phi$-coordinates as 
$$ \phi(x) 
\mapsto \phi(\phi^{-1}(g\sigma(x) g^{-1} N))
= \phi(p(g\sigma(x) g^{-1})) $$
and therefore smooth. 

Now Theorem C.2 applies and shows that there 
exists a unique structure of a
Lie group on $G/N$ for which there exists an open $0$-neighborhood 
in $U$ such that the map $\phi \: U \to G/N$ 
induces a diffeomorphism onto an open subset of $G/N$. 
\qed

\sectionheadline{Appendix D. The exact Inflation-Restriction Sequence}

In this section $G$ denotes a Lie group, $N \trile G$ a split normal
Lie subgroup (cf.\ Defi\-ni\-tion~C.4) 
and $A$ a smooth $G$-module. We write $q \: G \to G/N$ for the quotient map. 

\Definition D.1. (a) (Inflation and restriction) Restriction of cochains leads for each 
$p \in \N_0$ to a map 
$$ \tilde R \: C^p_s(G,A) \to C^p_s(N,A), $$
and since $R \circ d_G = d_N \circ R$, it follows that 
$\tilde R(B^p(G,A)) \subeq  B^p(N,A)$, 
$\tilde R(Z^p_s(G,A))  \subeq Z^p_s(N,A)$, so that $\tilde R$ induces a homomorphism 
$$ R \: H^p_s(G,A) \to H^p_s(N,A). $$

\par\nin (b) Since $N$ is a normal subgroup of $G$, the subgroup 
$$ A^N := \{ a \in A \: (\forall n \in N)\ n.a = a\} $$
is a $G$-submodule of $A$. If $A^N$ is a split Lie subgroup of
$A$, it inherits a natural structure of a smooth $G/N$-module
(Lemma~C.2) but we do not want to make this restrictive assumption. We
therefore define the chain complex $(C^*_s(G/N,A^N),d_{G/N})$ as the
complex whose cochain space 
$C^p_s(G/N, A^N)$ consists of those functions 
$f \: (G/N)^p \to A^N$ for which the pull-back 
$$q^*f \: G^p \to A^N, \quad 
 (q^*f)(g_1, \ldots, g_p) := f(q(g_1), \ldots, q(g_p)) $$
is an element of $C^p_s(G,A)$. 
With this definition we do not need a Lie group structure on the
subgroup $A^N$ of $A$. 
For a cochain $f \in C^p_s(G/N,A^N)$ we define 
$$ \tilde I := q^* \: C^p_s(G/N,A^N) \to C^p_s(G,A). $$
Then $(C^*_s(G/N,A^N), d_{G/N})$ becomes a chain complex with the
group differential from \break Lemma~B.3. 
Moreover, $q^* \circ d_{G/N} = d_G \circ q^*$, so that 
$q^*(B^p_s(G/N,A^N)) \subeq B^p_s(G,A)$, and 
$q^*(Z^p_s(G/N,A^N))  \subeq Z^p_s(G,A)$, showing that $q^*$ induces 
the so called {\it inflation map} 
$$ I \: H^p_s(G/N,A^N) \to H^p_s(G,A), \quad [f] \mapsto [q^*f]. 
\qeddis 

The restriction and inflation maps 
$$ C^p_s(G/N, A^N) \sssmapright{I} 
C^p_s(G,A) \sssmapright{R} 
C^p_s(N,A) $$ 
clearly satisfy $R \circ I = 0$, which is inherited by the corresponding maps 
$$ H^p_s(G/N,A^N) \sssmapright{I} H^p_s(G,A) \sssmapright{R} H^p_s(N,A).$$

\Lemma D.2. The restriction maps 
$\tilde R \: C^p_s(G,A) \to C^p_s(N,A)$
are surjective. 

\Proof. Since $N$ is a split Lie subgroup of $G$, there exists an open
$0$-neighborhood $U$ in a locally convex space $V$ and a smooth map 
$\phi \: U \to G$ with $\phi(0) = \1$ such that the map 
$$ \Phi \: N \times U \to G, \quad (n,x) \mapsto n \phi(x) $$
is a diffeomorphism onto an open subset $N \phi(U)$ of $G$. 

Let $f \in C^p_s(N,A)$. We extend $f$ to a function 
$\tilde f \: (N \phi(U))^p \to A$ by 
$$ \tilde f((n_1\phi(x_1), \ldots, n_p \phi(x_p)) := f(n_1,\ldots,
n_p). $$
Then clearly $\tilde f$ is smooth in an identity neighborhood and
vanishes if one argument $n_i \phi(x_i)$ is $\1$, because this implies
$x_i = 0$ and $n_i= \1$. Now we extend $\tilde f$ to a function on 
$G^p$ vanishing in all tuples $(g_1, \ldots, \1, \ldots, g_p)$. Then 
$\tilde f\in C^p_s(G,A)$ satisfies $\tilde R(\tilde f) = f$. 
\qed

Although the the inflation map $I$ is injective on cochains and $R$ is surjective 
on cochains, in general there are many cochains with trivial restrictions on 
$N$ which are not in the image of the inflation map. Therefore we do not have 
a short exact sequence of chain complexes, hence cannot expect a long exact sequence 
in cohomology. In this appendix we discuss what we still can say on the corresponding 
maps in low degree. It would be interesting to see if these results 
can also be obtained from a generalization of the Hochschild--Serre 
spectral sequence for Lie groups. As we shall see below, it is clear that the 
construction in [HS53a] has to be modified substantially for the locally smooth 
infinite-dimensional setting. 

\Lemma D.3. {\rm(a)} Each cohomology class in $H^p_s(G,A)$ annihilated
by $R$ can be represented by a cocycle in $\ker \tilde R$. 

\par\nin {\rm(b)} We have $B^p_s(N,A) \subeq \im(\tilde R)$ and 
therefore $[f] \in \im(R)$ is equivalent to $f \in \im(\tilde R)$. 

\Proof. (a) We may w.l.o.g.\ assume that $p \geq 1$. 
If $R[f] = 0$, then $\tilde R(f) = d_N \alpha$ for some 
$\alpha \in C^{p-1}_s(N,A)$. Let $\tilde\alpha \in C^{p-1}_s(G,A)$ be
an extension of $\alpha$ to $G$ (Lemma D.2). Then 
$f' := f - d_G \tilde\alpha$ restricts to 
$\tilde R(f) - d_N \alpha = 0$ and $[f'] = [f]$. 

\par\nin (b) For $\alpha \in C^{p-1}_s(G,A)$ we have 
$\tilde R(d_G \alpha) = d_N \tilde R(\alpha)$, so that 
$C^{p-1}_s(N,A) \subeq \im(\tilde R)$ implies that 
$\tilde R(B^p_s(G,A)) = B^p_s(N,A)$. 

For $f \in Z^p_s(N,A)$ it follows that 
$[f] \in \im(R)$ is equivalent to the existence of 
$\alpha \in B^{p-1}_s(N,A)$ with $f - d_N \alpha \in \im(\tilde R)$,
which implies that $f \in \im(\tilde R)$. 
\qed

\Lemma D.4. The coboundary operator $d_N$ is equivariant with respect to the action of 
$G$ on $C^p_s(N,A)$, $p \in \N_0$, given by 
$$ (g.f)(n_1, \ldots, n_p) 
:= g.f(g^{-1}n_1 g^{-1}, \ldots, g^{-1} n_p g). $$
In particular, this action leaves the space of cochains invariant and induces actions 
on the cohomology groups $H^p_s(N,A)$.  
\qed

The preceding lemma applies in particular to the case $N = G$, showing that 
the coboundary operator $d_G$ is equivariant for the natural action of $G$ on 
the spaces $C^p_s(G,A)$. 

\Definition D.5. In the following we need a refined concept of invariance of cohomology
classes in $H^p_s(N,A)$ under the action of the group $G$. 
We call $f \in Z^p_s(N,A)$ {\it smoothly cohomologically invariant}
if there exists a map 
$$ \theta \: G \to C^{p-1}_s(N,A) 
\quad \hbox{ with } \quad 
d_N(\theta(g)) = g.f - f \quad \hbox{ for all} \quad g \in G $$
for which the map 
$$ G \times N^p \to A, \quad (g,n_1,\ldots, n_{p-1}) \to 
\theta(g)(n_1,\ldots, n_{p-1}) $$
is smooth in an identity neighborhood of $G \times N^{p-1}$. 

We write $Z^p_s(N,A)^{[G]}$ for the set of smoothly cohomologically
invariant cocycles in the group $Z^p_s(N,A)$. If $f = d_N h$ for some 
$h \in C^{p-1}_s(N,A)$, then we may put 
$\theta(g) := g.h - h$ to find 
$$d_N(\theta(g)) = d_N(g.h- h) = g.d_N(h) - d_N(h) = g.f - f, $$
and the map 
$$ \eqalign{ G \times N^{p-1} 
&\to A, \cr
(g,n_1,\ldots, n_{p-1}) &\mapsto (g.h-h)(n_1,\ldots, n_{p-1}) 
= g.h(g^{-1}n_1 g, \ldots, g^{-1}n_{p-1}g) - h(n_1,\ldots, n_{p-1})
\cr} $$
is smooth in an identity neighborhood. This shows that 
$B^p_s(N,A) \subeq Z^p_s(N,A)^{[G]}$, and we define the space of 
{\it smoothly invariant cohomology classes} by 
$$H^p_s(N,A)^{[G]} := Z^p_s(N,A)^{[G]}/B^p_s(N,A). 
\qeddis 

For a generalization of the following fact to general $p$ for discrete
groups and modules we refer to [HS53a]. 

\Proposition D.6. Let $N \trile G$ be a split normal Lie
subgroup and $p \in \{0,1,2\}$. Then 
the restriction map $R$ maps $H^p_s(G,A)$ into $H^p_s(N,A)^{[G]}$. 
In particular 
$$ H^p_s(G,A) = H^p_s(G,A)^{[G]} \quad \hbox{ for } \quad
p=0,1,2. \leqno(D.1) $$

\Proof. In view of the $G$-equivariance of the restriction map 
$C^p_s(G,A) \to C^p_s(N,A)$, it suffices to prove the assertion in the
case $N = G$. 

For $p =0$ we have $C^0_s(G,A) = A$, and $Z^0_s(G,A) =
H^0_s(G,A) = A^G$ is the submodule of 
$G$-invariants. Clearly $G$ acts trivially on this space, so that
there is nothing to prove. 

For $p = 1$ and a cocycle $f \in Z^1_s(G,A)$ we have for $g,x \in G$: 
$$ \eqalign{ (g.f-f)(x) 
&= g.f(g^{-1}xg) - f(x) 
= g.(g^{-1}.f(xg) + f(g^{-1})) - f(x) 
= f(xg) + g.f(g^{-1}) - f(x) \cr
&= x.f(g) + f(x) - f(g) - f(x) 
= d_G(f(g))(x). \cr} $$
This shows that 
$$g.f - f = d_G(f(g)), \leqno(D.2) $$ 
so that $f \in Z^2_s(G,A)^{[G]}$ follows from the local smoothness of
$f$. 

For $p = 2$ and $f \in Z^2_s(G,A)$ we have 
$$ \eqalign{ 
&\ \ \ \ (g.f - f)(x,x') \cr
&= g.f(g^{-1}xg,g^{-1}x'g) - f(x,x') \cr
&= - f(g,g^{-1}xx'g) + f(g,g^{-1}xg) + f(xg, g^{-1}x'g) - f(x,x') \cr 
&= - f(g,g^{-1}xx'g) + f(g,g^{-1}xg) - f(x,g) + x.f(g,g^{-1}x'g) +
f(x,x'g) - f(x,x')  \cr
&= - f(g,g^{-1}xx'g) + f(g,g^{-1}xg) - f(x,g) + x.f(g,g^{-1}x'g) -
x.f(x',g)  + f(xx',g)\cr} $$
and the function 
$$ \theta(g) \: G  \to A, \quad 
\theta(g)(x) := f(g,g^{-1}xg) - f(x,g) $$
satisfies 
$$ \eqalign{ (d_G \theta(g))(x,x') 
&= x.\theta(g)(x') + \theta(g)(x) - \theta(g)(xx') \cr
&= x.f(g,g^{-1}x'g) - x.f(x',g) 
+f(g,g^{-1}xg) - f(x,g) - f(g,g^{-1}xx'g) + f(xx',g) \cr
&= (g.f - f)(x,x').\cr} $$
Since the function $G^2 \to A, (g,x) \mapsto \theta(g)(x)$ 
is smooth in an identity neighborhood of $G^2$, the assertion follows
for $p =2$. 
\qed

\Lemma D.7. For each $f \in Z^1_s(N,A)^{[G]}$ there exists 
$a \in C^1_s(G,A)$ with 
$$ d_N(a(g)) = g.f - f, \quad a(gn) = a(g) + g.f(n), \quad g \in G, n \in N. $$
Then $d_G a \in B^2_s(G,A)$ is $A^N$-valued and constant on 
$(N \times N)$-cosets, hence factors to a cocycle 
$\oline{d_G a} \in Z^2_s(G/N,A^N)$. 
The cohomology class $[\oline{d_G a}]$ does not depend on the choice of $f$ in $[f]$ and 
the function $a$, and we thus obtain a group homomorphism 
$$ \delta \: H^1_s(N,A)^{[G]} \to H^2_s(G/N,A^N), \quad [f] \mapsto
[\oline{d_G a}]. $$

\Proof. Since $N$ is a split Lie subgroup, there exists an open
$0$-neighborhood of some locally convex space $V$ and a smooth map 
$\phi \: U \to G$ with $\phi(0) = \1$ such that the multiplication map 
$$N \times U \to G, (x,n) \mapsto \phi(x)n $$
is a diffeomorphism onto an open subset of $G$. 
Let $E \subeq G$ be a set of representatives of 
the $N$-cosets containing $\phi(U)$, 
so that the multiplication map $E \times N \to G$ is bijective. 

The requirement $f \in Z^1_s(N,A)^{[G]}$ implies the existence of a
function $\alpha \in C^1_s(G,A)$ with $d_N(\alpha(g)) = g.f - f$. We now define 
$$ a \: G = EN \to A, \quad x\cdot n \mapsto \alpha(x) + x.f(n). $$
Then $a$ is smooth on an identity neighborhood because $E$ contains
$\phi(U)$. 
Since $f$ is a $1$-cocycle, we have for $x \in E$ and $n, n' \in N$ the relation 
$$ a(xnn') = a(x) + x.f(nn') = a(x) + x.f(n) + (xn).f(n') = a(xn) + (xn).f(n'), $$
which means that 
$$ a(gn) = a(g) + g.f(n), \qquad g \in G, n \in N. $$
In view of (D.2), we have for $n \in N$ the relation 
$n.f - f = d_N(f(n))$, so that 
$$ (xn).f - f = x.(n.f - f) + x.f -f 
= x.d_N(f(n)) + d_N(a(x)) = d_N(x.f(n) + a(x)) = d_N(a(xn)), $$
and hence 
$d_N(a(g)) = g.f -f$
for all $g \in G$. 

That the values of the function $d_G a$ lie in $A^N$ follows from 
$$ \eqalign{ d_N(a(g_1 g_2)) 
&= (g_1 g_2).f - f 
= g_1.(g_2.f - f) + g_1.f - f \cr
&= g_1.d_N(a(g_2)) + d_N(a(g_1)) = d_N(g_1.a(g_2) + a(g_1))\cr} $$
in $C^1_s(N,A)$. The coboundary $d_G a$ is a cocycle, hence an element
of $Z^2_s(G,A^N)$. 
We show that $d_G a$ is constant on the cosets of $N$. 
We have 
$$ \eqalign{ (d_G a)(g_1, g_2 n) 
&= g_1.a(g_2n) + a(g_1) - a(g_1 g_2n) \cr
&= g_1.a(g_2) + g_1 g_2.f(n) + a(g_1) - a(g_1 g_2) - g_1 g_2.f(n) 
= (d_G a)(g_1, g_2) \cr} $$
and 
$$ \eqalign{ (d_G a)(g_1 n, g_2) 
&= g_1 n.a(g_2) + a(g_1 n) - a(g_1 n g_2) \cr
&= g_1 n.a(g_2) + a(g_1) + g_1.f(n) - a(g_1 g_2 (g_2^{-1} n g_2)) \cr
&= g_1 n.a(g_2) + a(g_1) + g_1.f(n) - a(g_1 g_2) - (g_1 g_2).f(g_2^{-1} n g_2) \cr
&= g_1 n.a(g_2) + a(g_1) + g_1.f(n) - a(g_1 g_2)  - g_1.((g_2.f)(n)) \cr
&= (d_G a)(g_1, g_2) 
+ g_1.(n a(g_2) - a(g_2)) + g_1.f(n) -g_1.f(n) - g_1.(n.a(g_2) - a(g_2))\cr
&= (d_G a)(g_1, g_2) \cr} $$

We now define 
$$ \oline{d_G a} \: G/N \times G/N \to A^N, \quad (xN,yN) \mapsto 
(d_G a)(x,y). $$
Since $d_G a$ is a cocycle on $G$, the function $\oline{d_G a}$ 
is an element of $Z^2_s(G/N,A^N)$. 
It remains to show that the cohomology class of $\oline{d_G a}$ in
$H^2_s(G/N,A^N)$  
does not depend on the choices of $a$ and $f$. 
If $a' \in C^1_s(G,A)$ is another function with 
$$ d_N(a'(g)) = g.f - f, \quad a'(gn) = a'(g) + g.f(n), \quad g \in G, n \in N, $$
then $d_N(a'(g) - a(g)) = 0$ implies that 
$$ \beta(g) := a'(g) - a(g) \in A^N, \quad g \in G. $$
Moreover, 
$$ \beta(gn) = a'(gn) - a(gn) 
= a'(g) + g.f(n) - a(g) - g.f(n) 
= a'(g) - a(g) = \beta(g), $$
so that $\beta$ factors through a function $\gamma \: G/N \to A^N$, and we have 
$$ (d_{G/N}\gamma)(xN,yN) = x.\beta(y) - \beta(xy) + \beta(x) = (d_G\beta)(x,y)
= (d_G a - d_G a')(x,y). $$
Moreover, the fact that the quotient map $G \to G/N$ defines on $G$
the structure of a smooth $N$-principal bundle implies that 
$\gamma$ is smooth in an identity neighborhood of $G/N$. 
Hence the cocycle $\oline{ d_G a'}$ is an element of $Z^2_s(G/N,A^N)$
and satisfies $\oline{d_G a'} = \oline{d_G a} - d_{G/N}\gamma$, 
so that $[\oline{d_G a}] = [\oline{d_G a'}]$. 

Now suppose that $f' \in Z^1_s(N,A)$ satisfies $f' = f + d_N c$ for some $c \in A$. 
In view of the $G$-equivariance of the differential $d_N$, we have 
$$ g.(d_N c) - d_N c = d_N(g.c-c) \quad \hbox{ and } \quad 
(d_G c)(gn) = (d_G c)(g) + g.((d_G c)(n)), $$
so that the function 
$a' := a + d_G c$
satisfies 
$$ d_N(a'(g)) 
= d_N(a(g) + g.c - c) = g.f - f + g.d_N(c) - d_N(c) 
= g.f' - f', \quad 
a'(gn) = a'(g) + g.f'(n). $$
As $d_G c$ is  a cocycle, we have 
$d_G a' = d_G a$, so that we obtain in particular the same cocycles on $G/N$. 
\qed

With the preceding lemma, we can prove the exactness of the 
Inflation-Restriction Sequence: 

\Proposition D.8. Let $A$ be a smooth $G$-module 
and $N \trile G$ a split normal Lie subgroup. Then 
we have the following exact Inflation-Restriction Sequence: 
$$ \0 \to H^1_s(G/N,A^N) \sssmapright{I} H^1_s(G,A) \sssmapright{R}
H^1_s(N,A)^{[G]} \sssmapright{\delta} H^2_s(G/N,A^N) \sssmapright{I} H^2_s(G,A). $$

\Proof. (see [We95, 6.8.3] or [MacL63, pp.347--354] for the case of
abstract groups)   

\par\nin {\bf Exactness in $H^1_s(G/N,A^N)$:}  Let $\alpha \in Z^1_s(G/N, A^N)$. 
We have $[q^*\alpha] =0$ if and only if there exists an $a \in A$ with 
$\alpha(gN) = g.a -a$ for all $g \in G$. That this function is
constant on $N$-left cosets implies that 
$a \in A^N$, and hence that $\alpha = d_{G/N}a \in B^1_s(G/N,A^N)$. 
Therefore the inflation map $I$ is injective on $H^1_s(G/N, A^N)$. 

{\bf Exactness in $H^1_s(G,A)$:}  
That the restriction map $\tilde R$ maps into smoothly $G$-invariant cohomology classes 
follows from Proposition D.6 and the $G$-equivariance of $R$. 
The relation $R \circ I = 0$ is clear. 

To see that $\ker R \subeq \im I$,
let $f \in Z^1_s(G,A)$ vanishing on~$N$ (Lemma D.3). 
Then $f$ is constant on 
the $N$-cosets because 
$$ f(gn) = f(g) + g.f(n) = f(g), \quad g \in G, n \in N. $$
Moreover, 
$$ n.f(g) = f(ng) - f(n) = f(ng) = f(gg^{-1}ng) = f(g) $$
implies that $\im(f) \subeq A^N$. Hence $[f]$ is contained in the image of the 
inflation map $I$. 

{\bf Exactness in $H^1_s(N,A)^{[G]}$:} If $f \in Z^1_s(N,A)$ 
is the restriction of a $1$-cocycle 
$\alpha \in Z^1_s(G,A)$, then (D.2) implies 
$$ (g.f -f)(n) = (d_N(\alpha(g)))(n), $$
so that we may take $\alpha$ as the function $a$ in the definition of $\delta$. 
Then $d_Ga = d_G\alpha = 0$ because $\alpha$ is a cocycle, and hence 
$\delta([f]) = 0$. 

If, conversely, $\delta([f]) = 0$, then there exists $b \in
C^1_s(G/N,A^N)$ with $\oline{d_G a} = d_{G/N} b$, 
where $\oline{d_G a}(xN,yN) = (d_G a)(x,y)$ is defined as in Lemma~D.7. 
Then the function $a' := a - (b \circ q)$ satisfies 
$$ a'(gn) = a'(g) + g.f(n), \quad d_N(a'(g)) = g.f - f, \quad g \in G, n \in N, $$
and, in addition, 
$$ d_G a' = d_Ga - d_G(q^*b) = d_G a - q^*(d_{G/N} b) = q^*(\oline{d_G
a} - d_{G/N}b) = 0. $$
This means that $a' \in Z^1_s(G,A)$, so that $a'\res_N = a\res_N = f$ implies that 
$[f]$ is in the image of the restriction map $R$. 

{\bf Exactness in $H^2_s(G/N,A^N)$:} If $f \in Z^1_s(N,A)$ 
has a smoothly invariant cohomology class and 
$[\oline{d_G a}] = \delta([f])$ as in Lemma D.7, then the image of
$[\oline{d_G a}]$ in 
$Z^2_s(G,A)$ under $I$ is given by $d_G a = q^*\oline{d_G a}$, hence a coboundary. 

Suppose, conversely, that for $\alpha \in Z^2_s(G/N,A^N)$ the cocycle 
$q^*\alpha$ on $G$ is a coboundary and $\beta \in C^1_s(G,A)$
satisfies $q^*\alpha = d_G\beta$. Then $d_G\beta$ vanishes on $N$, so that 
$f := \beta\res_N$ is a cocycle. We have 
$$ \alpha(xN, yN) = x.\beta(y) - \beta(xy) +  \beta(x), \quad x,y \in G. $$
For $y \in N$ we obtain from $\alpha(xN,N) = \alpha(N, xN) = \{0\}$ the relation 
$$ \beta(gn) = \beta(g) + g.\beta(n) 
\quad \hbox{ and } \quad  \beta(ng) = \beta(n) + n.\beta(g). $$
For $g \in G$ and $n \in N$ we therefore have 
$$ \eqalign{ (g.f - f)(n) 
&= g.\beta(g^{-1}ng) - \beta(n) 
= \beta(ng) - \beta(g) - \beta(n) \cr
&= \beta(n) + n.\beta(g) - \beta(g) - \beta(n) 
=  n.\beta(g) - \beta(g) = d_N(\beta(g))(n). \cr} $$
This means that $[f]$ is smoothly $G$-invariant and that 
$\delta([f]) = [\alpha].$ 
\qed 

\Example D.9. The following example shows that the exact Inflation-Restriction 
sequence cannot be continued in an exact fashion by the restriction map 
$R \: H^2_s(G,A) \to H^2_s(N,A)^{[G]}$.

For that we consider the group $G := \R^2$, $N := \Z^2$, $G/N = \T^2$ and the 
trivial module $A = \T  = \R/\Z$. 
Then 
$$ H^2_s(G/N, A^N) = H^2_s(\T^2, \T) = \{0\}, \quad 
H^2_s(G, A) = H^2_s(\R^2, \T) \cong H^2_c(\R^2, \R) \cong \R, $$
and $H^2_s(N,A)^{[G]} = H^2(\Z^2,\T) \cong \T$. Now the assertion follows from the fact that 
the natural map 
$R \: H^2_s(\R^2,\T) \cong \R\to H^2_s(\Z^2,\R) \cong \T$ is not injective. 
It corresponds to restring an alternating $\T$-valued bilinear form to 
the lattice $\Z^2$. If the form is integral on this lattice, the corresponding 
extension of $\Z^2$ is abelian, hence trivial. 
\qed

\Remark D.10. If 
$A$ is a trivial $G$-module, then the connecting map has a simpler description. 
Then we have $H^1_s(N,A) = \Hom(N,A) = Z^1_s(N,A)$, and the condition that a homomorphism 
$f \: N \to A$ is invariant under $G$ means that 
it vanishes on the normal subgroup $[G,N]$ of $N$. 

The only condition on the function 
$a \: G \to A$ that we need to describe $\delta$ is 
$$ a(gn) = a(g) + f(n), \quad g \in G, n \in N. $$
Then the function 
$(d_G a)(x,y) = a(y) - a(xy) + a(x)$
is constant on $(N \times N)$-cosets and defines a $2$-cocycle in 
$Z^2_s(G/N,A)$. 
\qed

\Example D.11. (a) If $G$ is a Lie group, then its identity component
$G_0$ is a split normal subgroup and the quotient group $\pi_0(G)$ is
discrete. Therefore the Inflation-Restriction Sequence yields an exact sequence 
$$ \0 \to H^1(\pi_0(G),A^{G_0}) \sssmapright{I} H^1_s(G,A) \sssmapright{R}
H^1_s(G_0,A)^{[G]} \sssmapright{\delta} H^2(\pi_0(G),A^{G_0}) 
\sssmapright{I} H^2_s(G,A). $$

\par\nin (b) Assume that $A \cong \a/\Gamma_A$ for a discrete subgroup $\Gamma_A$ 
of the sequentially complete locally convex space $\a$. If $G$ is a connected Lie group, 
$q_G \: \tilde G \to G$
its universal covering and $\pi_1(G)$ its kernel, then 
$\pi_1(G)$ is discrete, hence a split Lie subgroup, and we obtain
for any smooth $G$-module $A$ the exact sequence 
$$ \0 \to H^1_s(G, A) \sssmapright{I} H^1_s(\tilde G,A) \sssmapright{R}
H^1_s(\pi_1(G),A)^{[G]} \sssmapright{\delta} H^2_s(G,A)
\sssmapright{I} H^2_s(\tilde G,A). $$
As $\pi_1(G)$ acts trivially on $A$ and $\pi_1(G)$ is central in
$\tilde G$, we have 
$$ H^1_s(\pi_1(G),A) = \Hom(\pi_1(G), A),  \quad 
H^1_s(\pi_1(G),A)^{[G]}  = H^1_s(\pi_1(G),A)^G  = \Hom(\pi_1(G),A^G).$$
In view of Corollary VII.3, we may identify 
$H^2_s(\tilde G,A)$ with the subgroup $\ker P_1$ of $H^2_c(\g,\a)$. 
On this subgroup the map $[\omega] \mapsto F_\omega$ 
given by the flux homomorphism 
defines a homomorphism 
$$ \tilde P_2 \: H^2_s(\tilde G,A) \to \Hom(\pi_1(G), H^1_c(\g,\a)) 
\cong  \Hom(\pi_1(G), H^1_s(\tilde G,A)) $$
whose kernel coincides with the image of $I$ (Theorem~VII.2). 
In Remark VI.8 we have seen that the image of 
$[\omega] \in H^2_s(\tilde G,A) \subeq H^2_c(\g,\a)$ in 
$H^2_s(\pi_1(G),A)$ is given by the commutator map of the corresponding 
central extension 
$$ C([\gamma], [\eta]) = - P(F_\omega([\gamma]))([\eta]),$$
where $P$ is defined in Proposition~III.4. From Example D.9 we know that 
the vanishing of $C$ does not imply the vanishing of $F_\omega$. 

Another interesting aspect of this observation is that, according to a result of 
H.~Hopf, there is an exact sequence 
$$ \0 \to H^2(\pi_1(G), A) \to 
H^2_{\rm sing}(G,A) \cong \Hom(H_2(G), A) \to \Hom(\pi_2(G), A) \to \0 $$
(cf.\ [ML78, p.5]). If $G$ is smoothly paracompact, then de Rham's Theorem 
holds ([KM97]) 
and the closed $2$-form $\omega^{\rm eq}$ defines a singular cohomology 
class in $H^2_{\rm sing}(G,\a) \cong \Hom(H_2(M), \a)$ and after composition 
with the map $q_A \: \a \to A$ a singular cohomology class $c_\omega \in H^2_{\rm sing}(G,A)$. 
The inclusion $\Pi_\omega \subeq \Gamma_A$ means that this class vanishes on 
the spherical cycles, i.e., the image of $\pi_2(G)$ in $H_2(G)$. Hence it determines 
a central extension of $\pi_1(G)$ by $A$ which is given by the commutator map 
$C \: \pi_1(G)^2 \to A$. If this map vanishes, then $c_\omega = 0$,
but Example D.9 shows that this does not imply the existence of a corresponding 
global group cocycle. If $G$ is simply connected, then $c_\omega$ vanishes if and 
only if $\omega$ integrates to a group cocycle, but in general this simple criterion 
fails. 
\qed

\Remark D.12. Let $f_N \in Z^1_s(N,A)^{[G]}$ and $f \in C^1_s(G,A)$ with
$$ f(gn) = f(g) + g.f_N(n), \quad d_N(f(g)) = g.f_N - f_N, \quad g \in G, n
\in N. $$
Then $\delta(f_N) = [\oline{d_{G} f}] \in Z^2_s(G/N,A^N)$ defines an
abelian extension of $G/N$ by $A^N$. We now describe this abelian
extension directly in terms of $f_N$. Here we assume that $A^N$ is a Lie
group and that any smooth map $X \to A$ with values in $A^N$ defines
a smooth map $X \to A^N$ (cf.\ Appendix C). 

Using the smooth action of $G$ on $A$, we can form the semi-direct
product Lie group $A \rtimes G$. Then we consider the map 
$$ \sigma \: G \to A \rtimes G, \quad g \mapsto (f(g),g). $$
In view of $f\res_N = f_N \in Z^1_s(N,A)$, the restriction 
$\sigma\res_N$ is a
homomorphism. Moreover, for $g,g' \in G$ we have 
$$ \sigma(g) \sigma(g') = (f(g) + g.f(g'), gg') 
\quad \hbox{ and } \quad 
\sigma(gg') = (f(gg'), gg'), $$
which implies that 
$$ \delta_\sigma(g,g') := \sigma(g)\sigma(g') \sigma(gg')^{-1} 
= ( (d_G f)(g,g'), \1) \in A^N \times \{\1\}. $$
Therefore the induced map 
$\oline\sigma \: G \to (A/A^N) \rtimes G$
is a group homomorphism, and the pull-back of the abelian extension 
$$  A^N \into A \rtimes G \onto (A/A^N) \rtimes G $$
is isomorphic to the abelian extension 
$\hat G := A^N \times_{d_G f} G$ defined by
$d_G f \in Z^2_s(G,A^N)$. Since $f$ vanishes on $N \times G$ and $G \times N$, 
the subset $\{0\} \times N$ is a normal subgroup of $\hat G$, and 
$\hat G/N \cong A^N \times_{\oline{d_G f}} G/N$. 
\qed

\sectionheadline{Appendix E. A long exact sequence for Lie group
  cohomology} 

Let $G$ be a Lie group and 
$$ \0 \to A_1 \sssmapright{q_1} A_2 \sssmapright{q_2} A_3 \to \0$$ 
be an extension of abelian Lie groups which
are smooth $G$-modules such that $q_1$ and $q_2$ are
$G$-equivariant. We assume that there exists a smooth section 
$\sigma \: A_3 \to A_2$ of $q_2$. Then the map 
$$ A_1 \times A_3  \to A_2, \quad (a,b) \mapsto a + \sigma(b) $$
is a diffeomorphism (not necessarily a group homomorphism). 
This assumption implies that the natural maps 
$$ C^p_s(G, A_1) \to C^p_s(G, A_2) \to C^p_s(G,A_3) $$
define a short exact sequence of chain complexes, hence induce a long
exact sequence in cohomology 
$$ \0 \to H^0_s(G,A_1) \to H^0_s(G,A_2) \to H^0_s(G,A_3) \to H^1_s(G,
A_1) \to \ldots $$
$$ \ldots \to H^{p-1}_s(G,A_3) 
\sssmapright{\delta} H^p_s(G,A_1) \to H^p_s(G,A_2) \to H^p_s(G,A_3) 
\sssmapright{\delta} H^{p+1}_s(G,A_1) \to \ldots $$
The connecting map $\delta \: H^p_s(G,A_3) \to H^{p+1}_s(G,A_1)$ is
constructed as follows. For $f \in Z^p_s(G,A_3)$ we first find 
$f_1 \in C^p_s(G,A_2)$ with $f = q_2 \circ f_1$. Then 
$0 = d_G f = q_2 \circ d_G f_1$ implies that 
$d_G f_1$ is $A_1$-valued, hence an element of $Z^{p+1}_s(G, A_1)$,
and then $\delta([f]) = [d_G f_1]$. 

For $p = 0$ we have $H^0_s(G,A) = A^G$, so that the exact sequence
starts with 
$$ A_1^G \into A_2^G \to A_3^G \to H^1_s(G, A_1) \to H^1_s(G, A_2) \to
\ldots. $$

\Remark E.1. A particularly interesting case arises if $A$ is a smooth $G$-module, 
$A_0$ its identity component and $\pi_0(A) := A/A_0$. Then $\pi_0(A)$
is discrete. Let us assume, in addition, that $G$ is connected. Then
$G$ acts trivially on the discrete group $\pi_0(A)$. We therefore
have an exact sequence 
$$ A_0^G \into A^G \to \pi_0(A) \sssmapright{\oline\theta_A} 
H^1_s(G, A_0) \to H^1_s(G, A) \to
H^1_s(G, \pi_0(A)) = \0, $$
where we use $Z^1_s(G,\pi_0(A)) \subeq C^\infty(G, \pi_0(A)) = \0$ 
(Lemma III.1) to see that $H^1_s(G,\pi_0(A))$ is trivial. 
Note that $\oline\theta_A$ is the characteristic homomorphism of the smooth 
$G$-module $A$, considered as a map into $H^1_s(G,A_0)$ which we may 
consider as a subspace of $H^1_c(\g,\a)$ (Definition~III.6). 
It follows in particular that the natural map 
$H^1_s(G, A_0) \to H^1_s(G, A)$
is surjective. 

Moreover, we obtain an exact sequence 
$$ \0 \to H^2_s(G,A_0) \to H^2_s(G,A) \to H^2_s(G,\pi_0(A)) 
\sssmapright{\delta} H^3_s(G,A_0) \to\ldots $$
Since $G$ is connected and $\pi_0(A)$ is a trivial module, 
the group $H^2_s(G,\pi_0(A))$ classifies the central extensions of $G$
by $\pi_0(A)$, which is parametrized by the abelian group 
$\Hom(\pi_1(G), \pi_0(A))$ (Theorem VII.2). 
This leads to an exact sequence 
$$ \0 \to H^2_s(G,A_0) \to H^2_s(G,A) \sssmapright{\gamma} \Hom(\pi_1(G), \pi_0(A)) \to
H^3(G, A_0), \leqno(E.1) $$
where $\gamma$ assigns to an extension of $G$ by $A$ the corresponding connecting 
homomorphism $\pi_1(G) \to \pi_0(A)$ in the long exact homotopy sequence. 
With the results of Section VII we have determined $H^2_s(G,A_0)$
in terms of the topology of $G$ and the Lie algebra cohomology space 
$H^2_c(\g,\a)$. To determine $H^2_s(G,A)$ in terms of $H^s(G,A_0)$ and
known data, one has to determine the image of $H^2_s(G,A)$ in 
$\Hom(\pi_1(G), \pi_0(A))$. Recall that Proposition~VI.4 shows that 
$$ F_{Df} = - \oline\theta_A \circ \gamma([f]) $$
holds for each $f \in Z^2_s(G,A)$. 

If $A$ is a trivial $G$-module, then the divisibility of $A_0$ implies
that $A \cong A_0 \times \pi_0(A)$ as Lie groups, and we thus obtain 
$$ H^2_s(G, A) 
\cong H^2_s(G,A_0) \times H^2_s(G, \pi_0(A)) 
\cong H^2_s(G,A_0) \times \Hom(\pi_1(G), \pi_0(A)). $$

For the universal covering $q_G \: \tilde G\to G$ we thus obtain an
isomorphism 
$$ H^2_s(\tilde G, A_0) \to H^2_s(\tilde G, A) $$
because $\pi_1(\tilde G)$ is trivial. 
\qed

\sectionheadline{Appendix F. Multiplication in Lie algebra and Lie group cohomology} 

In this appendix we collect some information concerning multiplication of Lie algebra 
and Lie group cocycles which is used in Section~IX. 

\subheadline{Multiplication of Lie algebra cochains} 

Let $U,V,W$ be topological modules of the topological Lie algebra $\g$ and 
$m \: U \times  V\to W, (u,v) \mapsto u \cdot v$ 
a $\g$-equivariant continuous bilinear map. Then we define a product 
$$ C^p_c(\g,U) \times C^q_c(\g,V) \to C^{p+q}_c(\g,W), \quad 
(\alpha, \beta) \mapsto \alpha \wedge \beta $$
by 
$$ (\alpha \wedge \beta)(x_1,\ldots, x_{p+q}) 
:= {1\over p!q!} \sum_{\sigma \in S_{p+q}} 
\sgn(\sigma) 
\alpha(x_{\sigma(1)}, \ldots, x_{\sigma(p)}) 
\beta(x_{\sigma(p+1)}, \ldots, x_{\sigma(p+q)}). $$
For $p = q = 1$ we have in particular 
$$ (\alpha \wedge \beta)(x,y) = \alpha(x) \cdot \beta(y) - \alpha(y) \cdot \beta(x). $$
In the following we write for a $p$-linear map $\alpha \: \g^p \to V$: 
$$ \Alt(\alpha)(x_1, \ldots, x_p) 
:= \sum_{\sigma \in S_p} \sgn(\sigma) \alpha(x_{\sigma(1)}, \ldots, x_{\sigma(p)}). $$
In this sense we have 
$$ \alpha \wedge \beta = {1 \over p! q!} \Alt(\alpha \cdot \beta), $$
where 
$$ (\alpha \cdot \beta)(x_1, \ldots, x_{p+q}) 
:= \alpha(x_1, \ldots, x_p) \cdot \beta(x_{p+1}, \ldots, x_{p+q}).$$ 

\Lemma F.1. For $\alpha \in C^p_c(\g,U)$ and $\beta \in C^q_c(\g,V)$ we have 
$$ d_\g(\alpha \wedge \beta) = d_\g\alpha \wedge \beta 
+ (-1)^p \alpha \wedge d_\g\beta. \leqno({\rm F.1}) $$

\Proof. First we verify that for $x \in \g$ the insertion map $i_x$ satisfies 
$$ i_x(\alpha \wedge \beta) = 
i_x\alpha \wedge \beta 
+ (-1)^p \alpha \wedge i_x\beta. \leqno({\rm F.2}) $$
For $p = 0$ or $q = 0$ this formula is a trivial consequence 
of the definitions. We may therefore assume $p,q \geq 1$. 
We calculate for $x_1,\ldots, x_{p+q} \in \g$: 
$$ \eqalign{
&\ \ \ \ i_{x_1}(\alpha \wedge \beta)(x_2, \ldots, x_{p+q}) 
= (\alpha \wedge \beta)(x_1, x_2, \ldots, x_{p+q}) \cr
&= {1\over p!q!} \sum_{\sigma \in S_{p+q}} \sgn(\sigma) 
\alpha(x_{\sigma^{-1}(1)}, \ldots, x_{\sigma^{-1}(p)}) 
\beta(x_{\sigma^{-1}(p+1)}, \ldots, x_{\sigma^{-1}(p+q)}) \cr
&= {1\over p!q!} \sum_{\sigma(1) \leq p} \ldots  
+ {1\over p!q!} \sum_{\sigma(1) > p} \ldots. \cr} $$ 
For $\sigma(1) \leq p$ we get 
$$ \eqalign{ 
\alpha(x_{\sigma^{-1}(1)}, \ldots, x_{\sigma^{-1}(p)}) 
&= (-1)^{\sigma(1) + 1} 
\alpha(x_1, x_{\sigma^{-1}(1)}, \ldots, \hat{x_1}, \ldots, x_{\sigma^{-1}(p)}) \cr
&= (-1)^{\sigma(1) + 1} 
(i_{x_1}\alpha)(x_{\sigma^{-1}(1)}, \ldots, \hat{x_1}, \ldots, x_{\sigma^{-1}(p)}), \cr}$$
which leads to 
$$ \eqalign{ 
&\ \ \ \ {1\over p!q!} \sum_{\sigma(1) \leq p} \ldots \cr 
&= {1\over p!q!} \sum_{i = 1}^p \sum_{\sigma(1) = i} \sgn(\sigma) 
(-1)^{i+1} (i_{x_1}\alpha)(x_{\sigma^{-1}(1)}, \ldots, \hat{x_1}, \ldots, x_{\sigma^{-1}(p)}) 
\beta(x_{\sigma^{-1}(p+1)}, \ldots, x_{\sigma^{-1}(p+q)}) \cr
&= {1\over p!q!} \sum_{i = 1}^p 
\Alt(i_{x_1}\alpha \cdot \beta)(x_2, \ldots, x_{p+q}) 
= {1\over (p-1)!q!} 
\Alt(i_{x_1}\alpha \cdot \beta)(x_2, \ldots, x_{p+q}) \cr
&= (i_{x_1}\alpha \wedge  \beta)(x_2, \ldots, x_{p+q}). \cr} $$
We likewise obtain 
$$ {1\over p!q!} \sum_{\sigma(1) > p} \ldots 
= (-1)^p (\alpha \wedge  (i_{x_1}\beta))(x_2, \ldots, x_{p+q}). $$
This proves (F.2). 
 
We now prove (F.1) by induction on $p$ and $q$. For $p = 0$ we have 
$$ (\alpha \wedge \beta)(x_1, \ldots, x_q) = \alpha \cdot \beta(x_1, \ldots, x_q)$$ 
and 
$$ \eqalign{ 
d_\g(\alpha\wedge \beta)(x_0, \ldots, x_q) 
&= \sum_{i =0}^q (-1)^i x_i.(\alpha\cdot \beta)(x_0, \ldots, \hat{x_i}, \ldots, x_q) \cr
&\ \ \ \ + \sum_{i < j} (-1)^{i+j} \alpha \cdot \beta([x_i, x_j], \ldots, \hat{x_i}, \ldots, \hat{x_j}, \ldots, x_q) \cr
&= \sum_{i =0}^q (-1)^i (x_i.\alpha) \cdot \beta(x_0, \ldots, \hat{x_i}, \ldots, x_q) 
+ \alpha \cdot (d\beta)(x_0, \ldots, x_q) \cr} $$
and 
$$ \eqalign{ (d_\g\alpha \wedge \beta)(x_0, \ldots, x_q) 
&= {1\over q!} \sum_{\sigma \in S_{q+1}} \sgn(\sigma) 
(d_\g \alpha)(x_{\sigma(0)}) \cdot 
\beta(x_{\sigma(1)}, \ldots, x_{\sigma(q)}) \cr
&= {1\over q!}\sum_{i=0}^q  \sum_{\sigma(0)= i} \sgn(\sigma) 
(x_i.\alpha) \cdot \beta(x_{\sigma(1)}, \ldots, x_{\sigma(q)}) \cr
&= {1\over q!}\sum_{i=0}^q (-1)^{i} (x_i.\alpha)\cdot 
\Alt(\beta)(x_0, \ldots, \hat{x_i}, \ldots, x_q) \cr
&= \sum_{i=0}^q (-1)^i (x_i.\alpha)\cdot \beta(x_0, \ldots, \hat{x_i}, \ldots, 
x_q). \cr} $$
This proves (F.1) for $p =0$. 
A similar argument works for $q = 0$. 
We now assume that $p,q \geq 1$ and that (F.1) hold for the pairs $(p-1,q)$ 
and $(p,q-1)$. Then we obtain with the Cartan formulas and (F.2) for 
$x \in \g$: 
$$ \eqalign{ 
&\ \ \ \ i_x (d_\g\alpha \wedge \beta + (-1)^p \alpha \wedge d_\g\beta)\cr
&= (i_x d_\g\alpha) \wedge \beta + (-1)^{p+1} d_\g \alpha \wedge i_x \beta 
+ (-1)^p i_x \alpha \wedge d_\g\beta 
+ \alpha \wedge i_x d_\g \beta \cr 
&= x.\alpha \wedge \beta - d_\g(i_x\alpha) \wedge \beta 
+ (-1)^{p+1} d_\g \alpha \wedge i_x \beta 
+ (-1)^p i_x \alpha \wedge d_\g\beta 
+ \alpha \wedge x.\beta - \alpha \wedge d_\g(i_x \beta) \cr 
&= x.(\alpha \wedge \beta) - d_\g(i_x \alpha \wedge \beta) 
+ (-1)^{p+1} d_\g(\alpha \wedge i_x \beta) \cr
&= x.(\alpha \wedge \beta) - d_\g(i_x(\alpha \wedge \beta)) 
= i_x(d_\g(\alpha \wedge \beta)). \cr} $$
Since $x$ was arbitrary, this proves (F.1). 
\qed

The preceding lemma 
implies that products of two cocycles are cocycles 
and that the product of a cocycle with a coboundary is a coboundary, so that we obtain 
bilinear maps 
$$ H^p_c(\g,U) \times H^q_c(\g,V) \to H^{p+q}_c(\g,W), \quad 
([\alpha], [\beta]) \mapsto [\alpha \wedge \beta] $$
which can be combined to a product 
$$ H^*_c(\g,U) \times H^*_c(\g,V) \to H^*_c(\g,W). $$

\subheadline{Multiplication of group cochains} 

Now let $U,V,W$ be smooth modules of the Lie group $G$ and 
$m \: U \times  V\to W, (u,v) \mapsto u \cdot v$ 
a $G$-equivariant biadditive continuous map. Then we define a product 
$$ C^p_s(G,U) \times C^q_s(G,V) \to C^{p+q}_s(G,W), \quad 
(\alpha, \beta) \mapsto \alpha \cup \beta, $$
where 
$$ (\alpha \cup \beta)(g_1,\ldots, g_{p+q}) 
:= \alpha(g_1, \ldots, g_p) \cdot (g_1\cdots g_p).\beta(g_{p+1}, \ldots, g_{p+q}) $$
(cf.\ [Bro82, p.110] up to the different signs which are caused by different signs 
for the group differential). 

\Lemma F.2. For $\alpha \in C^p_s(G,U)$ and $\beta \in C^q_s(G,V)$ we have 
$$ d_G(\alpha \cup \beta)  
= d_G \alpha \cup \beta  + (-1)^p \alpha \cup d_G \beta. $$

\Proof. For $g_0, \ldots, g_{p+q}\in G$ we have 
$$ \eqalign{ 
&\ \ \ \ d_G(\alpha \cup \beta)(g_0, \ldots, g_{p+q}) \cr 
&= g_0.(\alpha \cup \beta)(g_1, \ldots, g_{p+q}) 
+ \sum_{i=1}^{p+q} (-1)^i (\alpha \cup \beta)(g_0,\ldots, g_{i-1} g_i, \ldots, g_{p+q})\cr
&\ \ \ \ + (-1)^{p+q+1} (\alpha \cup \beta)(g_0,\ldots, g_{p+q-1}) \cr 
&= (g_0.\alpha(g_1, \ldots, g_p)) \cdot (g_0 \cdots g_p).\beta(g_{p+1}, \ldots, g_{p+q}) \cr
&+ \sum_{i=1}^p (-1)^i \alpha(g_0,\ldots, g_{i-1} g_i, \ldots, g_{p})
\cdot g_0 \cdots g_p.\beta(g_{p+1}, \ldots, g_{p+q}) \cr
&+ \sum_{i=p+1}^{p+q} (-1)^i \alpha(g_0,\ldots, g_{p-1}) 
\cdot g_0 \cdots g_{p-1}.\beta(g_p, \ldots, g_{i-1} g_i, \ldots, g_{p+q}) \cr
&\ \ \ \ + (-1)^{p+q+1} \alpha(g_0, \ldots, g_{p-1}) \cdot 
(g_0 \cdots g_{p-1}).\beta(g_p, \ldots, g_{p+q-1}) \cr 
&= (d_G \alpha)(g_0, \ldots, g_p) \cdot (g_0 \cdots g_p).\beta(g_{p+1}, \ldots, g_{p+q}) \cr
& + (-1)^p \alpha(g_0, \ldots, g_{p-1}) \cdot 
(g_0 \cdots g_{p}).\beta(g_{p+1}, \ldots, g_{p+q}) \cr 
&+ \alpha(g_0,\ldots, g_{p-1}) \cdot \cr
&\ \ \ \ g_0 \cdots g_{p-1}.\Big(\sum_{i=p+1}^{p+q} (-1)^i 
\beta(g_p, \ldots, g_{i-1} g_i, \ldots, g_{p+q}) 
+ (-1)^{p+q+1} \beta(g_p, \ldots, g_{p+q-1})\Big) \cr 
&= (d_G \alpha \cup \beta)(g_0, \ldots, g_{p+q}) 
 + (-1)^p (\alpha \cup d_G \beta)(g_0, \ldots, g_{p+q}). \cr} $$
\qed

Lemma F.2 implies that products of two cocycles are cocycles 
and that the product of a cocycle with a coboundary is a coboundary, so that we obtain 
biadditive maps 
$$ H^p_s(G,U) \times H^q_s(G,V) \to H^{p+q}_s(G,W), \quad 
([\alpha], [\beta]) \mapsto [\alpha \cup \beta]. $$

The following lemma shows that for Lie groups the multiplication of group and 
Lie algebra cochains is compatible with the differentiation map $D$. 

\Lemma F.3. If $G$ is a Lie group, $U$, $V$ and $W$ are smooth modules and 
$\beta \: U \times V \to W$ is continuous bilinear and equivariant, 
then we have for $\alpha \in C^p(G,U)$ and 
$\beta \in C^q(G,V)$ we have 
$$ D(\alpha \cup \beta) = D \alpha \wedge D\beta $$
in $C^{p+q}_c(\g,W)$. 

\Proof. In view of 
$D\alpha = \Alt(d^p \alpha(\1,\ldots, \1)),$
we get 
$$ \eqalign{ 
D\alpha \wedge D_\beta 
&= {1\over p!q!} \Alt(D\alpha \cdot D\beta)
= {1\over p!q!} \Alt(\Alt(d^p\alpha(\1,\ldots, \1)) \cdot \Alt(d^q \beta(\1,\ldots, \1)))\cr
&= \Alt(d^p\alpha(\1,\ldots, \1) \cdot d^q\beta(\1, \ldots, \1)), \cr}$$
so that it remains to see that 
$$ d^{p+q}(\alpha \cup \beta)(\1,\ldots, \1) 
= (d^p\alpha)(\1,\ldots, \1) \cdot (d^q\beta)(\1,\ldots, \1), $$
but this follows immediately from the normalization of the cocycles and 
the chain rule for jets, applied to the multiplication map $\beta$. 
\qed

\def\entries{

\[AK98 Arnold, V. I., and B. A. Khesin, ``Topological Methods in Hydrodynamics,'' 
Springer-Verlag, 1998 

\[AI95 de Azcarraga, J. A., and J. M. Izquierdo, ``Lie Groups, Lie
Algebras, Cohomology and some Applications in Physics,'' Cambridge
Monographs on Math. Physics, 1995 

\[Bi03 Billig, Y., {\it Abelian extensions of the group of
  diffeomorphisms of a torus}, Lett. Math. Phys.  {\bf 64:2}  (2003),  155--169 

\[Bou88 Bourbaki, N., ``General Topology,'' Chaps. I--IV, Springer
Verlag, 1988 

\[Bre93 Bredon, G.\ E., ``Topology and Geometry,'' Graduate Texts in
Mathematics {\bf 139}, Springer-Verlag, Berlin, 1993 

\[Bro82 Brown, K. S., ``Cohomology of Groups,'' Graduate Texts in Mathematics 
{\bf 87}, Springer-Verlag, 1982 

\[Bry90 Brylinski, J.-L., ``Loop Spaces, Characteristic Classes and Geometric 
Quantization,'' Progress in Mathematics, Birkh\"auser Verlag, Basel, 1990 

\[Ca52 Cartan, E., {\it La topologie des espaces repr\'esentifs de groupes de
Lie}, Oeuvres I, Gauthier--Villars, Paris, {\bf 2} (1952), 1307--1330 

\[CVLL98 Cassinelli, G., E. de Vito, P. Lahti and A. Levrero, {\it Symmetries of the quantum 
state space and group representations}, Reviews in Math. Physics {\it 10:7} 
(1998), 893--924 

\[Ch46 Chevalley, C., ``Theory of Lie Groups I,'' Princeton Univ.\ Press, 1946
 
\[CE48 Chevalley, C. and S. Eilenberg, {\it Cohomology theory of Lie groups and Lie 
algebras}, Transactions of the Amer. Math. Soc. {\bf 63} (1948), 85--124 
 
\[Dz92 Dzhumadildaev, A., {\it Central extensions of
infinite-dimensional Lie algebras}, Funct. Anal. Appl. {\bf 26:4}
(1992), 247--253 

\[vE55 van Est, W. T., {\it On the algebraic cohomology concepts in
Lie groups I,II}, Proc. of the Koninglijke Nederlandse Akademie van 
Wetenschappen {\bf A58} (1955), 225--233; 286--294  

\[EK64 van Est, W.\ T., and Th.\ J.\ Korthagen, {\it Non enlargible Lie
algebras}, Proc.\ Kon.\ Ned. Acad. v. Wet. A {\bf 67} (1964), 15--31 

\[Fu86 Fuks, D.B., ``Cohomology of Infinite Dimensional Lie Algebras,''
Contemp. Sov. Math., Consultants Bureau, New York, London, 1986 

\[GF68 Gelfand, I., and D. B. Fuchs, {\it Cohomology of the Lie algebra of vector fields 
on the circle}, Funct. Anal. Appl. {\bf 2:4} (1968), 342--343 

\[God71  Godbillon, C., ``El\'ements de Topologie Alg\'ebrique,''
Hermann, Paris, 1971 

\[Gl01 Gl\"ockner, H., {\it Infinite-dimensional Lie groups without completeness 
condition}, in ``Geometry and Analysis on Finite-
and Infinite-Dimensional Lie Groups,'' A.~Strasburger et al Eds., 
Banach Center Publications {\bf 55}, Warszawa 2002; 53--59 

\[HV04 Haller, S., and C. Vizman, {\it Nonlinear Grassmannians as coadjoint orbits}, 
Math. Annalen, to appear 

\[Hel78 Helgason, S., ``Differential Geometry, Lie Groups, and
Symmetric Spa\-ces", Acad. Press, London, 1978

\[HS53a Hochschild, G., and J.-P. Serre, {\it Cohomology of group
extensions}, Transactions of the Amer. Math. Soc. {\bf 74} (1953),
110--134 

\[HS53b ---, {\it Cohomology of Lie
algebras}, Annals of Math. {\bf 57:2} (1953), 591--603 

\[HoMo98 Hofmann, K.\ H., and S.\ A.\ Morris, ``The Structure of
Compact Groups,'' Studies in Math., de Gruyter, Berlin, 1998

\[Is96 Ismagilov, R. S., ``Representations of Infinite-Dimensional
Groups,'' Translations of Math. Monographs {\bf 152}, Amer. Math. Soc., 1996 

\[KM97 Kriegl, A., and P.\ Michor, ``The Convenient Setting of
Global Analysis,'' Math.\ Surveys and Monographs {\bf 53}, Amer.\
Math.\ Soc., 1997 

\[La99 Larsson, T. A., {\it Lowest-energy representations of non-centrally extended 
  diffeomorphism algebras}, Commun. Math. Phys. {\bf 201} (1999), 461--470 

\[Li74 Lichnerowicz, A. {\it Alg\`ebre de Lie des automorphismes infinit\'esimaux 
d'une structure unimodulaire}, Ann. Inst. Fourier {\bf 24} (1974), 219--266 

\[ML78 Mac Lane, S., {\it Origins of the cohomology of groups}, L'Enseignement math\'em. 
{\bf 24:2} (1978), 1--29 

\[MN03 Maier, P., and K. - H. Neeb,  {\it Central extensions of current groups\/}, 
Math. Annalen {\bf 326:2} (2003), 367--415 

\[MDS98 McDuff, D., and D. Salamon, ``Introduction to Symplectic
Topology,'' Oxford Math. Monographs, 1998 

\[Mi89 Mickelsson, J., ``Current Algebras and Groups,'' Plenum Press,
New York, 1989 

\[Mil83 Milnor, J., {\it Remarks on infinite-dimensional Lie groups},
Proc. Summer School on Quantum Gravity, B. DeWitt ed., Les Houches, 1983


\[Mo64 Moore, C. C., {\it Extensions 
and low dimensional cohomology theory of locally compact 
groups. I,II}, Trans. Amer. Math. Soc. {\bf 113} (1964), 40--63; 63---86  

\[Mo76 ---, {\it Group extensions and cohomology of locally compact groups. III, IV}, 
Trans. Amer. Math. Soc. {\bf 221:1} (1976), 1--33; 35--58 

\[Ne02 Neeb, K.-H., {\it Central extensions of infinite-dimensional
Lie groups}, Annales de l'Inst. Fourier 52:5 (2002), 1365--1442 

\[Ne02b ---, {\it Nancy Lectures on Infinite-Dimensional Lie Groups}, 
Preprint {\bf 2203}, TU Darmstadt, March 2002 

\[Ne03a ---, {\it Locally convex root graded Lie algebras}, 
Travaux math\'emathiques {\bf 14} (2003), 25--120.  

\[Ne03b ---, {\it Universal central extensions of Lie groups},
Acta Appl. Math. {\bf 73:1,2} (2002), 175--219 

\[Ne04 ---, {\it General extensions of infinite-dimensional Lie groups}, in
preparation 

\[NV03 Neeb, K.-H., and C. Vizman, {\it Flux homomorphisms and principal bundles over
infinite-dimensional manifolds}, Monatshefte f\"ur Math. {\bf 139} (2003), 309--333 

\[OR98 Ovsienko, V., and C. Roger, {\it Generalizations of Virasoro
  group and Virasoro algebra 
through extensions by modules of tensor densities on $S^1$}, 
Indag. Mathem. N. S. {\bf 9:2} (1998), 277--288 

\[Se02 Segal, G., Personal communications, March 2002. 

\[Sp66 Spanier, E. H., ``Algebraic Topology,'' Mc Graw Hill Series in Higher 
Math. 1966, Univ. of Cal., Berkeley 

\[Ste51 Steenrod, N., ``The Topology of Fibre Bundles,'' Princeton
University Press, Princeton, New Jersey, 1951 



\[Ti83 Tits, J., ``Liesche Gruppen und Algebren", Springer, New York, 
Heidelberg, 1983

\[Vi02 Vizman, C., {\it Geodesics on extensions of the Lie algebra of
vector fields on the circle}, in ``Geometry and Analysis on Finite-
and Infinite-Dimensional Lie Groups,'' A. Strasburger et al Eds., 
Banach Center Publications {\bf 55}, Warszawa 2002; 165--172

\[We95 Weibel, C. A., ``An introduction to homological algebra,''
Cambridge studies in advanced math. {\bf 38}, Cambridge Univ. Press,
1995 

}

\references
\lastpage 

\vfill\eject
\bye